\documentclass[12pt,a4paper,reqno]{amsart}
\usepackage[a4paper,top=3cm,bottom=3cm,inner=3cm,outer=3cm,includehead,includefoot]{geometry}
\usepackage{amsmath,amssymb,amsthm,calc,verbatim,enumitem,tikz,url,hyperref,mathrsfs,cite,graphicx,caption}
\usetikzlibrary{shapes.misc,calc,intersections,patterns,decorations.pathreplacing}
\hypersetup{colorlinks=true,citecolor=blue,linktoc=page}
\newtheorem{theorem}{Theorem}
\newtheorem{lemma}[theorem]{Lemma}

\newtheorem{obs}[theorem]{Observation}
\newtheorem{claim}[theorem]{Claim}

\newtheorem{conjecture}[theorem]{Conjecture}
\theoremstyle{definition}
\newtheorem{definition}[theorem]{Definition}
\newtheorem{remark}[theorem]{Remark}

\numberwithin{theorem}{section}
\setlist[itemize]{leftmargin=1cm}
\setlist[enumerate]{leftmargin=1.5cm}

\everydisplay{\thickmuskip=7mu}
\newcommand{\rob}[1]{{\color{black}{#1}}}   

\renewcommand{\leq}{\leqslant}
\renewcommand{\geq}{\geqslant}
\renewcommand{\le}{\leqslant}
\renewcommand{\ge}{\geqslant}
\renewcommand{\to}{\rightarrow}
\newcommand{\diam}{\operatorname{diam}}
\newcommand{\sh}{\operatorname{sh}}
\newcommand{\interior}{\operatorname{int}}

\let\epsilon\varepsilon
\let\eps\varepsilon

\def\H{\mathbb{H}}
\def\N{\mathbb{N}}
\def\Pr{\mathbb{P}}
\def\Q{\mathbb{Q}}
\def\R{\mathbb{R}}
\def\Z{\mathbb{Z}}
\def\C{\mathcal{C}}
\def\D{\mathcal{D}}
\def\J{\mathcal{J}}
\def\K{\mathcal{K}}
\def\S{\mathcal{S}}
\def\T{\mathcal{T}}
\def\U{\mathcal{U}}
\def\W{\mathcal{W}}
\def\Y{\mathcal{Y}}
\def\cZ{\mathcal{Z}}
\def\stab{\mathcal{S}}
\def\quasi{\mathcal{S}_Q}
\def\qquasi{\mathcal{Q}}
\def\stabb{\mathcal{S}_B}
\def\stabu{\mathcal{S}_U}
\def\hier{\mathcal{H}}
\def\root{v_{\text{root}}}
\def\weak{\kappa_0}
\def\ispan{I^\times}
\def\ifill{I}
\def\<{\langle}
\def\>{\rangle}
\def\alphabar{\bar{\alpha}}
\def\0{\mathbf{0}}
\def\ih{\mathrm{IH}}
\def\Dup{D^\uparrow}
\def\Ddown{D^\downarrow}

\title[Universality for two-dimensional critical cellular automata]{Universality for two-dimensional critical \\ cellular automata}

\author[B. Bollob\'as \and H. Duminil-Copin \and R. Morris \and P.J. Smith]{B\'ela Bollob\'as \and Hugo Duminil-Copin \and Robert Morris \and Paul Smith}

\address{Department of Pure Mathematics and Mathematical Statistics, Wilberforce Road, Cambridge, CB3 0WA, UK, and Department of Mathematical Sciences, University of Memphis, Memphis, TN 38152, USA}
\email{b.bollobas@dpmms.cam.ac.uk}
\address{D\'epartement de Math\'ematiques, Universit\'e de Gen\`eve, 2-4 Rue du Li\`evre, Case Postale 64, 1211 Gen\`eve 4, Switzerland}
\email{hugo.duminil@unige.ch}
\address{IMPA, Estrada Dona Castorina 110, Jardim Bot\^anico, Rio de Janeiro, 22460-320, Brazil}
\email{rob|psmith@impa.br}

\thanks{B.B.\ is partially supported by NSF grant DMS~1855745, R.M.\ by CNPq (Proc.~303275/2013-8), by FAPERJ (Proc.~201.598/2014), and by the ERC Starting Grant 680275 MALIG, and P.S.\ by a CNPq bolsa PDJ}

\date{\today}


\begin{document}

\begin{abstract}
We study the class of monotone, two-state, deterministic cellular automata, in which sites are activated (or `infected') by certain configurations of nearby infected sites. These models have close connections to statistical physics, and several specific examples have been extensively studied in recent years by both mathematicians and physicists. This general setting was first studied only recently, however, by Bollob\'as, Smith and Uzzell, who showed that the family of all such `bootstrap percolation' models on $\Z^2$ can be naturally partitioned into three classes, which they termed subcritical, critical and supercritical. 

In this paper we determine the order of the threshold for percolation (complete occupation) for every critical bootstrap percolation model in two dimensions. This `universality' theorem includes as special cases results of Aizenman and Lebowitz, Gravner and Griffeath, Mountford, and van Enter and Hulshof, significantly strengthens bounds of Bollob\'as, Smith and Uzzell, and complements recent work of Balister, Bollob\'as, Przykucki and Smith on subcritical models. 
\end{abstract}

\maketitle

\tableofcontents

\newpage

\section{Introduction}

An important and challenging problem in statistical physics, probability theory and combinatorics is to understand the typical global behaviour of so-called `lattice models', including cellular automata, percolation models, and spin models. Although these models are defined in terms of local interactions between the sites of the lattice, it is typically observed in simulations that, in fixed dimensions, the macroscopic behaviour of the models does not seem to depend on the precise nature of these local interactions. Indeed, since the breakthrough work of Kadanoff~\cite{Kad} and the development of the renormalization group framework by Wilson~\cite{Wilson}, this phenomenon of \emph{universality} has been widely predicted to occur throughout statistical physics (see, for example,~\cite{Fisher,LPS,Schramm}). Despite this, it has been proved rigorously in only a small handful of cases. One example of a model for which universality is partially understood is the Ising model, for which it was proved recently that the critical exponents exist and are equal on a large class of planar graphs~\cite{CS,GGM}. 

Cellular automata are interacting particle systems whose update rules are local and homogeneous. Since their introduction by von Neumann~\cite{vN} almost 50 years ago, many particular such systems have been investigated, but no general theory has been developed for their study, and for many simple examples surprisingly little is known. In this paper we develop such a general theory for monotone, two-dimensional cellular automata with random initial configurations, which may also be thought of as monotone versions of the Glauber dynamics of the Ising model with arbitrary local interactions. The study of this general class of models was initiated only recently by Bollob\'as, Smith and Uzzell~\cite{BSUgen}, although many special cases had been studied earlier, beginning with the work of Chalupa, Leath and Reich~\cite{CLR} in 1979. We refer to these models as \emph{bootstrap percolation}, but we emphasize that they are vastly more general than the family of models that usually bears this name.

The class of models we study is defined as follows. Fix $d \in \N$ and let $\U = \{ X_1,\ldots,X_m \}$ be an arbitrary finite collection of finite subsets of $\Z^d \setminus \{0\}$. We call $\U$ the \emph{update family} of the process, each $X \in \U$ an \emph{update rule}, and the process itself \emph{\,$\U$-bootstrap percolation}. Let the lattice $\Lambda$ be either $\Z^d$ or $\Z_n^d$ (the $d$-dimensional discrete torus). Now given a set $A \subset\Lambda$ of initially \emph{infected} sites, set $A_0 = A$, and define for each $t \geq 0$, 
\[
A_{t+1} = A_t \cup \big\{ x \in \Lambda \,:\, x + X \subset A_t \text{ for some } X \in \U \big\}.
\]
Thus, a site $x$ becomes infected at time $t+1$ if the translate by $x$ of one of the sets of the update family is already entirely infected at time $t$, and infected sites remain infected forever. The set $\bigcup_{t \geq 0} A_t$ of eventually infected sites is the \emph{closure} of $A$, denoted by $[A]$. We say that $A$ \emph{percolates} if $[A] = \Lambda$.

As mentioned above, this model was first introduced (in its full generality) only recently, in~\cite{BSUgen}, although various special cases were introduced and studied much earlier by several different authors; see for example \cite{CLR,DCvE,DCH,GG93,GG99}. Indeed, the general class of $\U$-bootstrap percolation models is easily seen to include as specific examples all previously studied bootstrap percolation models on lattice graphs. For example, the update family of the classical $r$-neighbour model on $\Z^d$, the most well-studied of all models \cite{AL,CerfCir,CerfManzo,Hol,BBM3D,BBDCM}, consists of the $\binom{2d}{r}$ $r$-subsets of the $2d$ nearest neighbours of the origin. The $r$-neighbour models are themselves examples of \emph{threshold models}, which, in their full generality, consist of the $r$-element subsets of an arbitrary finite set $Y \subset \Z^d\setminus\{0\}$. 

We say that an update family $\U$ is \emph{symmetric} if $X \in \U$ implies $-X \in \U$, so in particular a threshold model is symmetric if and only if the set $Y$ is centrally symmetric. With a single exception (the work of Mountford~\cite{MountDuarte} on the Duarte model, which exhibits a weaker notion of symmetry, see below), we are not aware of any results on non-symmetric update families before the work of~\cite{BSUgen}. In the absence of symmetry, many of the previously developed techniques appear to break down, and overcoming this obstacle is one of the main tasks of this paper. 

 
Motivated by applications to statistical physics (see, e.g.,~\cite{FSSIsing,HMT1,HMT2,MorrisGlauber,T22}), we shall study the global behaviour of the $\U$-bootstrap process acting on random initial sets. Specifically, let us say that a set $A \subset \Lambda$ is \emph{$p$-random} if each of the sites of $\Lambda$ is included in $A$ independently with probability $p$. The key question is that of how likely it is that a random set $A$ percolates on the lattice $\Lambda$; in particular, one would like to know how large $p$ must be before percolation becomes likely. The point at which this phase transition occurs is measured by the \emph{critical probability},
\[
p_c(\Lambda,\U) := \inf \Big\{ p \,:\, \Pr_p\big( A \text{ percolates in $\U$-bootstrap percolation on $\Lambda$}\big) \geq 1/2 \Big\},
\]
where $\Pr_p$ denotes the product probability measure on $\Lambda$ with density $p$.\footnote{Thus a $p$-random set is one chosen according to the distribution $\Pr_p$.}

For the $r$-neighbour model on $\Z_n^d$, with $d$ fixed and $n\to\infty$, a great deal is known about the critical probability. Up to a constant factor, the threshold was determined by Aizenman and Lebowitz~\cite{AL} for $r=2$, by Cerf and Cirillo~\cite{CerfCir} for $d=r=3$, and by Cerf and Manzo~\cite{CerfManzo} for all remaining $2\leq r\leq d$. The constant factor was later improved to a $1+o(1)$ factor by Holroyd~\cite{Hol} in the case $d=2$, by Balogh, Bollob\'as and Morris~\cite{BBM3D} for $d=3$, and by Balogh, Bollob\'as, Duminil-Copin and Morris~\cite{BBDCM} for all $d\geq 4$. The $r$-neighbour model has also attracted attention on lattices with the dimension $d$ tending to infinity (for example the hypercube)~\cite{BBMhigh,BBMmaj}, and on graphs other than lattices, including trees~\cite{BPP,BGHJP} and random graphs~\cite{BalPit,JLTV}.

For lattice models other than the $r$-neighbour model, considerably less is known. Among the exceptions are two-dimensional balanced, symmetric threshold models with star-neighbourhoods\footnote{These terms are defined below.}, for which the critical probability was determined up to a constant factor by Gravner and Griffeath~\cite{GG99}, and asymptotically by Duminil-Copin and Holroyd~\cite{DCH}. Some results about the critical probabilities of a rather limited number of so-called unbalanced models are also known; these were proved by Mountford~\cite{MountDuarte}, the authors of the present paper~\cite{BDMSDuarte}, van Enter and Hulshof~\cite{vEH}, Duminil-Copin and van Enter~\cite{DCvE}, all in two dimensions, and by van Enter and Fey~\cite{vEF} and Blanquicett~\cite{B19} in three dimensions.

For the remainder of the paper, with the exception of a brief discussion of higher dimensions in Section~\ref{se:highd}, we restrict our attention to the case $d=2$. As we shall see shortly, one of the key properties of the $\U$-bootstrap process is that its rough global behaviour depends only on the action of the process on discrete half-planes. In order to make this statement precise, let us introduce a little notation. For each $u \in S^1$, let $\H_u := \{x \in \Z^2 : \< x,u \> < 0 \}$ be the discrete half-plane whose boundary is perpendicular to $u$. We say that $u$ is a \emph{stable direction} if $[\H_u] = \H_u$ and we denote by $\stab = \stab(\U) \subset S^1$ the collection of stable directions. Note that a direction $u \in S^1$ is unstable if and only if there exists a rule $X \in \U$ such that $X \subset \H_u$.  

The following classification of two-dimensional update families was proposed by Bollob\'as, Smith and Uzzell~\cite{BSUgen}.

\begin{definition}\label{de:class}
An update family $\U$ is:
\begin{itemize}
\item \emph{subcritical} if every semicircle in $S^1$ has infinite intersection with $\stab$; \vspace{0.1cm}
\item \emph{critical} if there exists a semicircle in $S^1$ that has finite intersection with $\stab$, and if every open semicircle in $S^1$ has non-empty intersection with $\stab$; \vspace{0.1cm}
\item \emph{supercritical} if there exists an open semicircle in $S^1$ that is disjoint from $\stab$. 
\end{itemize}
\end{definition}

The justification of the above definition was completed in two stages. First, in their original paper, Bollob\'as, Smith and Uzzell~\cite{BSUgen} proved that the critical probabilities of supercritical families are polynomial, while those of critical families are polylogarithmic. Later, Balister, Bollob\'as, Przykucki and Smith~\cite{BBPS} proved that the critical probabilities of subcritical models are bounded away from zero. The combination of the results of~\cite{BSUgen} and~\cite{BBPS} may be summarized as follows\footnote{Our asymptotic notation is mostly standard; however, for the avoidance of ambiguity, precise definitions are given in Section~\ref{se:defs}.}:\smallskip
\begin{itemize}
\item if $\U$ is subcritical then $\displaystyle \liminf_{n\to\infty} \, p_c(\Z_n^2,\U)>0$; \vspace{0.15cm} 
\item if $\U$ is critical then $p_c(\Z_n^2,\U) = (\log n)^{-\Theta(1)}$; \vspace{0.35cm}
\item if $\U$ is supercritical then $p_c(\Z_n^2,\U) = n^{-\Theta(1)}$.\smallskip
\end{itemize}

In this paper we significantly strengthen the bounds of~\cite{BSUgen} by determining the threshold $p_c(\Z_n^2,\U)$ up to a constant factor for every critical update family. This result, which may be thought of as a universality statement for two-dimensional critical bootstrap percolation, was previously known only in the case of one very restrictive subclass of critical models~\cite{AL,GG99}, namely the symmetric, balanced threshold models, and just two other specific models~\cite{MountDuarte,vEH}.

The form of the threshold function depends on two properties of $\U$: the `difficulty' of $\U$, and whether or not $\U$ is `balanced'. In order to explain what these terms mean, first we need a quantitative measure of how hard it is to grow in each direction.

\begin{figure}[ht]
 \begin{minipage}{.45\textwidth}
  \centering
  \includegraphics[width=0.9\textwidth]{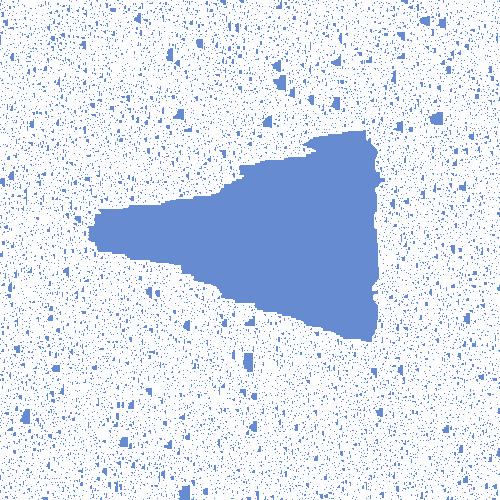} \\ \vspace{0.1cm} $(a)$
 \end{minipage}
 \begin{minipage}{.45\textwidth}
  \centering
  \includegraphics[width=0.9\textwidth]{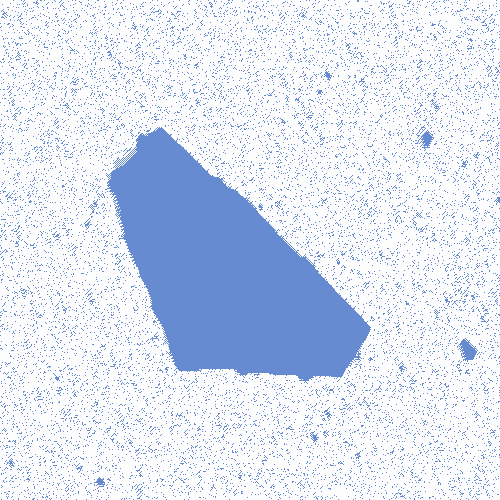} \\ \vspace{0.1cm} $(b)$
 \end{minipage}
 \captionsetup{singlelinecheck=false}
\caption[]{Two examples of critical models, paused during their evolution on $\Z^2$. In $(a)$, the Duarte model, an unbalanced model with drift; in $(b)$, a balanced critical model.}
\end{figure}

Let $\Q_1 \subset S^1$ denote the set of rational directions on the circle\footnote{That is, the set of all $u \in S^1$ such that the line $\{\lambda u : \lambda \in \R \}$ intersects the set $\Z^2 \setminus \{\0\}$.}, 
and for each $u \in \Q_1$, let $\ell_u^+$ be the subset of the line $\ell_u := \{x \in \Z^2 : \< x,u \> = 0 \}$ consisting of the origin and the sites to the right of the origin as one looks in the direction of $u$. Similarly, let $\ell_u^- := (\ell_u\setminus\ell^+_u) \cup \{0\}$ consist of the origin and the sites to the left of the origin. Note that the line $\ell_u$ is infinite for every $u \in \Q_1$. 

\begin{definition}\label{de:alpha}
Given $u \in \Q_1$, the \emph{difficulty} of $u$ is
\begin{equation}\label{eq:alpha}
\alpha(u) :=
\begin{cases}
\min\big\{ \alpha^+(u), \alpha^-(u) \big\} &\text{if } \alpha^+(u) < \infty \text{ and } \alpha^-(u) < \infty \\
\hfill \infty \hfill  & \text{otherwise,}
\end{cases}
\end{equation}
where $\alpha^+(u)$ (respectively $\alpha^-(u)$) is defined to be the minimum (possibly infinite) cardinality of a set $Z \subset \Z^2$ such that $[\H_u \cup Z]$ contains infinitely many sites of $\ell_u^+$ (respectively $\ell_u^-$). If $u \in S^1 \setminus \Q_1$, then $\alpha(u) := \infty$ if $u \in \stab$, and $\alpha(u) := 0$ otherwise.

We define the \emph{difficulty} of $\U$ to be
\begin{equation}\label{eq:alphaU}
\alpha = \alpha(\U) := \min_{C \in \C} \, \max_{u \in C} \, \alpha(u),
\end{equation}
where $\C$ denotes the collection of open semicircles of $S^1$.
\end{definition}

We would like to emphasize that $\alpha(u) < \infty$ only if \emph{both} $\alpha^+(u)$ and $\alpha^-(u)$ are finite (in which case it is their minimum). It moreover follows from simple properties of stable sets (see Section~\ref{se:simpleprops}) that $\alpha(u) > 0$ if and only if $u$ is a stable direction.

In Section~\ref{se:sketch} we discuss why these definitions of the difficulty of a direction under the action of $\U$ and of the difficulty of $\U$ itself are the natural ones. The final definition we need is as follows.

\begin{definition}\label{de:balanced}
A critical update family $\U$ is \emph{balanced} if there exists a closed semicircle $C$ such that $\alpha(u) \leq \alpha$ for all $u\in C$. It is said to be \emph{unbalanced} otherwise.
\end{definition}

The distinction between the open semicircles in the definition of $\alpha$ and the closed semicircles in the definition of balanced is subtle but important. It turns out that growth under the action of balanced critical families is completely two-dimensional, while growth under the action of unbalanced critical families is asymptotically one-dimensional. Despite this, the analysis of unbalanced families presents by far the greater number of difficulties.

The following theorem is the main result of this paper.

\begin{theorem}\label{th:main}
Let $\U$ be a critical two-dimensional bootstrap percolation update family.
\begin{enumerate}
\item If\/ $\U$ is balanced, then
\[
p_c\big( \Z_n^2,\U \big) = \Theta \bigg( \frac{1}{\log n} \bigg)^{1/\alpha}.
\]
\item If\/ $\U$ is unbalanced, then
\[
p_c\big( \Z_n^2,\U \big) = \Theta \bigg( \frac{(\log \log n)^2}{\log n} \bigg)^{1/\alpha}.
\]
\end{enumerate}
\end{theorem}

On the infinite lattice (where the critical probability is zero), one can state an essentially equivalent version of Theorem~\ref{th:main} in terms of the infection time of the origin. To be precise, given $A \subset \Z^2$, define
\[
\tau = \tau(A,\U) := \min\big\{ t \geq 0 \,:\, \0 \in A_t \big\}
\]
to be the time at which the origin is infected in the $\U$-bootstrap process on $\Z^2$ with $A_0 = A$. We write `with high probability' to mean `with probability tending to 1'.\footnote{We also say that ``$Z = \Theta\big( f(p) \big)$ with high probability", where $Z$ is a random variable, if there exist constants $c,C > 0$ (not depending on $p$) such that $cf(p) \leq Z \leq Cf(p)$ with high probability.}


\begin{theorem}\label{th:time}
Let $\U$ be a critical two-dimensional bootstrap percolation update family, and let $A$ be a $p$-random subset of $\Z^2$.
\begin{enumerate}
\item If $\U$ is balanced, then, with high probability as $p \to 0$,
\[
p^\alpha \log\tau = \Theta(1).
\]
\item If $\U$ is unbalanced, then, with high probability as $p \to 0$,
\[
p^\alpha \left(\log\frac{1}{p}\right)^{-2} \log\tau = \Theta(1).
\]
\end{enumerate}
\end{theorem}

The proof of Theorem~\ref{th:time} follows from the same method used to prove Theorem~\ref{th:main}, together with a (fairly standard) additional argument (see,~e.g.,~\cite{DCvE}), and therefore (since it is not our main focus) we will omit the details.

We noted earlier that various special cases of Theorem~\ref{th:main} have already been proved in the literature. The critical probability of the 2-neighbour model was established by Aizenman and Lebowitz~\cite{AL} using methods that have become central to the study of bootstrap percolation, including the `rectangles process' and the notion of a `critical droplet' (see Section~\ref{se:sketch} for details). Mountford~\cite{MountDuarte} determined the critical probability of the Duarte model, which is the unbalanced threshold model consisting of all two-element subsets of
\[
\big\{ (-1,0), (0,1), (0,-1) \big\}.
\]
His proof was based on martingale techniques, which makes it unique among proofs of this type of theorem. Gravner and Griffeath~\cite{GG99} generalized the result of Aizenman and Lebowitz to a class of balanced, symmetric threshold models, using somewhat non-rigorous methods. The critical probability of one further unbalanced model, namely the one consisting of all three-element subsets of
\[
\big\{ (-2,0), (-1,0), (0,1), (0,-1), (1,0), (2,0) \big\},
\]
was determined by van Enter and Hulshof~\cite{vEH}, correcting an assertion of Gravner and Griffeath~\cite{GG99}. Until now, the models studied by Mountford~\cite{MountDuarte} and by van Enter and Hulshof~\cite{vEH} were the only two unbalanced models whose critical probabilities were known, and they were, respectively, the unique such examples of `drift' and `non-drift' unbalanced models.\footnote{These terms are explained in Section~\ref{se:sketch}, but roughly speaking, the term `drift' refers to the phenomenon that occurs when $u\in\stab$ is such that $\alpha(u)=\infty$ but $\min\big\{\alpha^-(u),\alpha^+(u)\big\}<\infty$, which in certain cases causes the growth to be biased in one direction.}

One property that all of these previously studied models share, and one that simplifies the problem enormously, is that of symmetry. In all but the Duarte model, the symmetry is particularly strong, in that $X\in\U$ if and only if $-X\in\U$. The symmetry of the Duarte model is weaker (the useful property is that there exists a parallelogram of stable directions $\{u,-u,v,-v\}\subset\stab$), but it is enough to make a significant difference to the proof. An important aspect of the general models that we study -- perhaps \emph{the} most important aspect -- is the lack of any symmetry assumptions. Indeed, it is little exaggeration to say that the main task of this paper (as was that of~\cite{BSUgen}) is to handle the lack of symmetry, which causes all previously known techniques to break down.

In all of the above cases (namely, the 2-neighbour model of Aizenman and Lebowitz, the symmetric, balanced threshold models of Gravner and Griffeath, the Duarte model of Mountford, and the unbalanced model of van Enter and Hulshof) -- but in no other fundamentally different cases -- the critical probability has now been determined up to a $1+o(1)$ factor. These results are due to Holroyd~\cite{Hol}, Duminil-Copin and Holroyd~\cite{DCH}, the authors of the present paper~\cite{BDMSDuarte} and Duminil-Copin and van Enter~\cite{DCvE}, respectively, and in some cases, even sharper results are known~\cite{DEH,GravHolMor,HM}. Obtaining similarly sharp bounds for the general model is likely to be an important, but extremely difficult, direction for future research.

The organization of the rest of the paper is as follows. In Section~\ref{se:sketch} we give an outline of the proof, we introduce some notation, and we recall a number of basic facts about $\U$-bootstrap percolation from~\cite{BSUgen}. In Section~\ref{se:upperprelims} we lay the groundwork for the proofs of the upper bounds of Theorem~\ref{th:main}, which are then proved in Sections~\ref{se:upperbal} (balanced case) and~\ref{se:upperunbal} (unbalanced case). In Section~\ref{se:approxfill} we define three different notions of `approximately internally filled' sets and prove a number of deterministic properties of such sets. In Section~\ref{se:lowerbal} we deduce the lower bound in the balanced case. The hardest part of the proof is the lower bound in the unbalanced case, which is contained in Section~\ref{se:lowerunbal}. Finally, we end the paper with some open problems, including a discussion of the problem in higher dimensions.

\section{Outline of the proof}\label{se:sketch}

Let us begin by explaining why the definition~\eqref{eq:alpha} of $\alpha(u)$ is the right notion of the `difficulty' of growing in a direction $u\in S^1$. The key fact is that there is a sense (which is formalized in Lemma~\ref{le:line}) in which $\alpha(u)$ measures how hard it is to infect an entire new line in direction $u$, rather than merely an infinite subset of the line. More specifically, while the definition of $\alpha(u)$ only guarantees that there exist sets of $\alpha(u)$ sites that will infect infinitely many new sites on the line $\ell_u$ (with the help of $\H_u$), if $\alpha(u) < \infty$ then one can show that only boundedly many copies of this set are needed to infect the whole line. Moreover, if either $\alpha^-(u)$ or $\alpha^+(u)$ is infinite then no finite set of sites (together with $\H_u$) is sufficient to infect the whole of $\ell_u$.

Next let us see why the quantity $\alpha=\alpha(\U)$ defined in \eqref{eq:alphaU} is the constant one should expect to see in the exponent of the critical probability in Theorem~\ref{th:main}. In order to do this, we need the definition of a droplet, which is just a polygon in $\Z^2$. Droplets will be our means of controlling the growth of a set of infected sites.

\begin{definition}\label{def:droplets}
Let $\T \subset \Q_1$ be finite. A \emph{$\T$-droplet} is a non-empty set of the form 
\[
D = \bigcap_{u \in \T} \big( \H_u + a_u \big)
\]
for some collection $\{a_u \in \Z^2 : u \in \T\}$.
\end{definition}

Reinterpreted in terms of droplets, the definition of $\alpha$ in \eqref{eq:alphaU} is equivalent to the statement that there exist finite $\T$-droplets for some set $\T\subset\stab$ such that $\alpha(u)\geq\alpha$ for all $u\in\T$, but that the same is not true if $\alpha$ is replaced by any larger quantity. In other words, any finite set of infected sites is contained in a closed droplet such that a `cluster' of at least $\alpha$ sites is needed to create non-localized new infections. 
In the other direction, the condition that there exists an open semicircle $C\subset S^1$ such that every $u\in C$ has difficulty at most $\alpha$, which is implied by the definition of $\alpha$, means that there is an interval of directions having difficulty at most $\alpha$ just large enough for there to exist infinite sequences of nested droplets such that it is possible to grow between consecutive droplets using only sets of $\alpha$ sites. (Note that in the general model, unlike in symmetric bootstrap models, droplets do not necessarily grow in all directions.)

Before continuing with the outline of the proof, let us record two conventions that we use throughout the paper. First, $\U$ will always denote a fixed critical update family, unless explicitly stated otherwise. Using results from Section~\ref{se:simpleprops}, below, this is equivalent to the statement that $1\leq\alpha<\infty$. Second, $A$ will always denote a $p$-random subset of either $\Z^2$ or $\Z_n^2$. We emphasize that, since we will usually be working with droplets on a scale much smaller than $n$, most of the time we will be able to work in $\Z^2$, and the reader should assume throughout that we are in the infinite setting unless otherwise stated.


\subsection{Upper bounds}

The overall approach of the proofs of the upper bounds mirrors that of previous works (see for example~\cite{AL,GG99,vEH}). First we obtain a lower bound of the form $\exp\big(- p^{-\alpha + o(1)} \big)$ for the probability that a droplet at a particular intermediate scale (which is roughly $p^{-\Theta(1)}$) is (almost) internally filled, where `internally filled' is defined as follows. 

\begin{definition}\label{de:ifill}
A set $X \subset \Z^2$ is \emph{internally filled} by $A$ if $X \subset [X \cap A]$. The event that $X$ is internally filled by $A$ is denoted $\ifill(X)$.
\end{definition}

We remark that in the older bootstrap percolation literature this event was referred to as $X$ being `internally spanned' by $A$. However, following~\cite{BBM3D,BBDCM}, we will reserve that term for a different notion, see Definition~\ref{de:ispan} below.  

As alluded to before Definition~\ref{de:ifill}, we will in fact usually show that certain droplets are not quite \emph{exactly} internally filled, but \emph{almost} internally filled, where we use this terminology informally to mean that sites within distance $O(1)$ may be used to help fill the droplet. The resulting loss of independence is not a problem, because the events are increasing and we bound them using Harris's inequality.

The key step in the proof is a bound of the form
\[
\Pr_p\Big( D_m\subset\big[D_{m-1}\cup(D_{m+1}\cap A)\big] \Big) \geq \big(1-(1-p^\alpha)^{\Omega(m)}\big)^{O(1)},
\]
where $D_0\subset D_1\subset\cdots$ is a certain sequence of nested droplets. This bound corresponds to the intuition that it is enough to find somewhere along each side of the droplet a bounded number of sets of $\alpha$ sites contained in $A$. Once we have this bound, we then deduce that with high probability there exists an internally filled copy of this intermediate droplet in $\Z_n^2$, and that with high probability this droplet grows to infect the whole torus. 

All of what we have just said assumes to some extent that the family is balanced. If it is unbalanced then the droplets in the nested sequence $(D_m)_{m=0}^\infty$ are somewhat different: the sides (in the directions of growth) cannot grow linearly with $m$ (as in the balanced case), but instead all have the same length, and as a consequence the droplets are much less `regular' (for example, the initial droplet has width $\lambda$ and height $\lambda p^{-\alpha} \log(1/p)$, where $\lambda$ is a large constant, while in the balanced case it has constant size). The growth also features an extra step, in which an extremely long rectangular droplet grows a triangle of infected sites on its side.

Two key deterministic properties of the growth process are needed to make the above ideas work, for both balanced and unbalanced families. The first we have already discussed: the statement that a bounded number of sets of $\alpha$ sites are enough to infect an entire new line; we refer to this principle as `voracity', see Section~\ref{sec:voracity}. The second is the ability to grow to the corners of droplets, not just to within a bounded distance of the corners. This is in general not possible with $\T$ a subset of the stable set $\stab$. However, using the idea of `quasi-stability' introduced in~\cite{BSUgen}, one can show that it can be done if a certain set of unstable directions is included in $\T$, see Section~\ref{sec:quasi} for the details. 

\subsection{The lower bound for balanced families}

The lower bound for balanced update families is also not too difficult, but again requires refined versions of arguments from~\cite{BSUgen}. In order to sketch the proof, let us first briefly recall the argument of Aizenman and Lebowitz~\cite{AL} for the 2-neighbour model. Their key lemma states that if $A \subset \Z_n^2$ percolates, then for every $1 \leq k \leq n$, there exists an internally filled rectangle of semi-perimeter between $k$ and $2k$. Using the simple and well-known extremal result that such an internally filled rectangle contains at least $k$ initially infected sites, the bound follows from a straightforward calculation.

The key lemma of Aizenman and Lebowitz is proved via the so-called `rectangles process', which is an algorithm for determining the exact closure of a finite set under the 2-neighbour process. The algorithm proceeds by breaking down the bootstrap process into steps, each of which corresponds to the joining of two nearby rectangles into a larger rectangle. (Note that rectangles are closed under the 2-neighbour process.) One significant obstacle in the analysis of the general model is the lack of a corresponding \emph{exact} algorithm. Our solution is to use a process analogous to the rectangles process but rather more complicated. This process is an adaptation of the `covering algorithm' of Bollob\'as, Smith and Uzzell~\cite{BSUgen}, and we use it in order to prove lemmas corresponding to those of~\cite{AL}. Roughly speaking, we shall treat clusters of $\alpha$ nearby sites as seeds, cover each with a small $\stab$-droplet, and combine them pairwise into larger droplets if they are sufficiently close to interact in the $\U$-bootstrap process. The crucial deterministic property of the covering algorithm is that the remaining infections (those not in $\alpha$-clusters) contribute a negligible amount to the set of eventually infected sites; this is proved in Lemma~\ref{le:cover}.

We remark that we will not need to use the fact that $\U$ is balanced in the proof outlined above; for unbalanced families, however, the bound given by this approach is not tight. The main challenge of this paper, carried out in Sections~\ref{se:approxfill} and~\ref{se:lowerunbal}, will be to prove a stronger lower bound for unbalanced families. 

\subsection{The lower bound for unbalanced families}\label{se:sketchunbal}

The proofs of the previous three parts of the theorem are essentially refinements of established techniques. For this final part of the theorem, however, these techniques do not seem to be useful, and instead we introduce several substantial new ideas, including iterated hierarchies, the $u$-norm, and icebergs (see below). We mention these ideas only briefly in this section, focusing instead on the broad structure of the proof, and on some of the most important definitions. A much fuller outline of the proof is given at the start of Section~\ref{se:lowerunbal} (see also Section~\ref{se:approxfill}).

The first observation we make (see Lemmas~\ref{le:ustar} and~\ref{le:stabu}) is that there exist opposite stable directions $u^*$ and $-u^*$ that both have difficulty at least $\alpha+1$. We set 
$$\stabu = \big\{ u^*,-u^*,u^l,u^r \big\},$$
where $u^l$ and $u^r$ are stable directions on different sides of $u^*$, each of difficulty at least $\alpha$, and we consider only $\stabu$-droplets. (The `$U$' in $\stabu$ is meant to remind the reader that $\U$ is unbalanced.) Let us rotate our perspective so that $u^*$ is vertical, and write $h(D)$ and $w(D)$ for the height and width of an $\stabu$-droplet respectively.

As in the balanced case, first we need an approximate rectangles process, which will allow us to say that if a large droplet is internally filled then it must contain droplets at all scales that are `approximately internally filled'. The covering algorithm is no longer useful to us because it is too crude to capture the biased nature of the geometry of unbalanced models. Instead we use a second algorithm, the `spanning algorithm', which is an adaptation of an idea introduced by Cerf and Cirillo in~\cite{CerfCir} and subsequently developed in~\cite{CerfManzo,BBM3D,BBDCM,DCH}. The algorithm uses the following notion of connectedness and the subsequent notion of being `internally spanned', which is an approximation to being internally filled.

\begin{definition}\label{de:connected}
Let $\kappa$ be a sufficiently large constant, to be defined explicitly in~\eqref{eq:kappa}. Define a graph $G_\kappa$ with vertex set $\Z^2$ and edge set $E(G_\kappa) = \big\{ xy : \|x-y\|_2 \leq \kappa \big\}$. We say that a set $S \subset \Z^2$ is \emph{strongly connected} if it is connected in the graph $G_\kappa$.
\end{definition}

We will also say that two sets of sites $X$ and $Y$ are \emph{strongly connected to each other} (or, sometimes, that $X$ \emph{is strongly connected to} $Y$), if 
there exist elements $x \in X$ and $y \in Y$ with $\|x-y\|_2 \leq \kappa$, and hope that this will cause no confusion.

\begin{definition}\label{de:ispan}
Let $\T \subset \stab$ be finite. A $\T$-droplet $D$ is \emph{internally spanned} by $A$ if there is a strongly connected set $L \subset [D\cap A]$ such that $D$ is the smallest $\T$-droplet containing $L$. We will write $\ispan(D)$ for the event that $D$ is internally spanned.\footnote{Note that the event $\ispan(D)$ also depends on $\T$. However, we will only use this notation when $\T = \stabu$, and so we trust this will therefore not cause the reader any confusion.}
\end{definition}

As noted above, many previous authors have used the term `internally spanned' to mean (what we refer to as) `internally filled'. We reemphasize that our terminology (which follows~\cite{BBM3D,BBDCM}, and seems to us more natural) is different.

The spanning algorithm allows us to break down the formation of an internally spanned droplet into intermediate steps in the same way that the original rectangles process allowed Aizenman and Lebowitz~\cite{AL} to break down the formation of an internally filled droplet. Using the spanning algorithm we are able to say that if a large droplet is internally spanned, then it contains internally spanned droplets at all smaller scales. The scale we are particularly interested in is the `critical' scale, which for unbalanced models has the following specific meaning.

\begin{definition}\label{de:critical}
Let $\U$ be unbalanced and let $\xi>0$ be a small positive constant. An $\stabu$-droplet $D$ is \emph{critical} if either of the following conditions holds:\smallskip
\begin{enumerate}
\item[$(T)$] $w(D) \leq 3 p^{-\alpha-1/5}$ and $\frac{\xi}{p^\alpha}\log\frac{1}{p}\leq h(D)\leq\frac{3\xi}{p^\alpha}\log\frac{1}{p}$; \smallskip
\item[$(L)$] $p^{-\alpha-1/5}\leq w(D)\leq 3p^{-\alpha-1/5}$ and $h(D)\leq\frac{\xi}{p^\alpha}\log\frac{1}{p}$. 
\end{enumerate}
\end{definition}

Why might this be the right definition? It is certainly not surprising that the droplet should be long and thin; this is the nature of unbalanced growth, as suggested by the proof of the upper bound in Section~\ref{se:upperunbal}. The height $h=\frac{\xi}{p^\alpha}\log\frac{1}{p}$ is such that an initial rectangle of height $h$ and constant width will fail to grow sideways (that is, perpendicular to $u^*$) by a constant distance with probability roughly $p^{O(\xi)}$, and therefore one would expect the rectangle to grow sideways only to distance $p^{-O(\xi)}$. The width $w=p^{-\alpha-1/5}$ is such that the probability the rectangle grows to distance $w$ is sufficiently small to compensate for the number of choices for the initial rectangle. The reason for there being two types of critical droplet is that the spanning algorithm cannot control the width and the height of the critical droplet simultaneously.

In order to bound the probability that a critical droplet $D$ is internally spanned, we shall show that, if the droplet is of type $(T)$, then it is unlikely that $[D\cap A]$ contains a strongly connected set joining the $(-u^*)$-side of $D$ to the $u^*$-side, while if it is of type $(L)$, then instead it is unlikely that $[D\cap A]$ contains a strongly connected set joining the $u^l$-side to the $u^r$-side. (The $u$-side of a droplet is defined precisely below.) These events are called `vertical crossings' and `horizontal crossings' respectively.

There are several complications that occur while bounding the probabilities of such crossings. Consider first vertical crossings, and note that, since $\alpha(u^*) \geq \alpha + 1$, we have either $\min\big\{ \alpha^+(u^*),\alpha^-(u^*) \big\} \geq \alpha + 1$, or 
\begin{equation}\label{eq:unbal:driftcase}
\max\big\{ \alpha^+(u^*), \alpha^-(u^*) \big\} = \infty,
\end{equation}
and similarly for $-u^*$. Since the former case is much easier to handle, let us assume in this sketch that~\eqref{eq:unbal:driftcase} holds. If $u \in \{u^*,-u^*\}$ and $\max\{ \alpha^+(u), \alpha^-(u) \} = \infty$, then we call $u$ a \emph{drift direction}, and if either $u^*$ or $-u^*$ is a drift direction, then we say that the model \emph{exhibits drift}.

For concreteness, suppose that $\alpha^-(u^*) = \infty$ and $\alpha^+(u^*) = 1$. Since we have a pair $\{u^*,-u^*\}$ of opposite stable directions, we may partition the droplet $D$ into many smaller sub-droplets of the same width, and bound the probability that each is vertically crossed (possibly with help from above and below) independently, since these events depend on disjoint sets of infected sites. In order to bound these crossing probabilities, we need several new ideas. First, we need a method of controlling the range of the $\U$-bootstrap process assisted by a half-plane. We achieve this by introducing (in Section~\ref{se:iceberg}) a third algorithm for approximating the closure of a set of sites, which we call the `$u$-iceberg algorithm'. Second, we need a close-to-best-possible bound for the probability that certain smaller sub-droplets are internally spanned (following~\cite{BBM3D}, we call these sub-droplets `savers'). In order to obtain such a bound, we induct on the size of the droplet being crossed. (This means the proof for vertical crossings at a given scale depends on us having obtained sufficiently strong bounds for both vertical and horizontal crossings at the scale below.) Finally, we need to deal with the `stretched geometry' of drift models; we do so by introducing a family of norms (the `$u$-norms') that compress this geometry until it resembles Euclidean space, and we also introduce a new concept of (`weak') connectedness; see Sections~\ref{se:range} and~\ref{se:crossings} respectively for the details.

Now consider horizontal crossings, and observe that we no longer have symmetry, since $-u^l$ and $-u^r$ are in general not stable directions. This prevents us from partitioning into sub-droplets as with vertical crossings, and to overcome this obstacle we use the `hierarchy method', which was introduced by Holroyd in~\cite{Hol} and subsequently developed in~\cite{BBDCM,BBM3D,DCvE,GravHolMor}. We would like to emphasize that the reason for our use of hierarchies is different to that of all previous works: here, the reason is the lack of symmetry between $u^l$ and $u^r$ (which is also why we do not need them for vertical crossings); previously the reason has been to prove sharp thresholds for critical probabilities in symmetric settings.

In order to use hierarchies, we need three further ingredients: a bound on the probability that `seeds' (which are small sub-droplets) are internally spanned; a bound on the probabilities of certain ($p$ times shorter) horizontal crossing events; and a bound on the number of hierarchies with a given number of seeds. For these we use the induction hypothesis (once again, we need sufficiently strong bounds for both vertical and horizontal crossings at the scale below), and the method described above for vertical crossings (together with an additional argument, see Section~\ref{sec:growing:droplets}). Since our use of induction on the size of the droplet amounts to iterating the above argument $\alpha$ times, we refer to this as the `method of iterated hierarchies'.

\subsection{Basic facts about $\U$-bootstrap percolation}\label{se:simpleprops}

The $\U$-bootstrap process exhibits a number of particularly simple and elegant properties, some of which we now recall from~\cite{BSUgen}. We begin with a description of the stable set $\stab$. We write $[v,w]$ for the closed interval of directions between $v$ and $w$ taken anticlockwise starting from $v$, and we say that $[v,w]$ is \emph{rational} if $\{v,w\} \subset \Q_1$.

\begin{lemma}[Theorem 1.10 of \cite{BSUgen}]\label{le:stableset}
The stable set $\stab$ is a finite union of rational closed intervals of $S^1$. 
\end{lemma}

The converse to Lemma~\ref{le:stableset} is also true (and is part of Theorem~1.10 of~\cite{BSUgen}): if $\stab\subset S^1$ is any set consisting of a finite union of rational closed intervals, then there exists an update family $\U$ such that $\stab=\stab(\U)$. We shall not use this converse.

Another simple but important fact, proved in~\cite{BSUgen}, is that if $u$ is not stable then $\H_u$ grows to fill the whole of $\Z^2$.
 
\begin{lemma}[Lemma~3.1 of~\cite{BSUgen}]\label{lem:all:or:nothing}
If $u \not\in \stab$, then $[\H_u] = \Z^2$. 
\end{lemma}

Thus for every $u\in S^1$ we have the dichotomy $[\H_u] \in \big\{\H_u,\Z^2\}$. The following simple properties of directions of infinite difficulty, which were proved in~\cite{BSUgen}, will also be useful. For completeness, we sketch the proofs.

\begin{lemma}\label{le:alphainfinity}
If\/ $[v,w] \subset \stab$ is a connected component of\/ $\stab$ and $u \in [v,w] \cap \Q_1$, then 
\begin{equation}\label{eq:alphainfinity}
\alpha^+(u) < \infty \; \Leftrightarrow \; u = v \qquad \textup{and} \qquad \alpha^-(u) < \infty \; \Leftrightarrow \; u = w.
\end{equation}
In particular, $0 < \alpha(u) < \infty$ if and only if $u$ is an isolated point of $\stab$.
\end{lemma}

\begin{proof}
We will show first that $\alpha^+(v) < \infty$. To do so, observe that (since $[v,w]$ is a component of $\stab$) there exist unstable directions arbitrarily close to (and to the right of) $v$. Choose such a direction $v'$ sufficiently close to $v$, and choose $X \in \U$ such that $X \subset \H_{v'}$. Now, since the elements of $X$ all lie within a finite distance of the origin, and $v'$ was chosen sufficiently close to~$v$, it follows that $X \subset \H_v \cup \ell_v^-$, and therefore the set $[\H_v \cup Z] \cap \ell_v^+$ is infinite for any sufficiently large (finite) set $Z$ of consecutive sites of~$\ell_v$. By Definition~\ref{de:alpha}, this implies that $\alpha^+(v) < \infty$, as claimed.

In order to prove that $\alpha^+(u) = \infty$ for every $u \in (v,w] \cap \Q_1$, we will first show that (for each such $u$) there exists $u' \in (v,u)$ such that $\H_u \cup \H_{u'}$ is closed under the $\U$-bootstrap process. To do so, simply choose $u'$ closer to $u$ than any $u'' \in S^1 \setminus \{u\}$ perpendicular to a vector in the set
$$\Big\{ x - y \,:\, x,y \in \bigcup_{X \in \U} X \cup \{\0\}, \, x \ne y \Big\}.$$
That we can do this follows easily from the fact that $\U$ is a finite collection of finite sets. Now, suppose that there exists a rule $X \in \U$ such that $X \subset \H_u \cup \H_{u'}$. Since $u,u' \in \S$, there must exist $x,y \in X$ with $x \not\in \H_u$ and $y \not\in \H_{u'}$. But now $x - y$ is perpendicular to a vector in the interval $(u',u)$, which contradicts our choice of $u'$.

Now, let $u \in (v,w] \cap \Q_1$, and choose $u' \in (v,u)$ such that $\H_u \cup \H_{u'}$ is closed. Observe that for any finite set $Z \subset \Z^2$, there exists $y \in \ell_u$ such that $y + Z \subset \H_{u'}$, and therefore $[\H_u \cup (y + Z)] \subset \H_u \cup \H_{u'}$. It follows that $[\H_u \cup Z] \cap \ell_u^+$ is finite, and hence $\alpha^+(u) = \infty$, as required. The second claim in~\eqref{eq:alphainfinity} now follows by symmetry. 

For the final assertion, note that $\alpha(u) > 0$ if and only if $u \in \S$ (by Definition~\ref{de:alpha} and Lemma~\ref{lem:all:or:nothing}), and that if $u \not\in \Q_1$, then $\alpha(u) \in \{0, \infty\}$ (by Definition~\ref{de:alpha}) and $u$ is not an isolated point of $\stab$ (by Lemma~\ref{le:stableset}). On the other hand, if $u \in \S \cap \Q_1$, then by~\eqref{eq:alphainfinity} we have $\alpha(u) < \infty$ if and only if the component $[v,w]$ of $\S$ containing $u$ satisfies $v = w$, i.e., if and only if $u$ is an isolated point of $\S$.
\end{proof}

Let us note, for emphasis, that the proof above also implies the following lemma from~\cite{BSUgen}; we again give the details for completeness.

\begin{lemma}[Lemma~5.2 of~\cite{BSUgen}]\label{le:leftrightcluster}
If $u \in \Q_1$ and $\alpha^-(u) < \infty$, then there exists $X \in \U$ such that $X \subset \H_u \cup \ell_u^+$, and hence $\ell_u \subset [\H_u \cup \ell_u^+]$. 
\end{lemma}

\begin{proof}
By Lemma~\ref{le:alphainfinity}, there exist unstable directions arbitrarily close to (and to the left of) $u$. Choose such a direction $v$ sufficiently close to $u$, and choose $X \in \U$ such that $X \subset \H_v$. Since the elements of $X$ all lie within a bounded distance of the origin, it follows that $X \subset \H_u \cup \ell_u^+$, as required.
\end{proof}

We are now in a position to deduce the existence of opposite stable directions $u^*$ and $-u^*$ claimed earlier for unbalanced families $\U$.

\begin{lemma}\label{le:ustar}
Let $\U$ be an unbalanced critical update family. Then there exists $u^* \in \Q_1$ such that 
\[
\min\big\{ \alpha(u^*), \alpha(-u^*) \big\} \geq \alpha + 1.
\]
\end{lemma}

\begin{proof}
By the definition of $\alpha$, there exists an open semicircle $C \in \C$ such that $\alpha(u) \leq \alpha$ for every $u \in C$. Moreover, since $\U$ is critical we have $\alpha < \infty$. Thus, if one of the endpoints of $C$ has difficulty at most $\alpha$, then it is either unstable, or an isolated point of $\stab$, by Lemma~\ref{le:alphainfinity}. Hence, rotating $C$ slightly, we obtain a closed semicircle $C'$ such that $\alpha(u) \leq \alpha$ for all $u\in C'$. But this contradicts our assumption that $\U$ is unbalanced, and so both endpoints of $C$ have difficulty at least $\alpha + 1$. 
\end{proof}

\subsection{Definitions and notation}\label{se:defs}

In this subsection we collect for ease of reference various conventions, definitions and notation that we shall use throughout the paper. 

\subsubsection{Constants, and asymptotic notation}

All constants, including those implied by the notation $O(\cdot)$, $\Omega(\cdot)$ and $\Theta(\cdot)$, are quantities that may depend on $\U$ (and other quantities where explicitly stated) but not on $p$. The parameter $p$ will always be assumed to be sufficiently small relative to all other quantities. Our asymptotic notation is mostly standard, although we just remark that if $f$ and $g$ are positive real-valued functions of $p$, then we write $f(p)=\Omega\big(g(p)\big)$ if $g(p)=O\big(f(p)\big)$, and we write $f(p)=\Theta\big(g(p)\big)$ if both $f(p)=O\big(g(p)\big)$ and $g(p)=O\big(f(p)\big)$. Furthermore, if $c_1$ and $c_2$ are constants, then $c_1 \gg c_2 \gg 1$ means that $c_2$ is sufficiently large, and $c_1$ is sufficiently large depending on $c_2$, and $1 \gg c_1 \gg c_2 > 0$ means that $c_1$ is sufficiently small, and $c_2$ is sufficiently small depending on $c_1$. (This last piece of notation is somewhat non-standard; we trust it will not cause any confusion.)


\subsubsection{Measuring sizes and distances}\label{sec:measuring}



The unadorned norm $\|\cdot\|$ always denotes the Euclidean norm on $\R^2$, and $\langle\cdot,\cdot\rangle$ always denotes the Euclidean inner product. As remarked above, in Section~\ref{se:approxfill} we will define a family of norms on $\R^2$ called the `$u$-norms', which will be signified with a subscript $u$ thus: $\|\cdot\|_u$.

Now, for $u\in S^1$ and a compact set $K\subset \R^2$, define the \emph{$u$-projection} of $K$,
\begin{equation}\label{eq:uproj}
\pi(K,u) := \max\big\{ \< x-y,u\> : x,y\in K \big\}.
\end{equation}
Also, let
\[
\diam(K) := \max\big\{ \|x-y\| : x,y\in K \big\} = \max\big\{ \pi(K,u) : u\in S^1 \big\}
\]
be the \emph{diameter} of $K$. Owing to the biased nature of the geometry, in the unbalanced setting the diameter is usually not a useful measure of the size of $K$. Instead, we work with the \emph{height}
\[
h(K) := \pi(K,u^*) = \max\big\{ \< x-y,u^*\> : x,y\in K \big\},
\]
and the \emph{width}
\[
w(K) := \pi(K,w^\perp) = \max\big\{ \< x-y,w^\perp\> : x,y\in K \big\},
\]
where $\{u^*,-u^*\} \subset \stab$ is the pair of opposite stable directions with difficulty strictly greater than $\alpha$ given by Lemma~\ref{le:ustar}, and $w^\perp\in S^1$ is either of the two unit vectors that are orthogonal to~$u^*$. We will also make frequent use of the following constant, which we think of as being the `diameter' of $\U$:
\begin{equation}\label{eq:nu}
\nu = \nu(\U) := \max \Big\{ \|x - y\| : x,y \in X \cup \{\0\}, \, X \in \U \Big\}.
\end{equation}
We will define another constant $\rho$, which captures a different aspect of the ``range" of the $\U$-bootstrap process, in Section~\ref{se:approxfill}. 

Occasionally we shall want to talk about the distance between a site and a set of sites, or between two sets of sites. We use the following standard conventions:
\begin{align*}
\|x-Y\| &:= \min \big\{ \|x-y\| : y\in Y \big\}, \\
\text{and} \qquad \|X-Y\| &:= \min \big\{ \|x-y\| : x\in X,\, y\in Y \big\},
\end{align*}
whenever $X$ and $Y$ are finite subsets of $\Z^2$. We also use analogous conventions for other measures of distance, such as the `$u$-norms' and inner products.



\subsubsection{Subsets of the plane}


If $u \in \Q_1$, then the collection of non-empty discrete lines
\[
\Big\{ \big\{ x \in \Z^2 : \langle x-a,u \rangle = 0 \big\} : a \in \Z^2 \Big\}
\]
is a discrete set, naturally indexed by $\Z$. Thus, we may set $\ell_u(0) := \ell_u$ and (for each $j\in \Z$) let $\ell_u(j)$ denote the $j$th non-empty discrete line in the direction of $u$. 

For each $u \in S^1$ and $a \in \R^2$, we define the discrete half-planes
\[
\H_u(a) := \big\{x \in \Z^2 : \langle x-a, u \rangle < 0 \big\}.
\]
If $a \in \Z^2$ then we have $\H_u(a)=\H_u+a$, but this is false otherwise (since $\H_u \subset \Z^2$). Observe (cf. Definition~\ref{def:droplets}) that a $\T$-droplet is a non-empty set of the form 
\[
D = \bigcap_{u \in \T} \H_u(a_u)
\]
for some collection $\{a_u \in \Z^2 : u \in \T\}$. For each $u \in \T$, the \emph{$u$-side} of a $\T$-droplet $D$ is defined to be the set
\begin{equation}\label{def:uside}
\partial(D,u) := D \cap \ell_u(i),
\end{equation}
where $i$ is maximal subject to the set being non-empty. Finally, note that we can consider droplets (even those with diameter larger than $n$) as subsets of $\Z_n^2$ by taking all $x = (x_1,x_2) \in \Z_n^2$ such that $(x_1 + in,x_2+jn) \in D$ for some $i,j \in \Z$.

\subsection{Probabilistic lemmas}\label{se:prob}

We end the section by recalling the correlation inequalities of Harris~\cite{Harris}, and van den Berg and Kesten~\cite{BK}. For definitions of increasing events and disjoint occurrence, and for proofs of both inequalities, see~\cite{BRPerc}.

\begin{lemma}\label{le:harris}\emph{(Harris's inequality)}
If $\mathcal{A}$ and $\mathcal{B}$ are increasing events then
\[
\Pr_p(\mathcal{A} \cap \mathcal{B}) \geq \Pr_p(\mathcal{A}) \cdot \Pr_p(\mathcal{B}).
\]
\end{lemma}

We write $\mathcal{A} \circ \mathcal{B}$ for the event that $\mathcal{A}$ and $\mathcal{B}$ occur disjointly.

\begin{lemma}\label{le:bk}\emph{(The van den Berg--Kesten inequality)}
If $\mathcal{A}$ and $\mathcal{B}$ are increasing events then
\[
\Pr_p(\mathcal{A} \circ \mathcal{B}) \leq \Pr_p(\mathcal{A}) \cdot \Pr_p(\mathcal{B}).
\]
\end{lemma}

We shall apply Harris's inequality frequently throughout the paper, but the van den Berg--Kesten inequality only once, in the proof of Lemma~\ref{le:boundoverH}.

\section{Voracity and quasi-stability}\label{se:upperprelims}

In Section~\ref{se:sketch} we mentioned that there were two important deterministic concepts that we needed in order to make our upper bound proofs work. These were the notions of `voracious sets' and `quasi-stable directions'. In this section we introduce and develop these ideas, in preparation for the proofs of the upper bounds of Theorem~\ref{th:main} in the two sections to follow.

\subsection{Voracity and the infection of new lines}\label{sec:voracity}

We begin by studying sets of infected sites that are sufficient for stable half-planes to grow. Recall that $\U$ is a critical update family, fixed throughout, and that therefore $1 \leq \alpha < \infty$.

\begin{definition}
Let $u \in \Q_1$ with $\alpha(u) < \infty$, and let $Z \subset \Z^2$ be a set of size $|Z| = \alpha(u)$. If $[\H_u \cup Z]\cap\ell_u$ is infinite, then we say that $Z$ is \emph{voracious for $u$}.
\end{definition}

The definition of $\alpha(u)$ implies then there exists at least one voracious set for every $u \in \stab$ with $\alpha(u) < \infty$. We would like to show (see Lemma~\ref{le:line}, below) that a bounded number of voracious sets on the $u$-side of a (finite) droplet $D$ are sufficient to infect all but a bounded number of sites on the line adjacent to the $u$-side of $D$. The following definition will be useful.

\begin{definition}
A \emph{homothetic copy} of a set $S$ is a set 
$$Y = a + kS = \big\{ y \in \Z^2 : y = a + k b \textup{ for some } b \in S \big\}$$
for some $a \in \Z^2$ and non-zero $k \in \Z$.  
\end{definition}

Note that a homothetic copy of $\ell_u^+$ with $a \in \ell_u$ is an infinite subset of the line $\ell_u$. As a warm-up for the (slightly technical) finite setting, let's begin by proving the infinite version of the statement we require.

\begin{lemma}\label{lem:periodic:infinite}
Let $u \in \Q_1$ be such that $\alpha(u) < \infty$ and let $Z$ be voracious for $u$. Then $[\H_u\cup Z] \cap \ell_u$ contains a homothetic copy of $\ell_u^+$.
\end{lemma}

\begin{proof}
We may assume that $u$ is stable, since otherwise we are done by Lemma~\ref{lem:all:or:nothing}, and that $[\H_u\cup Z]$ contains infinitely many sites of $\ell_u^+$, since $Z$ is voracious. Since $Z$ is finite, there exists $a\in\Z^2$ such that $[\H_u\cup Z] \subset \H_u(a)$. Recall from~\eqref{eq:nu} the definition of $\nu$, and partition $\H_u(a) \setminus \H_u$ into disjoint congruent rectangles $\ldots,R_{-1},R_0,R_1,\ldots$  of width $s > 2\nu$, each of which is a translation of $R_0$, with $R_i$ immediately to the right of $R_{i-1}$ for each $i \in \Z$, and such that $Z \subset R_0$, noting that this is possible if $s$ is sufficiently large. Set $L_i = R_i \cap [\H_u \cup Z]$ for each $i \in \Z$. 

We claim that $L_{i+1} = R_{i+1} \cap [\H_u \cup L_i]$ for every $i \geq 0$. To see this, note first that $[\H_u \cup L_i] \subset [\H_u \cup Z]$, since $L_i \subset [\H_u \cup Z]$, so $R_{i+1} \cap [\H_u \cup L_i] \subset L_{i+1}$. For the other inclusion we argue by contradiction, considering the first time at which the process starting from $\H_u \cup Z$ infects a point of $\big( \bigcup_{j = i+1}^\infty R_j \big) \setminus [\H_u \cup L_i]$. Since $L_i$ is the \emph{final} configuration inside $R_i$, and $s > 2\nu$, it is easy to see that no such time can occur.

It follows that $L_{i+1}$ depends only on $L_i$, and therefore, since there are only finitely many possible configurations for $L_i$, there exist $j \ge 0$ and $r \ge 1$ such that $L_j \equiv L_{j+r}$ (meaning $L_{j+r}$ is obtained by translating $L_j$ by distance $rs$ perpendicular to $u$). Hence
\[
\big(L_j,\dots,L_{j+r-1}\big) \equiv \big(L_{j+r},\dots,L_{j+2r-1}\big) \equiv \big(L_{j+2r},\dots,L_{j+3r-1}\big) \equiv \cdots,
\]
and since $[\H_u \cup Z] \cap \ell_u^+$ is infinite, this is sufficient to prove the lemma.
\end{proof}

By taking suitable translates of the voracious set $Z$ in Lemma~\ref{lem:periodic:infinite}, it is clear that we can infect the whole of either $\ell_u^+$ or $\ell_u^-$. We can then use Lemma~\ref{le:leftrightcluster} to return back along the line and infect the rest of $\ell_u$.

Let us now turn to the finite setting, applicable to the growth of a new row on the side of a droplet. Recall that we denote by $\partial(D,u)$ the $u$-side of a $\T$-droplet $D$. 
In the lemma below we will also need the following notion: define the \emph{$u$-outside} of a $\T$-droplet $D$ to be the set $\partial_\circ(D,u)$ of points of $\ell_u(i+1)$ that lie within distance 1 of the convex hull of $\partial(D,u) = D \cap \ell_u(i)$. Let us also say that a set $Z$ lies \emph{above} the $u$-side of $D$ if its orthogonal projection onto the continuous line (perpendicular to $u$) through $\partial(D,u)$ is contained in the convex hull of $\partial(D,u)$. 

The following lemma says that a bounded number of voracious sets are sufficient (together with $D$) to infect all but a bounded number of sites of the $u$-outside of a droplet $D$, and moreover we may choose any suitable translate of each voracious~set. 

\begin{lemma}\label{le:line}
Let $\T \subset \Q_1$ be a finite set, let $u \in \T$ satisfy $\alpha(u) < \infty$, and let $Z$ be voracious for $u$. Then there exist $\mu > 0$, $r \in \N$, and $\0 \ne b \in \ell_u$, such that for every $\T$-droplet $D$, there exist $a_1,\dots,a_r \in \partial_\circ(D,u)$ such that the following holds.

Suppose that $k_1,\dots,k_r \in \Z$ are such that $Z + a_j + k_j b$ is above the $u$-side of $D$, and at distance at least $\mu$ from the corners of $D$, for every $1 \le j \le r$. Then the set
$$\big[ D \cup (Z+a_1+k_1b)\cup\dots\cup(Z+a_r+k_rb) \big]$$ 
contains all elements of $\partial_\circ(D,u)$ at distance at least $\mu$ from the corners.
\end{lemma}

\begin{proof}
First, note that, by translating $D$ by an element of $\Z^2$ if necessary, we may assume that $\partial_\circ(D,u)$ contains the origin. Now, choose a constant $\mu_0 > 0$ to be sufficiently large so that, for every $X \in \U$, we have $y + (X \cap \H_u) \subset D$ for every $\T$-droplet $D$, and every $y \in \partial_\circ(D,u)$ at distance at least $\mu_0$ from the corners. Note that (without loss of generality) we may assume that $[\H_u\cup Z]$ contains infinitely many sites on the line $\ell_u^+$. Define the sequence $\ldots,R_{-1},R_0,R_1,\ldots$ of rectangles (each of constant width $s > 2\nu$) as in the proof of Lemma~\ref{lem:periodic:infinite}, and set 
$$L_i = R_i \cap [ \H_u \cup Z ]$$ 
for each $i \in \Z$. Recall that $Z \subset R_0$, and define $t_0 = \min\big\{ t \ge 0 : L_0 \subset ( \H_u \cup Z )_t \big\}$, i.e., the number of steps of the $\U$-bootstrap process it takes to infect $L_0$, starting from $\H_u\cup Z$. Since `information' can only travel distance $\nu$ in one step of the process, it follows that if $R_0$ is at distance at least $t_0 \nu + \mu_0$ from the corners of $D$, then $L_0 \subset [ D \cup Z ]$. Next, for each $i \geq 1$, define 
$$t_i = \min\big\{ t \ge 0 : L_i \subset ( \H_u \cup L_{i-1} )_t \big\}$$
i.e., the number of steps of the $\U$-bootstrap process it takes to infect $L_i$, starting from $\H_u\cup L_{i-1}$. Note that $t_i$ is finite, and moreover, since the $L_j$ are periodic there exists a constant $T$ such that $t_i \le T$ for every $i \ge 1$. Therefore, if $L_{i-1} \subset [ D \cup Z ]$ and $R_i$ is at distance at least $T \nu + \mu_0$ from the corners of $D$, then $L_i \subset [ D \cup Z ]$. 

It follows that $[ D \cup Z ] \cap \partial_\circ(D,u)$ contains the intersection with $\partial_\circ(D,u)$ of a homothetic copy of $\ell_u^+$ with bounded difference. More precisely, there exists $a \in \partial_\circ(D,u)$ and $b \in \N$, where $\| a - Z \|$ and $b$ are both at most some constant depending on $\U$, $\T$, $u$ and $Z$, but not on $D$, such that $[ D \cup Z ]$ contains every element of $a + b \ell_u^+$ that is in $\partial_\circ(D,u)$, and at distance at least $T \nu + \mu_0$ from the corners of~$D$. 

Hence there exist $b \in \ell_u$, $r \in \N$, and $a_1,\dots,a_r \in \partial_\circ(D,u)$, such that the following holds: if $k_1,\dots, k_r \in \Z$ are such that the set $Z + a_j + k_j b$ is above the $u$-side of $D$ and sufficiently far from the corners of $D$ for every $1 \le j \le r$, then the set 
$$Y := \big[ D \cup (Z+a_1+k_1b)\cup\dots\cup(Z+a_r+k_rb) \big]$$ 
contains $\nu$ consecutive elements of $\partial_\circ(D,u)$ at distance at least $\mu_0$ from the corners. Now, by Lemma~\ref{le:leftrightcluster} there exist update rules $X^+$ and $X^-$ contained in $\H_u \cup \ell_u^+$ and $\H_u \cup \ell_u^-$ respectively. Note that $X^+ \setminus \H_u$ is contained in the first $\nu$ sites of $\ell_u^+$, and similarly for $X^- \setminus \H_u$ and $\ell_u^-$. Hence $Y$ in fact contains all elements of $\partial_\circ(D,u)$ that are at distance at least $\mu_0$ from the corners, as required.
\end{proof}

As a consequence of Lemma~\ref{le:line}, one would expect that a $\T$-droplet $D$ would `grow by one step in direction $u$' with probability at least
\[
\big( 1 - ( 1 - p^\alpha )^{\Omega(m)} \big)^{O(1)},
\]
where $m$ is the length of the side of $D$ corresponding to $u$. This is almost true; however, as the presence of the constant $\mu$ in Lemma~\ref{le:line} suggests, we have a problem near the corners of $D$: we may need sites not in $D$ but still below the (extended) $u$-side of $D$ in order to infect the last $O(1)$ sites. We resolve this problem using another idea from~\cite{BSUgen}: that of quasi-stable directions.

\subsection{Quasi-stability}\label{sec:quasi}

In many of the simpler bootstrap models, the droplets used as bases for growth are taken with respect to the set of stable directions. Droplets for the 2-neighbour model are rectangles -- or, put another way, they are taken with respect to the set $\stab=\{e_1,-e_1,e_2,-e_2\}$ of stable directions. Similarly, for balanced threshold models with symmetric star-neighbourhoods\footnote{We say that $Y$ is a symmetric star-neighbourhood if $x \in Y$ implies that that $-x \in Y$, and moreover that every vertex of $\Z^2$ on the straight line between $x$ and $-x$ is in $Y$.} droplets can be taken with respect to the set of stable directions, and the droplets  are therefore $2k$-gons consisting of pairs of parallel sides, for some $k \geq 2$. In this case $\stab$-droplets are suitable bases for growth because, when new infections spread in both directions along each edge of the droplet, the set that results is a new, slightly larger droplet. 

The same is not true in general: indeed, as noted above, we may fail to infect some of the sites near the corners of $D$ due to boundary effects. 
The solution to this problem, as used by Bollob\'as, Smith and Uzzell~\cite{BSUgen}, was to introduce a number of \emph{quasi-stable directions}, which are not stable directions, but which nevertheless are treated as such. Thus, droplets are taken with respect to a certain superset of the stable set. For a comprehensive discussion of quasi-stability, we refer the reader to Section~5.1 of~\cite{BSUgen}.

The next lemma is Lemma~5.3 of~\cite{BSUgen}. Since the lemma is so fundamental to the proofs of the upper bounds of Theorem~\ref{th:main}, we give the (short) proof here (which is also similar to that of Lemma~\ref{le:alphainfinity}) in full. Recall that $[u,v]$ denotes the interval of directions in $S^1$ between $u$ to $v$, taken anticlockwise starting from $u$. Given a set $\T \subset S^1$, we say that $u$ and $v$ are \emph{consecutive} elements of $\T$ if $u \neq v$ and $\T \cap [u,v] = \{u,v\}$ (note that the order of $u$ and $v$ matters in this definition). 

\begin{figure}[ht]
  \centering
  \begin{tikzpicture}[>=latex]
    \node[circle,fill,inner sep=0pt,minimum size=0.15cm] at (0,0) {};
    \draw (30:1) -- (-150:6) ++(-150:0.4) node {$\ell_v$};
    \draw (160:1) -- (-20:6) ++(-20:0.4) node {$\ell_u$};
    \draw[->] (-150:5) -- ++(120:1);
    \path (-150:5) ++(120:1.4) node {$v$};
    \draw[->] (-20:5) -- ++(70:1);
    \path (-20:5) ++(70:1.4) node {$u$};
    \draw[densely dashed] (-175:4) -- (5:4) ++(5:0.4) node {$\ell_w$};
    \draw[->] (-175:2) -- ++(95:1);
    \path (-175:2) ++(95:1.4) node {$w$};
    \draw[densely dashed] (175:0.7) -- (-5:2.5);
    \node at (-5:2.9) {$x$};
    \draw[->] (-5:2.5) -- ++(85:1);
    \path (-5:2.5) ++(85:1.4) node {$w'$};
    \node[circle,fill,inner sep=0pt,minimum size=0.15cm] at (-5:2.5) {};
    \node at (-85:2) {$\big(\H_u\cup\ell_u\big)\cap\big(\H_v\cup\ell_v\big)$};
  \end{tikzpicture}
  \caption{The proof of Lemma~\ref{le:quasi}. Since $w$ is unstable, there exists $X \in \U$ with $X\subset\H_w$; if $x\in X$ lies in the region between $\ell_w$ and $\ell_u$, then the direction $w'$ would be in $\qquasi$, by construction, which contradicts $u$ and $v$ being consecutive in $\stab\cup\qquasi$. Thus $X \subset (\H_u\cup\ell_u\big)\cap\big(\H_v\cup\ell_v\big)$, as required. Note that the figure is completely general: $[u,v]$ is at most a semicircle, by the definition of $\qquasi$.}\label{fi:quasi}
\end{figure}
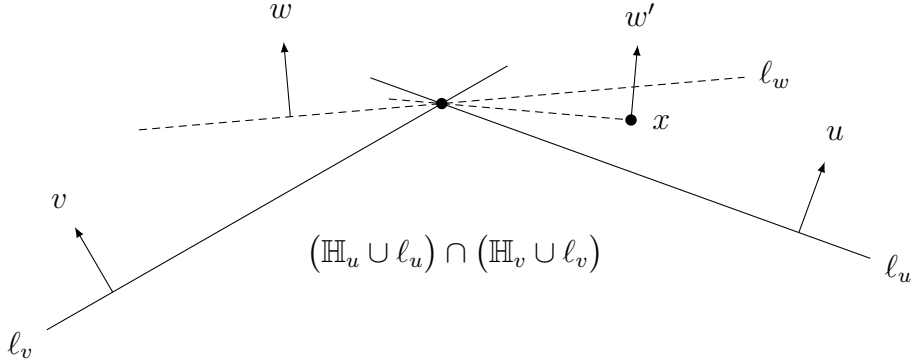

\begin{lemma}\label{le:quasi}
There exists a finite set $\qquasi \subset \Q_1$ such that for every pair $u,v$ of consecutive elements of $\stab \cup\qquasi$ there exists an update rule $X$ such that
\[
X \subset \big(\H_u\cup\ell_u\big)\cap\big(\H_v\cup\ell_v\big).
\]
Moreover, the same holds for any finite set of rational directions $\qquasi' \supset \qquasi$.
\end{lemma}

\begin{proof}
Form $\qquasi$ by taking the two unit vectors $u$ and $-u$ perpendicular to $x$ (considered as a vector) for every site $x\in X$ and every update rule $X\in\U$. Formally,
\[
\qquasi := \bigcup_{X\in\U} \bigcup_{x\in X} \big\{ u\in S^1 : \< u,x \> = 0 \big\}.
\]
Now let $\qquasi' \supset \qquasi$ be a finite set of rational directions, suppose that $u$ and $v$ are consecutive elements of $\stab \cup \qquasi'$ and let $w\in [u,v] \setminus \{u,v\}$. Since $w$ is not stable, there exists an update rule $X \subset \H_w$. Suppose the conclusion of the lemma fails, so
\[
X \not\subset \big( \H_u \cup \ell_u \big) \cap \big( \H_v \cup \ell_v \big).
\]
Then without loss of generality there exists $x \in X$ such that $\langle x,v\rangle < 0$ and $\langle x,u \rangle > 0$. But this implies that there exists $w' \in S^1$ perpendicular to $x$ with $w' \in [u,v] \setminus \{u,v\}$, contradicting the construction of $\qquasi$. (See Figure~\ref{fi:quasi}.)
\end{proof}

It follows immediately from the lemma that when droplets are taken with respect to finite supersets of $\stab\cup\qquasi$, there are rules that allow the droplets to grow along their sides all the way to the corners: \emph{droplets grow into droplets}.

\section{The upper bound for balanced families}\label{se:upperbal}

In this section we shall prove the following theorem, which is the upper bound of Theorem~\ref{th:main} for balanced families.

\begin{theorem}\label{th:upperbal}
Let $\U$ be critical and balanced. Then
\[
p_c(\Z_n^2,\U) = O\bigg(\frac{1}{\log n}\bigg)^{1/\alpha}.
\]
\end{theorem}

Recall that if $\U$ is balanced then there exists a closed semicircle $C \subset S^1$ such that $\alpha(u) \leq \alpha$ for all $u \in C$. Since $\alpha(u)<\infty$ for every $u\in C$, every stable direction $u\in C$ must be isolated, by Lemma~\ref{le:alphainfinity}. This implies (by Lemma~\ref{le:stableset}) the existence of a closed arc $C' = [u^-,u^+]$ such that $C \subsetneq C' \subsetneq S^1$ and such that $\alpha(u) \leq \alpha$ for all $u\in C'$. We may assume that the endpoints of $C'$ (the arc running anticlockwise from $u^-$ to $u^+$) are rational.

Let $\qquasi$ be the set of quasi-stable directions given by Lemma \ref{le:quasi} and set 
\[
\quasi := \big(\stab\cup\qquasi\cup\{u^-,u^+\}\big) \cap C' \qquad \text{and} \qquad \quasi' := \quasi \setminus \{u^-,u^+\}.
\]
These sets are finite, since $\qquasi$ and $\stab \cap C'$ are both finite by construction. Throughout this section droplets will be taken with respect to the set $\quasi$. 

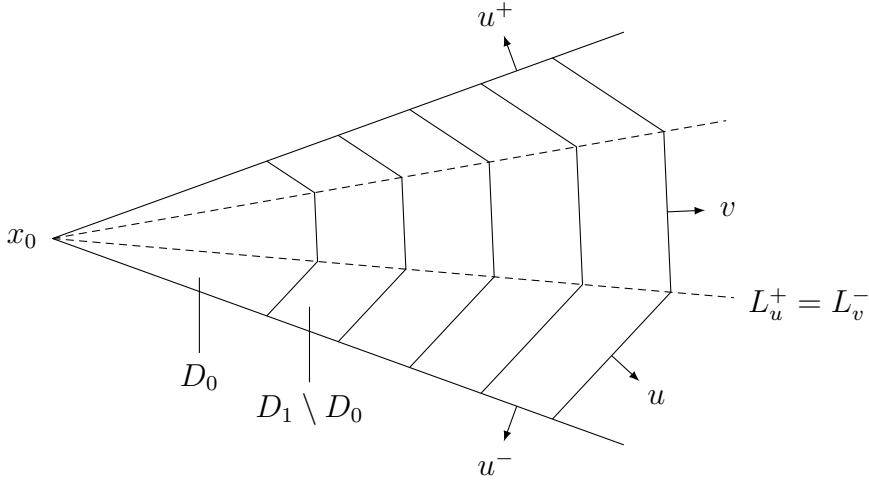
\begin{figure}[ht]
  \centering
  \begin{tikzpicture}[>=latex]
    \path [name path=C1] (10:7*7/6) circle (3);
    \path [name path=C2] (-5:7*7/6) circle (3);
    \path [name path=C3] (-20:7) circle (3);
    \pgfresetboundingbox
    \draw (0:0) -- (20:8);
    \draw [densely dashed] (0:0) -- (10:9);
    \draw [densely dashed] (0:0) -- (-5:9);
    \draw (0:0) -- (-20:8);
    \foreach \x in {3,4,...,7}
      \draw (20:\x) -- (10:\x*7/6) -- (-5:\x*7/6) -- (-20:\x);
    \draw (-15:2) -- ++(0,-1) node[below] {$D_0$};
    \draw (-15:3.5) -- ++(0,-1) node[below] {$D_1\setminus D_0$};
    \node at (-0.4,0) {$x_0$};
    \node at (-5:10) {$L_u^+=L_v^-$};
    \draw [->] (20:6.5) -- ++(110:0.5);
    \path (20:6.5) -- ++(110:0.8) node {$u^+$};
    \draw [->] (-20:6.5) -- ++(-110:0.5);
    \path (-20:6.5) -- ++(-110:0.8) node {$u^-$};
    \path [name intersections={of=C2 and C1,by=Q1}];
    \path [name intersections={of=C2 and C3,by=Q2}];
    \coordinate (P1) at ($(10:7*7/6) !0.5! (-5:7*7/6)$);
    \coordinate (P2) at ($(-5:7*7/6) !0.5! (-20:7)$);
    \draw [->] (P1) -- ($(P1)!0.5cm!(Q1)$);
    \draw [->] (P2) -- ($(P2)!-0.5cm!(Q2)$);
    \node at ($(P1)!0.8cm!(Q1)$) {$v$};
    \node at ($(P2)!-0.8cm!(Q2)$) {$u$};
  \end{tikzpicture}
  \caption{The sequence of droplets $D_0\subset D_1\subset D_2\subset\cdots$.}
  \label{fi:baldroplets}
\end{figure}

Choose a collection of vectors $\{a_u \in \Z^2 : u \in \quasi\}$ and sufficiently large positive constants $\{d_u > 0 : u \in\quasi' \}$ such that $d_u u \in \Z^2$, and the sequence of $\quasi$-droplets
\begin{equation}\label{eq:Dm}
D_m := \bigcap_{u\in\{u^-,u^+\}} \H_u(a_u) \cap \bigcap_{u\in\quasi'} \H_u(a_u+m d_u u)
\end{equation}
for $m = 0,1,2,\ldots$ have the following properties (see Figure~\ref{fi:baldroplets}):
\begin{enumerate}
\item $D_0$ is sufficiently large relative to the $d_u$;
\item for every $m \geq 0$ and every consecutive pair $u,v \in \quasi'$, the intersection\footnote{These are discrete lines and may have empty intersection, in which case we mean the intersection of the corresponding continuous lines. Alternatively, since both lines are rational we may ensure that they intersect in an element of $\Z^2$ by choosing a suitable dilation.} of the lines $\ell_u + a_u+m d_u u$ and $\ell_v + a_v + m d_v v$ lies on a (continuous) line $L_u^+ = L_v^-$;
\item the lines $L_u^+$ all intersect at the point $x_0\in\R^2$, which is also the intersection point of the sides of $D_0$ corresponding to $u^-$ and $u^+$.
\end{enumerate}

The key lemma in our proof of Theorem \ref{th:upperbal} will be the following bound on the probability that a droplet grows by a constant distance. We remark that the constants implicit in the statement depend on our choice of droplets, and hence on $\U$, but not on the probability $p$.

\begin{lemma}\label{le:upperkey}
Let $m\in\N$. Then
\[
\Pr_p\Big( D_m \subset \big[ D_{m-1} \cup (D_{m+1}\cap A) \big] \Big) \geq \big(1-(1-p^\alpha)^{\Omega(m)}\big)^{O(1)}.
\]
\end{lemma}

Having already established the deterministic lemmas of the previous section, the idea of the proof of Lemma~\ref{le:upperkey} is simple. In order to grow from $D_{m-1}$ to $D_m$ it is sufficient for a bounded number of events to occur, each event having failure probability at most $(1-p^\alpha)^{\Omega(m)}$. These events are all very loosely speaking of the form `there exists in $A$ a translate of a given set of $\alpha$ sites somewhere along one of the edges of the droplet'. Since any set of $\alpha$ sites is a subset of $A$ with probability $p^\alpha$, and since there are $\Omega(m)$ possible disjoint translates of that set, we obtain the desired bound on the probability. For completeness, we now give the details.\footnote{We also refer the reader to the proof of~\cite[Lemma~5.5]{BSUgen}, which is very similar.}

\begin{proof}[Proof of Lemma~\ref{le:upperkey}]
For each $u \in \quasi'$, fix a set $Z_u$ that is voracious for $u$, and recall that $|Z_u| \le \alpha$, and that the $u$-side of $D_{m-1}$ has length $\Omega(m)$. By Lemma~\ref{le:line}, applied with $\T = \quasi$ and $Z = Z_u$, it follows that there exist constants $\mu > 0$ and $r \in \N$ such that the following holds: there exist collections $\cZ_u^{(1)},\ldots,\cZ_u^{(r)}$ of disjoint translates of $Z_u$, each of size $\Omega(m)$, such that if at least one member of each collection is contained in $A$, then $\big[ D_{m-1} \cup (D_{m+1}\cap A) \big]$ contains all elements of $\partial_\circ(D_{m-1},u)$ that are at distance at least $\mu$ from the corners. Moreover, observe that the probability that $A$ contains at least one member of each collection is at least
$$\big( 1 - (1 - p^\alpha)^{\Omega(m)}\big)^{r}.$$
We can now repeat this step (on the $u$-side of the new infected set) a constant number of times, each time adding all of the next line along the $u$-side up to distance $O(1)$ from the corners (the distance from the corners increasing each time), until we reach the $u$-side of $D_m$ (recall that $D_0$ was chosen sufficiently large relative to the $d_u$). 

Repeating this entire process for each $u \in \quasi'$, it follows that with probability at least 
$$\big(1 - (1 - p^\alpha)^{\Omega(m)}\big)^{O(1)}$$
we are able to infect all of $D_m$ except for $O(1)$ sites near each corner. Finally, by Lemma~\ref{le:quasi}, the remaining elements of $D_m$ are infected (deterministically), one line at a time (working towards the corners), in order of increasing distance from~$D_{m-1}$. 
\end{proof}

The next step is to deduce from Lemma~\ref{le:upperkey} the following lower bound on the probability that a large droplet is almost internally filled.

\begin{lemma}\label{le:upperweak}
Let $m\in\N$. Then
\[
\Pr_p\Big( D_m \subset \big[ D_{m+1}\cap A \big] \Big) \geq \exp\Big(-O\big(p^{-\alpha}\big)\Big).
\]
\end{lemma}

\begin{proof}
Noting that all the events we are considering are increasing, it follows from Harris's inequality (Lemma~\ref{le:harris}) that
\begin{align*}
\Pr_p\Big( D_m \subset \big[ D_{m+1}\cap A \big] \Big) & \geq \Pr_p\bigg( I(D_0) \cap \bigcap_{k=1}^m \Big\{ D_k \subset \big[ D_{k-1} \cup (D_{k+1}\cap A) \big] \Big\} \bigg) \\
&\geq \Pr_p\big(I(D_0)\big) \prod_{k=1}^m \Pr_p\Big( D_k \subset \big[ D_{k-1} \cup (D_{k+1}\cap A) \big] \Big).
\end{align*}
Thus, by Lemma~\ref{le:upperkey}, we have
\begin{align*}
\Pr_p\Big( D_m \subset \big[ D_{m+1}\cap A \big] \Big) & \geq  p^{O(1)} \prod_{k=1}^\infty \left(1-(1-p^\alpha)^{\Omega(k)}\right)^{O(1)} \\
&\geq p^{O(1)} \exp\bigg(-O(1)\sum_{k=1}^\infty -\log\Big(1-e^{-\Omega(p^\alpha k)}\Big) \bigg) \\
&\geq p^{O(1)} \exp\left(-O\big(p^{-\alpha}\big) \int_0^\infty -\log\big(1-e^{-z}\big)\,dz \right) \\
&\geq \exp\Big(-O\big(p^{-\alpha}\big)\Big),
\end{align*}
where for the final inequality we used the fact that $\int_0^\infty -\log\big(1-e^{-z}\big)\,dz<\infty$.
\end{proof}

From here, the deduction of Theorem~\ref{th:upperbal} is straightforward.

\begin{proof}[Proof of Theorem~\ref{th:upperbal}]
Let $\lambda$ be a sufficiently large constant, and set
\[
p = \bigg( \frac{\lambda}{\log n} \bigg)^{1/\alpha}.
\]
As usual, $A$ is a $p$-random subset of $\Z_n^2$. We shall show that $[A] = \Z_n^2$ with high probability as $n \to \infty$, which is more than enough to prove the theorem.

To avoid some technical issues relating to wrap-around, instead of working in $\Z_n^2$ we will consider the process in a $C n \times Cn$ square $S \subset \Z^2$, where $C = C(\U)$ is a sufficiently large (integer) constant so that there exists $m \in \N$ and an $n \times n$ square $T \subset S$ with the following property: 
$$[n]^2 \subset x + D_m \qquad \text{and} \qquad x + D_{m+1} \subset S$$
for any $x \in T$. It is easy to see that a $p' := p/C^2$-random subset $A' \subset S$ is stochastically dominated by a $p$-random subset $A \subset \Z_n^2$ (simply project to the torus), so it will suffice to prove the claim in $S$, with a (different) sufficiently large constant $\lambda$.\footnote{We are grateful to the referee for suggesting this simple and elegant way of dealing with the `wrap-around' on the torus.}

To prove the result in $S$, we `sprinkle' the initially infected sites in two rounds; that is, we take $A^{(1)}$ and $A^{(2)}$ to be independent $p$-random subsets of $S$, and redefine the set of initially infected sites to be $A' = A^{(1)} \cup A^{(2)}$ (losing another factor of $2$). We use the first round of sprinkling to find an almost internally filled droplet $x + D_k$ for some $x \in T$, where $k := (\log n)^3$, and the second round to show that this copy of $D_k$ grows (with high probability) to fill $x + D_m$, and hence $[n]^2$.

Let us define the following events:
\[
E := \bigcup_{x \in T} \Big\{x + D_k \subset \big[(x+D_{k+1})\cap A^{(1)}\big]\Big\}
\]
is the event that $x+D_k$ is `almost internally filled' by $A^{(1)}$ for some $x \in T$, and 
\[
F(x) := \Big\{ x + D_m \subset \big[ (x + D_k) \cup A^{(2)} \big] \Big\}
\]
is the event that $x + D_k$ grows (with the help of $A^{(2)}$) to fill $x + D_m$, and hence $[n]^2$. The events $E$ and $F(x)$ are independent, so by the comments above, it will suffice to show that $\Pr_p(E^c) = o(1)$ and that $\Pr_p\big(F(x)^c\big) = o(1)$ for each fixed $x \in T$. 

To bound $\Pr_p(E^c)$, observe first that there exists a collection of $\Omega\big(n^2/k^2\big)$ sites $x \in T$ such that the sets $x + D_{k+1}$ are pairwise disjoint. By Lemma~\ref{le:upperweak}, and since $\lambda$ is sufficiently large, it follows that
\[
\Pr_p\big(E^c\big) \leq \Big(1 - \exp\big( - O(p^{-\alpha}) \big) \Big)^{\Omega(n^2/k^2)} \leq \exp\Big(- n^{2 + o(1)} e^{ - O(\log n) / \lambda} \Big) = o(1).
\] 

To bound $\Pr_p\big( F(x)^c \big)$, recall that $x + D_{m+1} \subset S$ for each $x \in T$. By Lemma~\ref{le:upperkey} and Harris's inequality, it follows that, for any $x \in T$,
\begin{align*}
\Pr_p\big( F(x) \big) & \geq \prod_{\ell = k}^m \Big( 1 - \big( 1 - p^\alpha \big)^{\Omega(\ell)} \Big)^{O(1)} \geq \Big( 1 - \big(1 - p^\alpha \big)^{\Omega(k)} \Big)^{O(n)} \\
& \geq \exp\Big( -O\big( e^{-\Omega(p^\alpha k)}\cdot n \big) \Big) = 1 - o(1),
\end{align*}
as claimed. By the comments above, this completes the proof of the theorem.
\end{proof}

\section{The upper bound for unbalanced families}\label{se:upperunbal}

In this final section on upper bounds we prove the following general theorem, which in particular implies the upper bound in Theorem~\ref{th:main} for unbalanced families.

\begin{theorem}\label{th:upperunbal}
Let $\U$ be a critical update family. Then
\[
p_c(\Z_n^2,\U) = O\bigg(\frac{(\log\log n)^2}{\log n}\bigg)^{1/\alpha}.
\]
\end{theorem}

The theorem does not require the hypothesis that $\U$ is unbalanced, although of course it is only under this assumption that the result is tight up to the implicit constant. It may be helpful in this section to think of $\U$ as being unbalanced, even though this is not strictly necessary.

By the definition of $\alpha$ (and Lemma~\ref{le:stableset}), there exists an open semicircle $C\subset S^1$, with rational endpoints, such that $\alpha(u)\leq\alpha$ for all $u \in C$. Let $w^\perp$ be the midpoint of $C$, let $u^*$ and $-u^*$ be the left and right endpoints of $C$ respectively, and observe that $\alpha^+(u^*) < \infty$, 
by Lemma~\ref{le:alphainfinity}. Thus, by Lemma~\ref{le:leftrightcluster}, there exists $X \in \U$ such that $X \subset \H_{u^*} \cup \ell_{u^*}^-$, and hence there exists a finite set of consecutive sites $Z \subset \ell_{u^*}$ such that $\ell_{u^*}^+ \subset [ \H_{u^*} \cup Z ]$. Define $\alpha^*$ to be the (minimum) size of such a set $Z$.

Let $\qquasi$ be the set of quasi-stable directions given by Lemma~\ref{le:quasi}, and set
\[
\quasi := \big((\stab\cup\qquasi)\cap C\big) \cup \{u^*,-u^*,-w^\perp\} \quad \text{and} \quad \quasi' := (\stab\cup\qquasi)\cap C.
\]
As in the previous section, both of these sets are finite. In this section all droplets will be $\quasi$-droplets. Since the growth process will predominantly take place in directions parallel to the vectors $w^\perp$ and $u^*$, to simplify the notation we rotate the lattice $\Z^2$ so that $u^*$ is directed vertically upwards.
The discrete rectangle with opposite corners $(a,b)$ and $(c,d)$ is thus defined to be
\[
R\big((a,b),(c,d)\big) := \Big\{ xw^\perp + yu^* \in \Z^2 \,:\, x,y\in\R,\; a \leq x \leq c \text{ and } b \leq y \leq d \Big\}.
\] 
The sequences of droplets will be defined in terms of the following quantities:
\[
m_1(p) := \frac{\lambda_1}{p^\alpha}\log\frac{1}{p}, \quad m_2(p) := p^{-\lambda_2}, \quad m_3(p) = p^{2\alpha^*}m_2(p) \quad \text{and} \quad m_4(p) := \lambda_1 n,
\]
where $\lambda_1 \gg \lambda_2 \gg 1$ are sufficiently large positive constants and $n=n(p)$, to be specified later (see~\eqref{eq:tdef:unbalanced}), satisfies
\begin{equation}\label{eq:unbal:logn:cond}
\log n \leq p^{-\lambda_2/2}.
\end{equation}

\begin{figure}[ht]
  \centering
  \begin{tikzpicture}[>=latex]
    \draw (0,0) rectangle (0.5,3) (0.5,0) rectangle (8,3) (8,0) rectangle (8.5,6);
    \draw (11,0) -- (8,0) -- (8,6) -- (11,6);
    \draw [dashed] (11,0) -- (12,0) (11,6) -- (12,6) (8,3) -- (9.5,3);
    \draw (0,3) -- (8.505,6.195);  \draw (8.5,6) -- (8.5,6.195);
    \begin{scope}
      \clip (0.5,0) rectangle (12,3);
      \foreach \x in {0,0.5,...,2}
        \draw (\x,3) -- (\x+0.8,2) -- (\x+0.9,1) -- (\x,0);
    \end{scope}
    \begin{scope}
      \clip (8.5,0) rectangle (12,6);
      \foreach \x in {7,7.5,...,9}
        \draw (\x,6) -- (\x+1.6,4) -- (\x+1.8,2) -- (\x,0);
    \end{scope}
    \draw [->] (5,1.7) -- node [below] {1} (6,1.7);
    \draw [->] (5.5,3.5) -- node [right] {2} (5.5,4.5);
    \draw [->] (11,3) -- node [below] {3} (12,3);
    \draw [->] (1,5.5) -- (1,6.3) node [above] {$u^*$};
    \draw [->] (1,5.5) -- (1.8,5.5) node [right] {$w^\perp$};
    \draw (0.25,2.5) -- (0.25,3.8) node [above] {$R_0$};
    \node at (3.8,1.5) {$R_1$};
    \node at (11.5,1.2) {$R_3$};
    \draw (8.25,0.6) -- (7.75,0.6) node [left] {$R_2$};
    \node at (3.8,3.7) {$T$};
    \draw [<->] (0,-0.4) -- node [below] {$\lambda_1$} (0.5,-0.4);
    \draw [<->] (0,-0.2) -- node [below] {$m_2(p)$} (8,-0.2);
    \draw [<->] (-0.2,0) -- node [left] {$m_1(p)$} (-0.2,3);
    \draw [<->] (7.8,3) -- node [left] {$m_3(p)$} (7.8,6);
    \draw (1.1,2.25) -- (1.1,3.8) node [above] {$D_0^{(1)}$};
    \draw (9.1,5.25) -- (10,5.25) node [right] {$D_0^{(3)}$};
  \end{tikzpicture}
  \caption{The growth mechanism in the unbalanced setting.}\label{fi:unbaldroplets}
\end{figure}
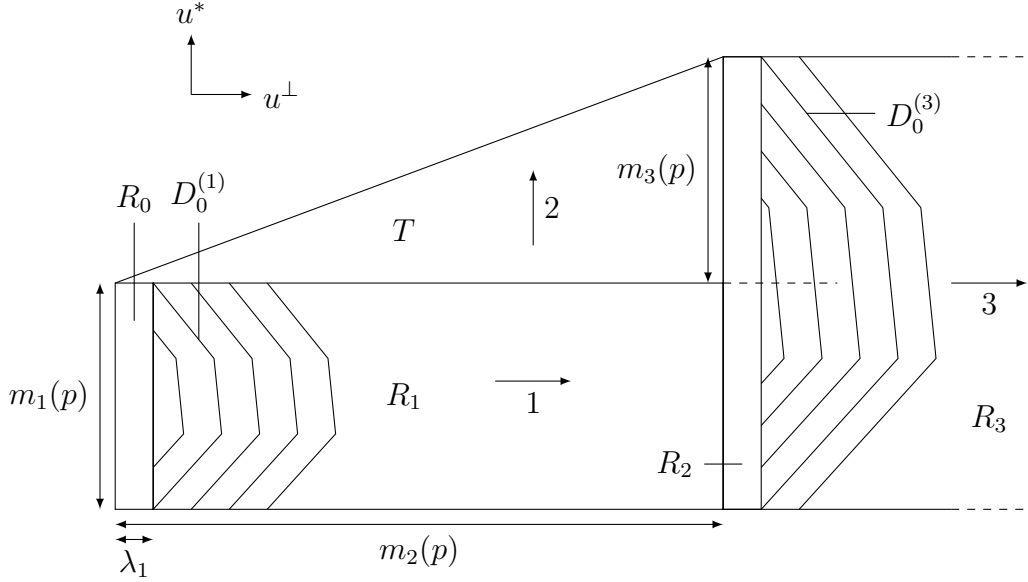

Let
\begin{gather*}
R_0 := R\Big((0,0),\big(\lambda_1,m_1(p)\big)\Big), \quad R_1 := R\Big((0,0),\big(2m_2(p),m_1(p)\big)\Big), \\
R_2 := R\Big(\big(m_2(p),0\big),\big(m_2(p)+\lambda_1,m_1(p)+m_3(p)\big)\Big), \\
\text{and} \quad R_3 := R\Big(\big(m_2(p),0\big),\big(m_4(p),m_1(p)+m_3(p)\big)\Big)
\end{gather*}
be rectangles, and let
\[
T := \Big\{ xw^\perp + y u^* \in\Z^2 \,: \, 0 \leq x \leq m_2(p) + \lambda_1 \text{ and } 0 \leq y-m_1(p) \leq p^{2\alpha^*} x \Big\}
\]
be a triangle; see Figure~\ref{fi:unbaldroplets}.

For technical reasons, we also need to use the rectangles
\begin{gather*}
R_1' := R\Big((0,0),\big(4m_2(p),m_1(p)\big)\Big) \\
\text{and} \quad R_3' := R\Big(\big(m_2(p),0\big),\big(2m_4(p),m_1(p)+m_3(p)\big)\Big).
\end{gather*}
which are roughly twice as long as $R_1$ and $R_3$ respectively.

Figure~\ref{fi:unbaldroplets} illustrates the growth mechanism we use to prove Theorem~\ref{th:upperunbal}. It comes in five stages, and, as in the previous section, we use sprinkling to maintain independence between the different stages.
\begin{itemize}
\item \emph{Stage $0$.} We find a copy of $R_0$ contained in $A$.
\item \emph{Stage $1$.} The infection spreads in the direction $w^\perp$ from $R_0$ and fills the rectangle $R_1$. This occurs in a similar way to growth in balanced models, except that the rows are not increasing in size.
\item \emph{Stage $2$.} The infection spreads in the direction $u^*$ from $R_1$, using infected sites in the triangle~$T$ to fill $R_2$.
\item \emph{Stage $3$.} Exactly as in Stage 1, the infection spreads in the direction $w^\perp$ from $R_2$ to fill $R_3$.
\item \emph{Stage $4$.} $R_3$ is a `strip' that wraps around the torus and returns to its starting point (since $u^*$ and $w^\perp$ are rational directions). The infection now spreads in direction $u^*$ from $R_3$ (like in Stage 2) to infect the rest of $\Z_n^2$.
\end{itemize}

Our task is to make the above sketch precise. We postpone the proof of the following key lemma until later in the section.

\begin{lemma}\label{le:upperunbalkey}
The event
\[
\Big\{ R_3\subset \big[ R_0\cup\big((R_1' \cup T \cup R_3')\cap A\big)\big] \Big\}
\]
occurs with high probability as $p\to 0$.
\end{lemma}

From here, the deduction of Theorem~\ref{th:upperunbal} is relatively straightforward, the main complication being the growth (in Stage 4) from $R_3$ to $\Z_n^2$. 

\begin{proof}[Proof of Theorem~\ref{th:upperunbal}]
The proof is similar to the proof of Theorem~\ref{th:upperbal} for balanced families. Let $\lambda \gg \lambda_1$ be sufficiently large, and set
\begin{equation}\label{eq:tdef:unbalanced}
p = \bigg(\frac{\lambda(\log\log n)^2}{\log n}\bigg)^{1/\alpha}.
\end{equation}
We shall show that $[A] = \Z_n^2$ with high probability as $n \to \infty$.

As before, we sprinkle in two rounds, each round using probability $p$ (which, also as before, is permissible, if a slight abuse of notation), and denote by $A^{(1)}$ and $A^{(2)}$ the sites infected in each round. There are (crudely) at least $n$ choices of $x \in \Z_n^2$ such that the sets $x + R_0$ are disjoint, and the probability that $x + R_0 \not\subset A^{(1)}$ for all such $x$ is at most
\begin{equation}\label{eq:unbal:R0calc}
\big(1-p^{O(m_1(p))}\big)^n.
\end{equation}
Noting that
\[
p^{O(m_1(p))} = \exp\Bigg( - \frac{O(1)}{p^\alpha}\bigg(\log\frac{1}{p}\bigg)^2 \Bigg) \geq \frac{1}{\sqrt{n}}
\]
since $\lambda$ is sufficiently large, it follows that~\eqref{eq:unbal:R0calc} tends to 0 as $n \to \infty$.

Now fix $x$ such that $x+R_0 \subset A^{(1)}$, and in fact without loss of generality let us assume that $x = \0$. By Lemma~\ref{le:upperunbalkey}, 
\[
\Pr_p\Big( R_3\subset \big[ R_0\cup\big((R_1' \cup T \cup R_3')\cap A^{(2)}\big)\big] \Big) = 1-o(1),
\]
since the condition on $n$ in~\eqref{eq:unbal:logn:cond} holds with our definition of $p$. Now, since $u^*$ and $w^\perp$ are rational directions and $\lambda_1$ is sufficiently large, it follows that $R_3$ is a `strip' that wraps around $\Z_n^2$ a bounded number of times before returning to its starting point. Recall that $\alpha^*$ consecutive sites of $\ell_{u^*}$, together with $\H_{u^*}$, are sufficient to infect all sites of $\ell_{u^*}$ to the right. Thus, in order to infect the remaining sites in $\Z_n^2$, it is enough for the following event to occur: every line in $\Z_n^2$ parallel to $w^\perp$ contains $\alpha^*$ consecutive sites. Indeed, since the height of $R_3$ is larger than $\nu$, this would ensure that the remaining lines above the strip $R_3$ are infected one-by-one. Since each line has length at least $\Omega(n)$, and there are at most $O(n)$ lines, the probability that this event fails is at most
\[
O(n) \cdot \big(1-p^{\alpha^*}\big)^{\Omega(n)} = o(1),
\]
and this completes the proof of the theorem.
\end{proof}

We have reduced our task to that of proving Lemma~\ref{le:upperunbalkey}. As in the previous section, given the framework of voracity and quasi-stability from Section~\ref{se:upperprelims}, the idea of the proof is simple. In Stage 1 of the process, the probability of advancing a constant number of steps is
\[
\big(1-(1-p^\alpha)^{\Omega(m_1(p))}\big)^{O(1)} \leq (1-p^{\Omega(\lambda_1)})^{O(1)}.
\]
Since $\lambda_1 \gg \lambda_2$, the set should grow to fill $R_1$, and for similar reasons, the infection spreads out rightwards from $R_2$ to fill $R_3$. Both of these steps are almost the same as the corresponding part of the proof for balanced models. The growth upwards from $R_1$ through $T$ to fill $R_2$ is a little different. Since the infection might only spread rightwards when advancing row-by-row in the $u^*$ direction, each set of $\alpha^*$ consecutive infected sites we find when growing upwards through $T$ from $R_1$ must lie to the right of the previous set. Nevertheless, the probability of filling $T$ (except possibly for a small number of sites near the diagonal) is at least
\[
\big( 1 - (1 - p^{\alpha^*})^{\Omega(p^{-2\alpha^*})} \big)^{O(m_3(p))} = 1-o(1).
\]
We now give the details, which are similar to those in Section~\ref{se:upperbal}.
 

We define two sequences of droplets as in~\eqref{eq:Dm}, except with $u^+ = u^*$ and $u^- = -u^*$ (so the corresponding lines are now parallel), and with $-w^\perp$ added to the set of quasi-stable directions. Specifically, for each $m \in \Z$ and each $i \in \{1,3\}$, define the $\quasi$-droplets
\begin{equation}\label{eq:Fm}
D^{(i)}_m := R_i \cap \bigcap_{u \in\quasi'} \H_u\big(a_u + m d_u u \big)
\end{equation}
for some $\{a_u \in \Z^2 : u \in \quasi'\}$ and sufficiently large constants $\{d_u > 0 : u \in\quasi' \}$ such that $d_u u \in \Z^2$, and the following hold (see Figure~\ref{fi:unbaldroplets}): 
\begin{itemize}
\item $R_{i-1} \subset D^{(i)}_0$;
\item for every consecutive pair $u,v \in \quasi'$, there exists a horizontal line $L_u^+ = L_v^-$ (that is, one parallel to $w^\perp$) that intersects $R_i$, such that for every $m \in \Z$, the intersection of $\ell_u + a_u+m d_u u$ and $\ell_v + a_v + m d_v v$ lies on $L_u^+ = L_v^-$;
\item for each $u \in \quasi'$ and each $m\in\N$, the $u$-side of $D_m^{(i)}$ has size $\Omega\big(m_i(p)\big)$.
\end{itemize}
Note that we shall also need to use $D^{(i)}_m$ for those negative values of $m$ for which the droplet is non-empty, as well as for positive values of $m$.

The following lemma is essentially Lemma~\ref{le:upperkey} applied to the droplets $D^{(i)}_m$, and so the proof is omitted. 

\begin{lemma}\label{le:upperunbalstep}  
Let $i \in \{1,3\}$ and $m \in \Z$. Then
\begin{equation}
\Pr_p\Big( D^{(i)}_m \subset \big[ R_{i-1} \cup D^{(i)}_{m-1} \cup ( D^{(i)}_{m+1} \cap A ) \big] \Big) \geq \big(1-(1-p^\alpha)^{\Omega(m_i(p))}\big)^{O(1)}. \tag*{\qedsymbol}
\end{equation}
\end{lemma}

We now complete the proof of Lemma~\ref{le:upperunbalkey}. Note that, since the rectangle $R_3'$ only wraps around $\Z_n^2$ a bounded number of times, we may work entirely in $\Z^2$.

\begin{proof}[Proof of Lemma~\ref{le:upperunbalkey}]
We shall show that the event
\[
\Big\{ R_1 \subset \big[ R_0 \cup (R_1' \cap A)\big] \Big\} \cap \Big\{ R_2 \subset \big[ R_1 \cup (T \cap A) \big] \Big\} \cap \Big\{ R_3 \subset \big[ R_2 \cup (R_3' \cap A) \big] \Big\}
\]
occurs with high probability as $p\to 0$, which clearly implies the lemma. We begin by deducing from Lemma~\ref{le:upperunbalstep} that the first and third parts of this event occur with high probability as $p \to 0$. To see this, let $i \in \{1,3\}$ and observe that there exists $m \in \N$ such that $R_i \subset D_m^{(i)} \subset D_{m+1}^{(i)} \subset R_i'$, if $p$ is sufficiently small, where $m = O\big(m_{i+1}(p)\big)$. It follows from Lemma~\ref{le:upperunbalstep} that
\begin{multline*}
\Pr_p\Big( D^{(i)}_m \subset \big[ R_{i-1} \cup \big( D^{(i)}_{m+1} \cap A \big) \big] \Big) \geq \Big( 1 - (1 - p^\alpha )^{\Omega(m_i(p))} \Big)^{O(m_{i+1}(p))}\\
\geq \exp\Big( - O\big(m_{i+1}(p)\big) \cdot \exp\Big( - \Omega\big( m_i(p) \cdot p^\alpha \big) \Big) \Big) = 1 - o(1)
\end{multline*}
as $p \to 0$. Indeed, we have $\exp\big( -\Omega(m_1(p) \cdot p^\alpha) \big) = p^{\Omega(\lambda_1)} = o\big(1/m_2(p)\big)$, and 
\[
\exp\Big( - \Omega\big( m_3(p) \cdot p^\alpha \big) \Big) = \exp\Big( - \Omega\big( p^{-\lambda_2+2\alpha^*+\alpha} \big) \Big) < \frac{1}{n^2} = o\bigg(\frac{1}{m_4(p)}\bigg),
\]
where we used our assumptions that $\lambda_1 \gg \lambda_2 \gg 1$ and $\log n \leq p^{-\lambda_2/2}$.

It remains to show that the event
\[
\Big\{ R_2 \subset \big[ R_1 \cup (T \cap A) \big] \Big\}
\]
occurs with high probability as $p \to 0$. To do so, consider the set $U_i$ of the leftmost $p^{-2\alpha^*}$ sites of $T \cap \ell_{u^*}(i)$ for each line $\ell_{u^*}(i)$ that intersects $T$. Now, suppose that, for every such line, the middle $p^{-2\alpha^*} / 3$ sites of $U_i$ contain a set of $\alpha^*$ consecutive sites of $A$. Then $R_2 \subset \big[ R_1 \cup (T \cap A) \big]$, by the definition of $\alpha^*$, and since the height of $R_1$ is larger than $\nu$. But this has probability at least 
\[
\big(1-(1-p^{\alpha^*})^{p^{-2\alpha^*}/3\alpha^*}\big)^{m_3(p)} \geq \exp\Big( - p^{- \lambda_2} \exp\big( - p^{ - \alpha^* / 2}) \big) \Big) = 1-o(1),
\]
as required.
\end{proof}

\section{Approximately internally filled sets}\label{se:approxfill}

In this section we shall lay the groundwork for the proofs of the lower bounds in Theorem~\ref{th:main} by defining and proving basic properties of three of our key tools: the \emph{covering}, \emph{spanning} and \emph{iceberg} algorithms. These should all be thought of as ways of using droplets to approximate the closure of $A$ under the $\U$-bootstrap process.

The covering algorithm, which we introduce in Section~\ref{se:covering}, replaces the rectangles process in the balanced case, and allows us to find in $\Z_n^2$ (if $A$ percolates) a droplet of size about $\log n$ containing $\Omega(\log n)$ disjoint, strongly connected subsets of $A$ of size~$\alpha$. For unbalanced models, we use the spanning algorithm, introduced in Section~\ref{se:spanning}, to find an internally spanned critical droplet and to construct an iterated sequence of `hierarchies' for this droplet. For models with drift (that is, unbalanced models that satisfy~\eqref{eq:unbal:driftcase}), we will in addition require the iceberg algorithm, which we will introduce in Section~\ref{se:iceberg}, in order to bound the range of the $\U$-bootstrap process in certain directions with the help of half-planes. 

Having completed the proofs of the upper bounds of Theorem~\ref{th:main}, we no longer have any need for quasi-stable directions. In fact, henceforth all droplets will be assumed to be taken with respect to one of two specific finite sets of stable directions. 
The existence of these sets is verified in the next two lemmas. For $u\in \Q_1$, let
\begin{equation}\label{eq:alphabar}
\alphabar(u) := \min \big\{ \alpha^+(u),\alpha^-(u) \big\},
\end{equation}
so $\alphabar(u)=\alpha(u)$ if and only if $\alpha^+(u)$ and $\alpha^-(u)$ are either both finite or both infinite. 

\begin{lemma}\label{le:stabb}
If $\,\U$ is a critical update family, then there exists a finite set $\stabb \subset \Q_1$ such that:
\begin{enumerate}
\item $\alphabar(u) \geq \alpha$ for every $u\in\stabb$; and\smallskip
\item $\stabb \cap C \neq \emptyset$ for every open semicircle $C \subset S^1$.
\end{enumerate}
\end{lemma}


We will be most interested in $\stabb$ when dealing with balanced update families. We remark that condition~$(ii)$ is equivalent to the origin lying in the interior of the convex hull of $\stabb$, and also to $\stabb$-droplets being finite.

\begin{proof}[Proof of Lemma~\ref{le:stabb}]
Observe first that, by~\eqref{eq:alphaU} in Definition~\ref{de:alpha}, there exists a finite set $\T \subset \Q_1$, satisfying condition~$(ii)$, such that $\alpha(u) \geq \alpha$ for all $u \in \T$. Now, recall from Lemma~\ref{le:alphainfinity} that $\alpha(u) = \alphabar(u)$ unless $u$ is an endpoint of a non-trivial interval of $\stab$. However, if this is the case for some $u \in \T$, then there exist vectors $u' \in \Q_1$ with $\alphabar(u') = \infty$ arbitrarily close to $u$. Choosing such a $u'$ sufficiently close to $u$, and replacing $u$ by $u'$, we see that condition~$(ii)$ still holds. Repeating this for each $u \in \T$ with $\alpha(u) \ne \alphabar(u)$, we obtain a set with the desired properties.
\end{proof}

It is easy to see that we may in fact take $\stabb$ to have size 3, except when $|\stab|=4$ and $\stab = -\stab$ (i.e., $\stab$-droplets are parallelograms), in which case we may take $\stabb$ to have size 4. We shall not need this observation, however. 

We next show that a suitable collection of stable directions exists when $\U$ is unbalanced; the properties we need in this case are somewhat different. 

\begin{lemma}\label{le:stabu}
If $\,\U$ is an unbalanced critical update family, then there exists a finite set $\stabu =\{u^*,-u^*,u^l,u^r\} \subset \Q_1$ such that the following hold:
\begin{enumerate}
\item $u^l$ lies in the open semicircle to the left of $u^*$, and $u^r$ lies in the open semicircle to the right of $u^*$;\smallskip
\item $\min\big\{ \alpha(u^*),\alpha(-u^*) \big\} \geq \alpha+1$;\smallskip
\item $\min\big\{ \alphabar(u^l),\alphabar(u^r) \big\} \geq \alpha$.
\end{enumerate}
\end{lemma}

\begin{proof}
Choose $u^*$ satisfying condition~$(ii)$ using Lemma~\ref{le:ustar}, and then $u^l$ and $u^r$ satisfying the remaining two conditions using Definition~\ref{de:alpha}. In particular, note that if one of the open semicircles bounded by $u^*$ and $-u^*$ contains an interval of stable directions then we may choose any interior point of this interval, and if not, then $\alpha(u) = \alphabar(u)$ for every $u$ in the open semicircle. 
\end{proof}

As mentioned before Lemma~\ref{le:stabb}, we shall henceforth fix, for each critical update family, a set $\stabb$, 
and if $\,\U$ is unbalanced also a set $\stabu$, with the above properties. We also make the following definition, which we will use extensively. 

\begin{definition}\label{def:DK}
Let $\T \subset \Q_1$ be a finite set of rational directions that intersects each open semicircle of $S^1$. Given 
a finite set $K \subset \Z^2$, we will write $D_\T(K)$ for the unique minimal $\T$-droplet\footnote{Note that this is obtained by taking a tangent line to $K$ in each direction of $\T$.} containing $K$. We will also write $D_B(K) = D_{\stabb}(K)$, and also $D_U(K) = D_{\stabu}(K)$ when $\,\U$ is unbalanced.
\end{definition}

We will frequently use the following simple but key `subadditivity' lemma. 

\begin{lemma}\label{lem:subadditivity}
Let $\T \subset \Q_1$ be a finite set of rational directions that intersects each open semicircle of $S^1$. For any pair of intersecting $\T$-droplets $D_1$ and $D_2$,
\begin{equation}\label{eq:subadd:diam}
\diam\big( D_\T(D_1\cup D_2) \big) \leq \diam(D_1) + \diam(D_2) + O(1).
\end{equation}
Moreover, 
\begin{equation}\label{eq:subadd:pi}
\pi\big(D_\T(D_1\cup D_2),u\big) \leq \pi(D_1,u) + \pi(D_2,u) + O(1)
\end{equation}
for every $u \in S^1$.
\end{lemma}

\begin{proof}
For $i = 1,2$, let $\hat{D}_i$ be the smallest continuous $\T$-droplet containing $D_i$, so 
$$\hat{D}_i = \bigcap_{u \in \T} \big\{ x \in \R^2 : \< x,u \> \le a_u(i) \big\}$$
with each $a_u(i) \in \R$ minimal subject to the condition $\< x,u \> \le a_u(i)$ for all $x \in D_i$. Let $\hat{D}$ be the smallest continuous $\T$-droplet containing $\hat{D}_1 \cup \hat{D}_2$, and note that
$$\hat{D} = \bigcap_{u \in \T} \big\{ x \in \R^2 : \< x,u \> \le \max\{ a_u(1), a_u(2) \} \big\}.$$
We claim that if $\0 \in \hat{D}_1 \cap \hat{D}_2$, then 
\begin{equation}\label{eq:Minkowski}
\hat{D} \subset \hat{D}_1 + \hat{D}_2,
\end{equation}
where $\hat{D}_1 + \hat{D}_2$ is the Minkowski sum of $\hat{D}_1$ and $\hat{D}_2$. 


To prove~\eqref{eq:Minkowski}, observe first that $\hat{D}_1 + \hat{D}_2$ is convex (since $\hat{D}_1$ and $\hat{D}_2$ are convex), and is therefore a continuous $\T$-droplet, since each of its sides is either the sum of a side and a vertex of $\hat{D}_1$ and $\hat{D}_2$ (in some order), or the sum of two parallel sides (note that the sum of two non-parallel sides would not be one-dimensional). Now simply note that, for each $u \in \T$, the $u$-side of $\hat{D}_1 + \hat{D}_2$ (if it has one) is at least as far from the origin as the $u$-sides of $\hat{D}_1$ and $\hat{D}_2$, since $\0 \in \hat{D}_1 \cap \hat{D}_2$. 

In order to deduce~\eqref{eq:subadd:pi}, observe that $\pi(X + Y,u) = \pi(X,u) + \pi(Y,u)$ for any (compact) $X,Y \subset \R^2$. Now, applying~\eqref{eq:Minkowski}, it follows that
\begin{align*}
\pi\big(D_\T(D_1\cup D_2),u\big) & \, \le \pi(\hat{D},u) \le \pi\big( \hat{D}_1 + \hat{D}_2,u\big) \\
& \, = \pi(\hat{D}_1,u) + \pi(\hat{D}_2,u) \leq \pi(D_1,u) + \pi(D_2,u) + O(1).
\end{align*}
Applying~\eqref{eq:subadd:pi} with $u$ such that $\diam\big( D_\T(D_1\cup D_2) \big) = \pi\big(D_\T(D_1\cup D_2),u\big)$, we also obtain~\eqref{eq:subadd:diam}. 
\end{proof}

Recall that in Section~\ref{se:defs} we defined $\nu = \nu(\U)$ to be the diameter of $\U$,
\[
\nu = \max \Big\{ \|x - y\| \,:\, x,y \in X \cup \{\0\}, \, X \in \U \Big\}.
\]
We will also need the following additional measure of the range of an update family:
\begin{equation}\label{eq:rho}
\rho = \rho(\U) := \sup \Big\{ \| y - Z \| \,:\, |Z| = \alpha - 1,\; y \in \big[\H_u\cup Z\big]\setminus\H_u,\; u\in\stabb \Big\},
\end{equation}
where the supremum is taken over all choices of $u$, $y$ and $Z$ satisfying the stated conditions. In order to prove that $\rho$ is finite for every 
update family $\U$, we will need the following extremal lemma from~\cite{BSUgen}.

\begin{lemma}[Lemma~4.7 of~\cite{BSUgen}]\label{lem:BSU:extremal}
For any finite set $Z \subset \Z^2$, the closure $[Z]$ is contained in a collection of disjoint $\S_B$-droplets\footnote{Here $\S_B$ could be replaced by any finite set $\T \subset \S$ such that $\T$-droplets are finite.}, each of diameter $O(|Z|)$. 
\end{lemma}

To prove Lemma~\ref{lem:BSU:extremal}, simply place an $\S_B$-droplet on each element of $Z$, and then recursively unite any pair that lie within distance $\nu$ of each other by replacing them by the smallest $\S_B$-droplet that contains both (cf.~Definition~\ref{de:cover} and Lemma~\ref{le:extremal}, below). It is now not too difficult to deduce that $\rho$ is finite; 
we record this important fact as the following lemma.

\begin{lemma}\label{le:rhofinite}
$\rho < \infty$ for every critical update family $\U$. 
\end{lemma}

\begin{proof}
For each $u \in \Q_1$, set $\rho_0(u) = 0$ and for each $i \in \N$ with $i < \alphabar(u) \le \infty$, define
\[
\rho_i(u) := \sup \Big\{ \| y - Z \| : |Z| = i, \, y \in \big[ \H_u \cup Z \big] \setminus \H_u \Big\}.
\]
We shall prove by induction that each $\rho_i(u)$ is finite. Note that, since $\stabb$ is finite and $\alphabar(u) \geq \alpha$ for every $u\in\stabb$, by Lemma~\ref{le:stabb}, this will be sufficient to prove the lemma. Let $1 \le i < \alphabar(u)$, and assume that $\rho_{i-1}(u)$ is finite. Let $a > 0$ be a sufficiently large constant. There are various cases to deal with.

First, suppose that $\H_u(a u) \cap Z = \emptyset$. By Lemma~\ref{lem:BSU:extremal}, we have $\| y - Z \| = O(|Z|)$ for every $y \in [Z] = \big[ \H_u \cup Z \big] \setminus \H_u$, as required, where the last equality holds since $a$ is sufficiently large. So let $Z \subset \Z^2$ with $|Z| = i$ and $\H_u(a u) \cap Z \neq \emptyset$. Now, if
\[
Z \not\subset \H_u\big( 3ia u \big),
\]
then there exists a strip perpendicular to $u$ of width $2a$ that contains no element of $Z$, and separates some element of $Z$ from $\H_u(a u)$. Since $Z$ contains at least one element of $\H_u(a u)$, it follows that there are elements of $Z$ on both sides of this strip, and these two sets of elements cannot interact, by the induction hypothesis (and since $a$ is sufficiently large). It follows that we are done in this case.

We may therefore assume that $Z \subset \H_u\big( 3ia u \big)$. Now, if $[\H_u \cup Z] \setminus \H_u$ is infinite then it must contain an infinite number of elements of some line $\ell_u(j) \subset \H_u(3ia u) \setminus \H_u$, in which case there exists a translate $Z'$ of $Z$ such that $[\H_u \cup Z'] \cap \ell_u$ is infinite. But this contradicts our assumption that $i < \alphabar(u)$, so in fact $[\H_u \cup Z] \setminus \H_u$ is finite. 

Finally, observe that if $Z$ can be partitioned into two non-empty subsets at least distance $2\rho_{i-1}(u) + \nu$ apart, then they do not interact with each other (by the induction hypothesis), in which case we are again done. But now there are only a bounded number of choices for $Z$ (up to translation by an element of $\ell_u$), and for each of these $[\H_u \cup Z] \setminus \H_u$ is finite, so it follows that $\rho_i(u)$ is finite, as required.
\end{proof}

The constant $\kappa$ in Definition~\ref{de:connected} of a strongly connected set will be defined as follows:
\begin{equation}\label{eq:kappa}
\kappa = \kappa(\U) := 2\rho + 3\nu.
\end{equation}
Recall that $x,y \in \Z^2$ are said to be strongly connected to each other if $\| x - y \| \leq \kappa$. 

\begin{remark}\label{remark:wraparound}
To simplify the presentation, we will work on the infinite lattice $\Z^2$. However, it will be clear that the algorithms and lemmas below can be easily modified to the setting of the torus $\Z_n^2$, modulo some (easily resolved, but distracting\footnote{For example, we could simply set $D(K) = \Z_n^2$ for any set $K$ of diameter larger than $n$.}) technical issues that arise when the droplets have diameter $\Theta(n)$. Since the droplets in our applications (see Sections~\ref{se:lowerbal} and~\ref{se:lowerunbalend}) will all have diameter $(\log n)^{O(1)}$, we can reassure the concerned reader that these technical issues will not arise in practice. 
\end{remark}


\subsection{The covering algorithm}\label{se:covering}

In this subsection, all droplets will be $\stabb$-droplets. 
We will define the collection of \emph{$\alpha$-covers} of a finite set $K$ (see Definitions~\ref{de:cover} and~\ref{def:cover2}), and use this definition to prove two key lemmas: an `Aizenman--Lebowitz lemma' (see Lemma~\ref{le:AL}), which says that an $\alpha$-covered droplet contains $\alpha$-covered droplets of all intermediate sizes, and an extremal lemma (see Lemma~\ref{le:extremal}), which says that an $\alpha$-covered droplet contains many disjoint `$\alpha$-clusters'. The proofs of both lemmas are straightforward applications of the covering algorithm.

The key complication arising from the algorithm is that an $\alpha$-cover of a set $K$ does not necessarily contain the closure of $K$ under the $\U$-bootstrap process. However, an approximate version of this statement is true, and this is proved in Lemma~\ref{le:cover}. Roughly speaking, the lemma says that one can obtain (a superset of) the closure $[K]$ from an $\alpha$-cover of $K$ via only `local' modifications. 

We define an \emph{$\alpha$-cluster} to be any strongly connected set of $\alpha$ sites. These will be our basic building blocks in the covering algorithm. Recall that 
$D_B(K)$ denotes the unique minimal $\stabb$-droplet containing $K$.

\begin{definition}[\emph{The $\alpha$-covering algorithm}]\label{de:cover}
Let $\U$ be a critical update family. Suppose that we are given:
\begin{itemize}
\item $K$, a finite set of infected sites in $\Z^2$; \smallskip 
\item $B_1,\dots,B_{k_0}$, a maximal collection 
of disjoint $\alpha$-clusters in $K$;\smallskip
\item $\D^0 = \{D_1^0,\dots,D_{k_0}^0\}$, a collection of copies of a fixed, sufficiently large $\stabb$-droplet $\hat{D}$, such that $B_j \subset D_j^0$ for each $j \in \{1,\dots,k_0\}$.
\end{itemize}
Set $t := 0$ and repeat the following steps until STOP: 
\begin{itemize}
\item[1.] If there are two droplets $D_i^t,D_j^t \in \D^t$ and an $x \in \Z^2$ such that the set
\begin{equation}\label{eq:mergedroplets}
D_i^t \cup D_j^t \cup (x+\hat{D})
\end{equation}
is strongly connected, then set
\[
\D^{t+1} := \big( \D^t \setminus \{D_i^t,D_j^t \} \big) \cup \big\{ D_B(D_i^t \cup D_j^t) \big\},
\]
and set $t := t + 1$.
\item[2.] Otherwise set $T := t$ and STOP.
\end{itemize}
The output of the algorithm is the family $\D := \{D_1^T,\ldots,D_k^T\}$, where $k = k_0 - T$. 
\end{definition}

Thus, at each step of the algorithm, we take two nearby droplets in our collection, and replace them by the smallest $\stabb$-droplet containing their union. Let us fix from now on a sufficiently large $\stabb$-droplet $\hat{D}$ as in the covering algorithm. In particular, in Lemma~\ref{le:cover} we shall need that $\hat{D}$ contains a ball of radius $2\alpha\kappa$.

\begin{definition}\label{def:cover2}
We say that $\D = \{D_1,\ldots,D_k\}$ is an \emph{$\alpha$-cover} of a finite set $K \subset \Z^2$ if $\D$ is a possible output of the $\alpha$-covering algorithm with input $K$. 

We say that a (finite) droplet $D$ is \emph{$\alpha$-covered} by $A$ if the single droplet $\D = \{D\}$ is an $\alpha$-cover of $K$ for some subset $K \subset D \cap A$.
\end{definition}

We will show (see Lemmas~\ref{le:AL},~\ref{le:extremal} and~\ref{lem:cover:torus}) that if $A$ percolates then there exist $\alpha$-covered droplets at all scales, and that an $\alpha$-covered droplet of diameter $k$ must contain at least $\Omega(k)$ disjoint $\alpha$-clusters. It will then be straightforward to deduce the lower bound in Theorem~\ref{th:main} for balanced update families (see Section~\ref{se:lowerbal}).

The first important property of the $\alpha$-covering algorithm is given by the following lemma. We call this result an `Aizenman--Lebowitz lemma for $\alpha$-covered droplets', since the corresponding result for the 2-neighbour process was first proved in~\cite{AL}. Let $\lambda$ be a sufficiently large constant, depending on $\hat{D}$.

\begin{lemma}\label{le:AL}
Let $D$ be an $\alpha$-covered droplet. Then for every $\lambda \leq k \leq \diam(D)$ there exists an $\alpha$-covered droplet $D'\subset D$ such that $k \leq \diam(D') \leq 3k$.
\end{lemma}

\begin{proof}
The lemma is an immediate consequence of two simple observations: that the droplets $D_i^t \in \D^t$ obtained during the $\alpha$-covering algorithm are all $\alpha$-covered, and that at each step of the algorithm,
\[
\max\big\{ \diam(D_i^t) : D_i^t \in \D^t \big\}
\]
at most triples in size, provided that this maximum is at least an absolute constant (depending on $\hat{D}$). 

To prove the first observation, simply run the algorithm on a suitable subset of $D_i^t \cap A$, using the same $\alpha$-clusters and combining them in the same order. To prove the second, observe that if droplets $D_i^t$ and $D_j^t$ are united in step $t$ of the algorithm, then by definition there exists $x \in \Z^2$ such that the distance between $D_i^t$ and $x + \hat{D}$, and that between $D_j^t$ and $x + \hat{D}$, are at most $\kappa$. By Lemma~\ref{lem:subadditivity}, it follows that
$$\diam\big(D_B(D_i^t\cup D_j^t)\big) \leq \diam\big(D_i^t\big) + \diam\big(D_j^t\big) + O(1),$$
where the implicit constant depends on $\kappa$ and $\hat{D}$. This proves the second observation, and completes the proof of the lemma.
\end{proof}

The algorithm also allows us to prove the following extremal result, which says that the number of initial $\alpha$-clusters in an $\alpha$-covered droplet must be at least linear in the diameter of the droplet. This lemma is the reason why we use the $\alpha$-covering algorithm in the balanced setting, rather than the spanning algorithm (Definition~\ref{de:spanalg}), for which there is no correspondingly strong extremal lemma.

\begin{lemma}[Extremal lemma for $\alpha$-covered droplets]\label{le:extremal}
Let $D$ be an $\alpha$-covered droplet. Then $D \cap A$ contains $\Omega\big(\diam(D)\big)$ disjoint $\alpha$-clusters.
\end{lemma}

\begin{proof}
The algorithm begins with $k_0$ disjoint $\alpha$-clusters, and ends with $\D = \{D\}$. At each step of the algorithm the number of droplets is reduced by 1, and the sum of the diameters of the droplets increases by at most a constant, by Lemma~\ref{lem:subadditivity}. Hence
\[
\diam(D) \leq k_0\diam(\hat{D}) + O(k_0),
\]
and so $k_0 = \Omega\big(\diam(D)\big)$, as required.
\end{proof}

It remains to show that if $A$ percolates, then there exist large $\alpha$-covered droplets. The main step is Lemma~\ref{le:cover}, below, which shows that an $\alpha$-cover $\D$ of a set $K$ is a reasonable approximation of the closure $[K]$. The basic idea is simple: since all $\alpha$-clusters are contained in some droplet of $\D$, the remaining `dust' of infected sites -- that is, the set $K \setminus (D_1 \cup \cdots \cup D_k)$ -- should contribute only locally to the set of eventually infected sites. We remark that a simplified version of the covering algorithm was used in~\cite{BSUgen}, not requiring Lemma~\ref{le:cover}, and in most cases resulting in non-optimal bounds. 


Before stating Lemma~\ref{le:cover}, let us note that we can replace half-planes by $\stabb$-droplets in the definition~\eqref{eq:rho} of $\rho$. 

\begin{lemma}\label{le:rhodroplets}
Let 
$D$ be an $\stabb$-droplet, 
and let $Y \subset \Z^2$ have size at most $\alpha - 1$. Then $\|x - Y\| \leq \rho$ for all $x \in [D\cup Y]\setminus D$.
\end{lemma}

\begin{proof}
Let $\{a_u\in\Z^2:u\in\stabb\}$ be a collection of vectors such that
\[
D = \bigcap_{u\in\stabb} \H_u(a_u),
\]
and let $x\in[D\cup Y]\setminus D$. Since $x\notin D$, there exists $u\in\stabb$ such that $x\notin\H_u(a_u)$. But $|Y\setminus\H_u(a_u)| \le \alpha-1$, and
\[
x \in \big[\H_u(a_u)\cup Y\big] \setminus \H_u(a_u),
\]
so $\|x-Y\|\leq\rho$ by the definition of $\rho$.
\end{proof}

We are now ready to prove the key property of $\alpha$-covers. 

\begin{lemma}\label{le:cover}
Let $\,\U$ be a critical update family, let $K \subset \Z^2$ be a finite set, let $\D = \{D_1,\ldots,D_k\}$ be an $\alpha$-cover of $K$, and set $Y := K \setminus \big( D_1 \cup \cdots \cup D_k \big)$. Then 
$$\| x - Y \| \leq \rho$$ 
for every $x \in [K] \setminus \big(D_1 \cup \cdots \cup D_k \big)$.  
\end{lemma}


\begin{proof}
Observe first that we may assume that $\alpha \geq 2$, since otherwise the sets $Y$ and $[K] \setminus (D_1 \cup \dots \cup D_k)$ are empty (the latter because the covering algorithm stopped). We will prove a slightly stronger statement: setting
\[
X = \bigcup_{D \in \D} D \qquad \text{and} \qquad Z = \big[ X \cup Y \big],
\]
we will show that the same conclusion holds with $[K]$ replaced by $Z$.

To begin, we partition $Y$ into a collection $Y_1,\dots,Y_s$ of maximal strongly connected components, so in particular if $y \in Y_i$ and $z \in Y_j$ for some $i\neq j$, then
\begin{equation}\label{eq:coverYY}
\| y - z \| > \kappa = 2\rho + 3\nu.
\end{equation}
(Note that the sets $Y_i$ are uniquely defined.) By the definition of an $\alpha$-cover, we must have $|Y_i|\leq \alpha - 1$, and hence $\diam(Y_i) \leq (\alpha-2)\kappa$, for every $i \in [s]$.

For the clarity of what follows, we shall forget the labelling of the elements of $\D$ given in the statement of the lemma, so that we may reuse the notation $D_i$. Since $\hat{D}$ contains a ball of radius $2\alpha\kappa$, we may assume that
\begin{equation}\label{eq:coverDD}
\|D-D'\| > 2\alpha\kappa
\end{equation}
for every distinct pair $D,D' \in \D$. Thus, for each $i \in [s]$ there is at most one droplet $D_i \in \D$ such that $Y_i$ and $D_i$ are strongly connected to each other, since $\kappa + \diam(Y_i) + \kappa \leq \alpha\kappa$. (Of course we may have $D_i=D_j$ for distinct $i$ and $j$.) If there is no such $D_i$ then set $D_i = \emptyset$. In particular, if $D\in\D$ and $D \neq D_i$, then
\begin{equation}\label{eq:coverDY}
\|Y_i-D\| > \kappa.
\end{equation}
Now set $Y_i' := [D_i \cup Y_i] \setminus D_i$ for each $i \in [s]$, so that, by Lemma~\ref{le:rhodroplets}, we have
\begin{equation}\label{eq:coverYx}
\|x-Y_i\| \leq \rho
\end{equation}
for all $x \in Y_i'$ and all $i \in [s]$.

We claim that
\begin{equation}\label{eq:coverclaim}
X \cup Y_1' \cup \dots \cup Y_s' = Z;
\end{equation}
that is, that the set on the left-hand side is closed. Since the left-hand side of~\eqref{eq:coverclaim} may be re-written as
\[
X \cup [D_1 \cup Y_1] \cup \dots \cup [D_s \cup Y_s],
\]
and each of these sets is closed individually (in the case of $X$ this follows from~\eqref{eq:coverDD} and the fact that $\kappa > \nu$), it is enough to show that if $x \in Y_i'$, and either $y \in D$ for some $D \in \D$ with $D \neq D_i$, or $y \in Y_j'$ with $i \ne j$, then $x$ and $y$ are not close enough to interact; that is, $\|x-y\| > \nu$. Indeed, if $x \in Y_i'$ and $y \in D$ for some $D \in \D$ with $D \neq D_i$, then 
\[
\|x-y\| \geq \|Y_i-D\| - \|x-Y_i\| > \kappa - \rho \ge 2\nu
\]
by~\eqref{eq:coverDY} and~\eqref{eq:coverYx}. On the other hand, if $x \in Y_i'$ and $y \in Y_j'$, with $i \neq j$, then
\[
\|x-y\| \geq \|Y_i-Y_j\| - \|x-Y_i\| - \|y-Y_i\| > \kappa - 2\rho = \nu,
\]
by~\eqref{eq:coverYY} and~\eqref{eq:coverYx}. Thus,~\eqref{eq:coverclaim} holds.

We are now done, since we have shown that if $x \in [K] \setminus X$, then $x \in Y_i'$ for some $i \in [s]$, and we know that any such $x$ satisfies $\|x-Y_i\| \leq \rho$.
\end{proof}

In order to prove (in Section~\ref{se:lowerbal}) the lower bound in Theorem~\ref{th:main} for balanced update families, we will need only Lemma~\ref{le:extremal} and the following straightforward consequence of Lemmas~\ref{le:AL} and~\ref{le:cover}.

\begin{lemma}\label{lem:cover:torus}
Let $\,\U$ be a critical update family, and let $A \subset \Z_n^2$. If $[A] = \Z_n^2$, then there exists an $\alpha$-covered $\S_B$-droplet $D$ with 
\begin{equation}\label{eq:between:logs}
\log n \leq \diam(D) \leq 3 \log n.
\end{equation}
\end{lemma}

\begin{proof}
Run the $\alpha$-covering algorithm on $\Z_n^2$, with initial set $A$. If at some point we obtain an $\alpha$-covered droplet $D$ with $\diam(D) \geq \log n$, then choose the first such droplet, and observe that it satisfies~\eqref{eq:between:logs}, by the proof of Lemma~\ref{le:AL}. (Alternatively, choose any such droplet, and apply Lemma~\ref{le:AL} to it.)

So suppose that the $\alpha$-covering algorithm stops without creating any droplets of diameter larger than $\log n$, and let $\D = \{D_1,\ldots,D_k\}$ be the output of the algorithm. Setting $Y := A \setminus \big( D_1 \cup \cdots \cup D_k \big)$, and applying Lemma~\ref{le:cover}, it follows that
$$\| x - Y \| \leq \rho$$ 
for every $x \in [A] \setminus \big(D_1 \cup \cdots \cup D_k \big)$. Now, recall that $\| D_i - D_j \| \ge 2\alpha\kappa$ for each $i \ne j$ (since $\hat{D}$ was chosen sufficiently large), that $\kappa = 2\rho + 3\nu$, and that each strongly connected component of $Y$ has size at most $\alpha - 1$. Since $\diam(D) < \log n$ for each $D \in \D$, it follows that $[A] \ne \Z_n^2$, which contradicts our assumption.  
\end{proof}

\subsection{The spanning algorithm: unbalanced families}\label{se:spanning}

Next we describe our second analogue of the rectangles process, which will be a key tool in our analysis of unbalanced models. Throughout this subsection we assume that $\U$ is unbalanced\footnote{This is not strictly speaking necessary: unlike in the previous subsection, the results here hold for $\T$-droplets for any $\T \subset \stab$ such that $\T$-droplets are finite. Nevertheless, the only applications of the results in this subsection will be to unbalanced families, and it is useful to fix the set $\T$.} and that droplets are taken with respect to $\stabu$. Recall that we write $D_U(K)$ to denotes the smallest $\S_U$-droplet containing $K$. We also remind the reader that we define the algorithm in $\Z^2$ to avoid some (unimportant) technical details relating to strongly connected sets of diameter $\Theta(n)$. 

Recall from Section~\ref{se:sketch} that an $\stabu$-droplet $D$ is said to be internally spanned by~$A$ if there exists a strongly connected set $L \subset [D \cap A]$ such that $D_U(L) = D$. Given a finite set $K$ of infected sites, the output of the spanning algorithm is a minimal collection $\D$ of internally spanned $\stabu$-droplets whose union contains $K$. At each step of the algorithm we maintain a partition $\K^t = \{K^t_1, \ldots, K^t_k\}$ of $K$ such that each set $[K_j^t]$ is strongly connected. 

\begin{definition}[\emph{The spanning algorithm}]\label{de:spanalg}
Let $K = \{x_1,\ldots,x_{k_0}\} \subset \Z^2$ be a set of infected sites. Set $\K^0 := \{K_1^0,\ldots,K_{k_0}^0 \}$, where $K_j^0 := \{x_j\}$ for each $1\leq j \leq k_0$. Set $t := 0$, and repeat the following steps until STOP: 
\begin{itemize}
\item[1.] If there are two sets $K_i^t,K_j^t \in \K^t$ such that the set
\begin{equation}\label{eq:mergespanningsets}
\big[ K_i^t \big] \cup \big[ K_j^t \big]
\end{equation}
is strongly connected, then set
\[
\K^{t+1} := \big( \K^t \setminus \{K_i^t,K_j^t \} \big) \cup \big\{ K_i^t \cup K_j^t \big\},
\]
and set $t := t + 1$. 
\item[2.] Otherwise set $T := t$ and STOP.
\end{itemize}
The output of the algorithm is the \emph{span of $K$}, 
\[
\< K \> := \big\{ D_U\big( [K_1^T] \big),\ldots,D_U\big( [K_k^T] \big) \big\},
\]
where $k = k_0 - T$. 
\end{definition}

The following lemma shows that the span is well-defined, and provides an alternative description of the span of a set $K$. 

\begin{lemma}\label{le:spanoutput}
For every finite set $K$, we have 
\begin{equation}\label{eq:spanofK}
\< K \> = \big\{D_U(L_1),\ldots,D_U(L_k)\big\},
\end{equation}
where $L_1,\ldots,L_k$ are the strongly connected components of $[K]$.
\end{lemma}

\begin{proof}
We shall show that the sets $[K_i^T]$ are precisely the strongly connected components of $[K]$. Indeed, it follows from~\eqref{eq:mergespanningsets} (and a simple induction on $t$) that $[K_i^t]$ is strongly connected for every $t \in [T]$ and $1\leq i \leq k_0 - t$, and it follows that no two sets $[K_i^T]$ and $[K_j^T]$ are strongly connected to each other, since the algorithm stopped at step $T$. We therefore have $[K] = \bigcup_{i = 1}^k [K_i^T]$, since, by~\eqref{eq:kappa}, the constant $\kappa$ (from the definition of strongly connected) is greater than $\nu$ (the diameter of $\U$), and so no site can be infected by two or more of these sets.  
\end{proof}

It is now easy to deduce that we can use the spanning algorithm to determine whether or not $D$ is internally spanned.

\begin{lemma}\label{le:span}
An $\stabu$-droplet $D$ is internally spanned if and only if $D \in \< D \cap A \>$. 
\end{lemma}

\begin{proof}
Applying Lemma~\ref{le:spanoutput} to $K = D \cap A$, we see that $D \in \< D \cap A \>$ if and only if $D_U(L) = D$ for some strongly connected component $L$ of $[D \cap A]$. But $[D \cap A] \subset D$, since $\stabu \subset \stab$, and so this is equivalent to the event that $D$ is internally spanned.
\end{proof}

We can now prove the `Aizenman--Lebowitz lemma for internally spanned droplets', which is the spanning analogue of Lemma~\ref{le:AL} for $\alpha$-covered droplets. For the applications we shall need a slightly more general statement than before. Recall that $\pi(D,u)$ denotes the size of the projection of $D$ in the direction $u$. We shall once again use $\lambda$ to denote a sufficiently large constant.

\begin{lemma}\label{le:ALspan}
Let $D$ be an internally spanned $\stabu$-droplet, and let $u \in S^1$. Then for every $\lambda \leq k \leq \pi(D,u)$, there exists an internally spanned $\stabu$-droplet $D' \subset D$ with 
$$k \leq \pi(D',u) \leq 3k.$$
\end{lemma}

\begin{proof}
Apply the spanning algorithm to $K = D \cap A$ and observe that, for every $t \leq T$ and every $1\leq i \leq k_0 - t$, the $\stabu$-droplet $D_U([K_i^t])$ is internally spanned by $A$, since $K_i^t \subset D_U([K_i^t]) \cap A$, and the set $[K_i^t]$ is strongly connected.

We claim that 
\[
\max\big\{ \pi\big( D_U([K_i^t]), u \big) : K_i^t \in \K^t \big\}
\]
at most triples in size at each step, provided that this maximum is at least an absolute constant. Indeed, by Lemma~\ref{lem:subadditivity}, we have
$$\pi\big(D_U(D_1\cup D_2),u\big) \leq \pi(D_1,u) + \pi(D_2,u) + O(1),$$
for any pair of droplets $D_1$ and $D_2$ that are within distance $O(1)$ of one another, and that $D_U([Y]) = D_U(Y)$ for any set $Y$, since $\stabu \subset \stab$. The lemma now follows easily, as in the proof of Lemma~\ref{le:AL}.
\end{proof}

We can now deduce an extremal lemma which, while much weaker than the corresponding lemma for $\alpha$-covered droplets (Lemma~\ref{le:extremal}), is in fact tight up to the implicit constant. This fact underlines how much we are `giving away' in assuming only that our droplets are spanned (rather than filled). Nevertheless, this lemma will be sufficient to prove the base case (Lemma~\ref{le:smalldroplet} below) of the main induction argument (Lemma~\ref{le:induction}) for unbalanced models in Section~\ref{se:lowerunbal}.

\begin{lemma}\label{le:extremalspan}\emph{(Extremal lemma for internally spanned droplets.)}
Let $D$ be an internally spanned $\stabu$-droplet. Then $|D\cap A| = \Omega\big(\diam(D)\big)$.
\end{lemma}

\begin{proof}
As in the proof of the previous lemma, we apply the spanning algorithm with $K=D\cap A$. The algorithm starts with $k_0$ sets containing the individual elements of $D\cap A$, and it finishes with a collection
\[
\<D\cap A\> = \big\{ D_U\big( [K_1^T] \big),\ldots,D_U\big( [K_k^T] \big) \big\}
\]
such that $D\in\<D\cap A\>$. At each step of the algorithm the number of sets in the collection decreases by 1, and the sum of the diameters of the minimal droplets containing those sets increases by at most a constant, by Lemma~\ref{lem:subadditivity}. Hence,
\[
\diam(D) \leq \sum_{i=1}^k \diam\big(D_U([K_i^T])\big) \leq k_0 \diam\big(D_U([K_1^0])\big) + O(k_0) = O(k_0),
\]
which implies that $k_0 = \Omega\big(\diam(D)\big)$, as required.
\end{proof}

Using Lemma~\ref{le:extremalspan}, we can deduce a non-trivial bound on the probability that a very small droplet is internally spanned. As noted above, this will form the base case of our induction argument in Lemma~\ref{le:induction}.

Recall that $\ispan(D)$ denotes the event that the $\stabu$-droplet $D$ is internally spanned, and that $w(D)$ and $h(D)$ denote its width and height respectively, as defined in Section~\ref{se:defs}.

\begin{lemma}\label{le:smalldroplet}
For every $\eta > 0$, there exists $\delta > 0$ such that the following holds. Let $D$ be an $\stabu$-droplet such that
\[
\min\big\{w(D),h(D)\big\} \leq p^{-1+\eta}.
\]
Then
\[
\Pr_p\big(\ispan(D)\big) \leq p^{\delta \max\{w(D), \, h(D)\}}.
\]
\end{lemma}

\begin{proof}
Let us write $m(D):=\min\big\{w(D),h(D)\big\}$ and $M(D):=\max\big\{w(D),h(D)\big\}$. Suppose the $\stabu$-droplet $D$ is internally spanned. Then by Lemma~\ref{le:extremalspan}, $D\cap A$ must contain at least $\Omega\big(M(D)\big)$ sites. The probability that this occurs is at most
\[
\binom{O\big(w(D)\cdot h(D)\big)}{\delta' \cdot M(D)} p^{\delta' M(D)} \leq \big( O(1) \cdot m(D) \cdot p \big)^{\delta' M(D)} \leq p^{\delta M(D)},
\]
for some $\delta,\delta'>0$, as required.
\end{proof}

The final lemma of this subsection will be used in Section~\ref{se:hierarchies} as part of an induction argument to prove the existence of `good and satisfied hierarchies' for internally spanned droplets (see Definitions~\ref{de:hierarchy}--\ref{def:satisfied} and Lemma~\ref{le:goodsat} for details).

\begin{lemma}\label{le:penultimate}
Let $K\subset\Z^2$, with $2\leq|K|<\infty$, be such that $[K]$ is strongly connected. Then there exists a partition $K=K_1\cup K_2$ into non-empty (disjoint) sets such that each of the sets $[K_1]$, $[K_2]$ and $[K_1]\cup[K_2]$ is strongly connected.
\end{lemma}

\begin{proof}
Run the spanning algorithm on $K$ and consider the penultimate step. Since $[K]$ is strongly connected, and therefore $\< K\>=\big\{ D_U\big([K]\big)\big\}$, by Lemma~\ref{le:spanoutput}, we have
\[
\K^{T-1} = \big\{ K_1,K_2 \big\}
\]
for some $K_1\subsetneq K$ and $K_2\subsetneq K$ such that $K=K_1\cup K_2$. By their construction in the spanning algorithm, both $[K_1]$ and $[K_2]$ are strongly connected, and since $K_1$ and $K_2$ combine at the final step, so too is $[K_1] \cup [K_2]$.
\end{proof}

\subsection{The iceberg algorithm: unbalanced families with drift}\label{se:iceberg}

Our third algorithm will play a crucial role in the proof for update families that exhibit drift. Assume that $\U$ is unbalanced and let $\{u^*,-u^*\} \subset \stabu$ be the pair of stable directions given by Lemma~\ref{le:stabu} (and originally by Lemma~\ref{le:ustar}), so in particular
\[
\min\big\{ \alpha(u^*), \alpha(-u^*) \big\} \geq \alpha + 1. 
\]
Recall from the discussion after~\eqref{eq:unbal:driftcase} that we say a model \emph{exhibits drift} if 
$$\max\big\{ \alpha^+(u^*), \alpha^-(u^*), \alpha^+(-u^*), \alpha^-(-u^*) \big\} = \infty.$$
If $\U$ does not exhibit drift then we shall not need the results proved in this section, so let us assume that $\alpha^-(u^*) = \infty$.




When our droplet is growing in direction $u^*$ in a model with drift, it will tend to form a triangle, as in Section~\ref{se:upperunbal}. In order to control the growth in this direction, we therefore need to `give away' this triangle (in fact, a slightly larger one), and bound the growth outside it. The point of the algorithm defined in this subsection is exactly to control this outwards growth using `icebergs', defined as follows.

Since $\alpha^-(u^*) = \infty$, there exists a non-trivial interval $[u^*,u_0]$ such that $\alpha^-(u) = \infty$ for every $u \in [u^*,u_0]$, by Lemma~\ref{le:alphainfinity}. Fix such a $u_0$ sufficiently close to $u^*$ in the following sense: we choose $u_0$ to be closer to $u^*$ than any $v \in S^1 \setminus \{u^*\}$ perpendicular to $x - y$, where $x,y \in \bigcup_{X \in \U} X \cup \{\0\}$ and $x \ne y$. That we can choose such a $u_0$ follows easily from the fact that $\U$ is a finite collection of finite sets. Finally, choose $u_1 \in (u^*,u_0)$ arbitrarily.

\begin{definition}\label{de:iceberg}
Let $u \in (u^*,u_1]$. A \emph{$u$-iceberg} is any non-empty set $J$ of the form
\[
J = \big( \H_{u_0}(a) \cap \H_{u^*}(b) \big) \setminus \H_u,
\]
where $a,b \in \R^2$. If $X$ is a finite set of sites such that $X\not\subset\H_u$, then denote by $J_u(X)$ the smallest $u$-iceberg such that $X\subset \H_u\cup J_u(X)$.
\end{definition}

\begin{figure}[ht]
  \centering
  \begin{tikzpicture}[>=latex,scale=1.1]
    \draw (0,0) -- (12,0) (0,0) -- (12,3);
    \draw (3,0.75) -- (4,2.25) -- (9,2.25);
    \draw [<->,densely dashed] (3,2.75) -- node [above] {$O\big(\gamma/\sigma(u)\big)$} (9,2.75);
    \draw [<->,densely dashed] (1.5,0.75) -- node [left] {$O(\gamma)$} (1.5,2.25);
    \node at (4.5,1.7) {$J$};
    \draw [->] (11,2.75) -- ($(11,2.75)!0.5cm!(10,6.75)$);
    \node at ($(11,2.75)!0.8cm!(10,6.75)$) {$u$};
    \draw [->] (11,0) -- (11,0.5);
    \node at (11,0.8) {$u^*$};
    \draw [->] (3.5,1.5) -- ($(3.5,1.5)!0.5cm!(2,2.5)$);
    \node at ($(3.5,1.5)!0.9cm!(2,2.5)$) {$u_0$};
    \draw (3.5,0) arc (0:atan(0.25):3.5);
    \pgfmathparse{atan(0.125)};
    \node at (\pgfmathresult:4.4) {$\sigma(u)$};
  \end{tikzpicture}
  \caption{A $u$-iceberg $J$, together with bounds on its width and height given by Lemma~\ref{le:extremaliceberg}.}\label{fi:iceberg}
\end{figure}
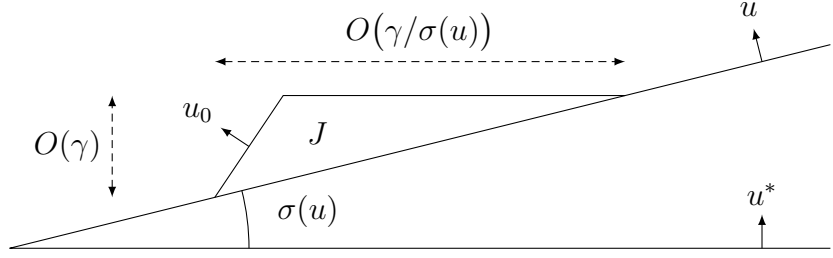

Thus a $u$-iceberg is a discrete triangle whose sides are perpendicular to $-u$, $u^*$ and $u_0$; see Figure~\ref{fi:iceberg}. The role of $u_1$ is to ensure that that the angle between the $(-u)$-side and the $u_0$-side of a $u$-iceberg is uniformly bounded away from zero: this will be important in Lemma~\ref{le:extremaliceberg}. We make a simple but key observation, which follows easily from the definition of~$u_0$ (cf.~the proofs of Lemmas~\ref{le:alphainfinity} and~\ref{le:quasi}).

\begin{lemma}\label{le:icebergclosed}
If $J$ is a $u$-iceberg, then $\H_u \cup J$ is closed.
\end{lemma}

\begin{proof}
It suffices to show that there does not exist a rule $X \in \U$ such that either $X \subset \H_u \cup \H_{u_0}$ or $X \subset \H_u \cup \H_{u^*}$. Note that, since $\{u,u^*,u_0\} \subset \S$, we cannot have $X \subset \H_u$, $X \subset \H_{u_0}$, or $X\subset \H_{u^*}$. Hence, if $X \subset \H_u \cup \H_{u_0}$ then there exist $x,y \in X$ with $x \in \H_u \setminus \H_{u_0}$ and $y \in \H_{u_0} \setminus \H_u$, and similarly, if $X \subset \H_u \cup \H_{u^*}$ then there exist $x,y \in X$ with $x \in \H_u \setminus \H_{u^*}$ and $y \in \H_{u^*} \setminus \H_u$. In either case, $x - y$ is perpendicular to a vector in the interval $(u^*,u_0)$, contradicting our choice of~$u_0$.
\end{proof}

We are now ready to introduce the iceberg algorithm, which is a modified version of the covering algorithm allowing sites to be infected with the help of $\H_u$. At each step of the algorithm we have a collection $\W^t$ of $\S_U$-droplets and $u$-icebergs; we either take a droplet near $\H_u$ and replace it by the smallest $u$-iceberg containing it, or we take two nearby sets in our collection, and replace them by either the smallest $u$-iceberg containing their union (if they are sufficiently close to $\H_u$), or by the smallest droplet containing their union (otherwise).

\begin{definition}[\emph{The $u$-iceberg algorithm}]\label{de:icebergalg}
Let $\U$ be an unbalanced update family that exhibits drift, let $u^*$, $u_0$ and $u_1$ be as defined above, and let $u \in (u^*,u_1]$. Suppose we are given:
\begin{itemize}
\item $K = \big\{ x_1,\ldots,x_{k_0} \big\} \subset \Z^2 \setminus \H_u$, a finite set of infected sites;\smallskip
\item $\W^0 = \{W_1^0,\dots,W_{k_0}^0\}$, a collection of copies of a fixed, sufficiently large $\stabu$-droplet $\hat{D}_U$, such that $x_j \in W_j^0$ for each $j \in \{1,\dots,k_0\}$.
\end{itemize}
Set $t := 0$ and repeat the following steps until STOP: 
\begin{itemize}
\item[1.] If there is a droplet $W_i^t \in \W^t$ and an $x \in \Z^2$ such that the set
\[
W_i^t \cup (x + \hat{D}_U) \cup \H_u
\]
is strongly connected, then set
\[
\W^{t+1} := \big( \W^t \setminus \{W_i^t \} \big) \cup \big\{ J_u(W_i^t) \big\},
\]
and set $t := t + 1$. \smallskip
\item[2.] If not, but there are two sets\footnote{Note that in this case at least one of the sets must be an iceberg.} $W_i^t,W_j^t \in \W^t$ and an $x \in \Z^2$ such that the sets
\[
W_i^t \cup W_j^t \cup (x+\hat{D}_U) \quad \text{and} \quad W_i^t \cup W_j^t \cup (x+\hat{D}_U) \cup \H_u
\]
are both strongly connected, then set
\[
\W^{t+1} := \big( \W^t \setminus \{W_i^t,W_j^t \} \big) \cup \big\{ J_u(W_i^t \cup W_j^t) \big\},
\]
and set $t := t + 1$.\smallskip
\item[3.] If not, but there are two droplets $W_i^t,W_j^t \in \W^t$ and an $x \in \Z^2$ such that the set
\[
W_i^t \cup W_j^t \cup (x+\hat{D}_U)
\]
is strongly connected, then set
\[
\W^{t+1} := \big( \W^t \setminus \{W_i^t,W_j^t \} \big) \cup \big\{ D_U(W_i^t \cup W_j^t) \big\},
\]
and set $t := t + 1$.\smallskip
\item[4.] Otherwise set $T := t$ and STOP.
\end{itemize}
The output of the algorithm is the family $\W := \{W_1^T,\ldots,W_k^T\}$.
\end{definition}

\begin{definition}
Let $u \in (u^*,u_1]$. We say that $\W = \{W_1,\ldots,W_k\}$ is a \emph{$u$-iceberg cover} of a finite set $K$ if $\W$ is a possible output of the $u$-iceberg algorithm with input $K$. We say that a $u$-iceberg (or $\stabu$-droplet) $J$ is \emph{$u$-iceberg covered} if $\W = \{ J \}$ is a $u$-iceberg cover of $K$ for some subset $K \subset J \cap A$.
\end{definition}

Before continuing, let us note that $u$-iceberg covers are closed.

\begin{lemma}\label{lem:iceberg:covers:are:closed}
Let $u \in (u^*,u_1]$, let $K \subset \Z^2$ be a finite set, and let $\W$ be a $u$-iceberg cover of $K$. Then the set
\[
\H_u \cup \bigcup_{W \in \W} W
\]
is closed under the $\U$-bootstrap process. 
\end{lemma}

\begin{proof}
Since the algorithm has terminated, and $\hat{D}_U$ is sufficiently large, no two elements of $\W$ are strongly connected to each other, and any element of $\W$ strongly connected to $\H_u$ must be a $u$-iceberg. The lemma now follows by Lemma \ref{le:icebergclosed}. 
\end{proof}

We can now prove our extremal result for icebergs; the lemma is illustrated in Figure~\ref{fi:iceberg}. Let $\sigma(u)$ denote the angle (in radians) between $u$ and $u^*$.

\begin{lemma}[Extremal lemma for $u$-iceberg covers]\label{le:extremaliceberg}
Let $u \in (u^*,u_1]$, let $J$ be a $u$-iceberg covered $u$-iceberg \rob{or $\stabu$-droplet}, and let $\gamma = |J \cap A|$. Then
\[
w(J) = O\big( \gamma / \sigma(u) \big) \quad \text{and} \quad h(J) = O(\gamma),
\]
where the implicit constants depend on\/ $\U$, and on the fixed directions $u^*$, $u_0$ and $u_1$, but not on $J$, $\gamma$ or $u$.
\end{lemma}

\begin{proof}
Suppose first that $J$ is an $\stabu$-droplet, and recall that if $J$ is $u$-iceberg covered, then there exists $K \subset J \cap A$ and a way of running the $u$-iceberg algorithm with input $K$ such that $\{J\}$ is the output. Note that in this case every step of the algorithm must have been of type 3, since there are no icebergs in the output. By Lemma~\ref{lem:subadditivity}, in each such step the sum of the widths and heights of the droplets in $\W^t$ goes up by at most a constant, and hence $w(J) + h(J) = O(\gamma)$. 

Let us therefore assume that $J$ is a $u$-iceberg, and observe that\footnote{The implicit constants here and subsequently may increase as the angle between $u$ and $u_0$ decreases, but this angle is bounded below by a function of $u_0$ and $u_1$ only, and the statement of the lemma permits such a dependence.} 
$$w(J) = O\big( h(J) / \sigma(u) \big),$$
since the angle between $u$ and $u_0$ is bounded away from zero. It will therefore suffice to prove that $h(J) = O(\gamma)$.

To do so, observe first that in the $u$-iceberg algorithm (applied with $K = J \cap A$), droplets can be replaced by icebergs, but not vice versa. It follows that
\begin{equation}\label{eq:iceT}
T \le 2\gamma,
\end{equation}
since at all but at most $\gamma$ steps of the algorithm, $|\W^t|$ is reduced by 1.  

Now, let $\D^t$ and $\J^t$ denote, respectively, the collections of droplets and $u$-icebergs in $\W^t$, so $\W^t=\D^t\cup\J^t$. Let $\lambda > 0$ be a sufficiently large constant, depending on the angle between $u$ and $u_0$. We claim that, for each $t \leq T$,
\begin{equation}\label{eq:iceclaim2}
\lambda \sum_{D^t\in\D^t} \big( \sigma(u) \cdot w(D^t) + h(D^t) \big) + \sum_{J^t\in\J^t} h(J^t) = O(t+\gamma) = O(\gamma), 
\end{equation}
The second equality follows from~\eqref{eq:iceT}; we will prove the first by induction on $t$. Note that the inequality is trivial when $t=0$, so it will be enough to prove that the left-hand side of~\eqref{eq:iceclaim2} increases by at most $O(1)$ at each step. When two $u$-icebergs are replaced by another $u$-iceberg (as in step 2 of the $u$-iceberg algorithm), or when two droplets are replaced by another droplet (as in step 3 of the algorithm), this follows\footnote{More precisely, for $u$-icebergs this follows from the proof of Lemma~\ref{lem:subadditivity}, since $\S_U$ is a finite set of rational directions that intersects every open semicircle.} from Lemma~\ref{lem:subadditivity}, since the sums individually increase by at most $O(1)$. On the other hand, when a droplet $D^t$ is replaced by a $u$-iceberg, or a droplet $D^t$ and a $u$-iceberg are replaced by a $u$-iceberg, the first sum decreases by $\lambda \big( \sigma(u) w(D^t) + h(D^t) \big)$ and the second sum increases by at most $O\big( \sigma(u) w(D^t) + h(D^t) + 1 \big)$. Since $\lambda$ was chosen sufficiently large, it follows that~\eqref{eq:iceclaim2} holds, and hence, since $\W^T = \{J\}$, 
\[
h(J) = O(\gamma),
\]
as required.
\end{proof}

\section{The lower bound for balanced families}\label{se:lowerbal}

In this section we complete the proof of the lower bound in Theorem~\ref{th:main} for balanced\footnote{The proof also applies to unbalanced update families, but the stronger bound required for the theorem will be proved for those families in Section~\ref{se:lowerunbal}.}  update families. The proof is a straightforward consequence of the $\alpha$-covering algorithm, and the lemmas proved in Section~\ref{se:covering}.

\begin{theorem}\label{th:lower}
Let\/ $\U$ be a critical update family. Then
\[
p_c(\Z_n^2,\U) = \Omega\bigg( \frac{1}{\log n} \bigg)^{1/\alpha}
\]
\end{theorem}

\begin{proof}
Let $A$ be a $p$-random subset of $\Z_n^2$, where 
\[
p = \bigg( \frac{\eps}{\log n} \bigg)^{1/\alpha},
\]
for some sufficiently small constant $\eps = \eps(\U) > 0$. We shall show that, with high probability as $n \to \infty$, the $\U$-bootstrap closure of $A$ is not equal to $\Z_n^2$. 


Indeed, by Lemma~\ref{lem:cover:torus}, if $[A] = \Z_n^2$ then there exists an $\alpha$-covered droplet $D$ with 
$$\log n \leq \diam(D) \leq 3\log n,$$ 
and by Lemma~\ref{le:extremal} it follows that $D \cap A$ contains at least $\delta\log n$ disjoint $\alpha$-clusters, for some constant $\delta = \delta(\U) > 0$. 

Noting that $D$ contains $O\big( \diam(D)^2 \big) = O( \log n )^2$ distinct $\alpha$-clusters, it follows that the probability $D$ is $\alpha$-covered is at most
\[
\binom{O(\log n)^2}{\delta \log n} p^{\alpha \delta \log n} \leq \Big(O\big( p^\alpha \log n \big)\Big)^{\delta \log n} \leq \frac{1}{n^3},
\]
since $\epsilon$ is sufficiently small. Finally, since there are at most $n^2 (\log n)^{O(1)}$ choices of the droplet $D$ having diameter at most $3\log n$, it follows that
\[
\Pr_p\big( [A] = \Z_n^2 \big) \le n^2 \cdot (\log n)^{O(1)} \cdot \frac{1}{n^3} = o(1),
\]
as required.
\end{proof}

Note that we actually proved a stronger result than that stated in Theorem~\ref{th:lower}: it follows from the proof above that if $\eps = \eps(\U) > 0$ is a sufficiently small constant, 
\[
p = \bigg( \frac{\eps}{\log n} \bigg)^{1/\alpha},
\]
and $A$ is a $p$-random subset of $\Z_n^2$, then with high probability all strongly connected components of $[A]$ have diameter $O(\log n)$.

\section{The lower bound for unbalanced families}\label{se:lowerunbal}

In this section we shall prove the following theorem, and hence complete the proof of Theorem~\ref{th:main}.

\begin{theorem}\label{th:lowerunbal}
Let $\U$ be an unbalanced critical update family. Then
\[
p_c(\Z_n^2,\U) = \Omega\bigg( \frac{(\log\log n)^2}{\log n} \bigg)^{1/\alpha}.
\]
\end{theorem}

Throughout the section we assume that $\U$ is unbalanced, and that droplets are taken with respect to the set $\stabu=\{u^*,-u^*,u^r,u^l\}$ as in Lemma~\ref{le:stabu}, where
\begin{equation}\label{eq:choosing:stab:U}
\min\big\{ \alpha(u^*),\alpha(-u^*) \big\} \geq \alpha+1 \qquad \text{and} \qquad \min\big\{ \alphabar(u^l),\alphabar(u^r) \big\} \ge \alpha,
\end{equation}
and $u^l$ and $u^r$ are contained in opposite semicircles separated by $u^*$ and $-u^*$, with $u^l$ to the left and $u^r$ to the right of $u^*$. We also let $\xi>0$ be a sufficiently small constant (which will depend in particular on the constant $\delta(2\alpha+1)$ defined below; see~\eqref{eq:constants1}), and we fix
\begin{equation}\label{def:eta}
\eta := \frac{1}{10\alpha}.
\end{equation}

The main step in the proof of Theorem~\ref{th:lowerunbal} is an upper bound on the probability that a critical droplet is internally spanned. Recall from Definition~\ref{de:critical} that an $\S_U$-droplet $D$ is said to be \emph{critical} if its dimensions satisfy either
\begin{enumerate}
\item[$(T)$] $w(D) \leq 3p^{-\alpha-1/5}$ and $\frac{\xi}{p^\alpha}\log\frac{1}{p}\leq h(D)\leq\frac{3\xi}{p^\alpha}\log\frac{1}{p}$, or \smallskip
\item[$(L)$] $p^{-\alpha-1/5}\leq w(D)\leq 3p^{-\alpha-1/5}$ and $h(D) \leq \frac{\xi}{p^\alpha}\log\frac{1}{p}$.
\end{enumerate}
The key bound we shall prove is that there exists $\delta > 0$ such that if $D$ is a critical droplet then
\begin{equation}\label{eq:critunbal}
\Pr_p\big(\ispan(D)\big) \le \exp\left( -\frac{\delta}{p^\alpha}\left(\log\frac{1}{p}\right)^2 \right).
\end{equation}
(Recall again that $\ispan(D)$ is the event that the $\stabu$-droplet $D$ is internally spanned.) The proof of~\eqref{eq:critunbal} will only be given towards the end of this section, in Lemma~\ref{le:critunbal}. We build up to that proof gradually via an induction argument, at each step of which we bound the probability that droplets of certain (increasingly large) sizes are internally spanned.

During the course of this section we shall use a large number of constants, with various dependencies. The main constants we shall use are $\delta(\beta)$, for $2\leq\beta\leq 2\alpha+1$, and the constants $\delta$ (which will appear in Lemma~\ref{le:critunbal}), $\xi$ (from the definition of a critical droplet), and $\epsilon$ (which will be used in the proof of Theorem~\ref{th:lowerunbal} in Section~\ref{se:lowerunbalend}). These will be chosen so that
\begin{equation}\label{eq:constants1}
1 \gg \delta(2) \gg \dots \gg \delta(2\alpha+1) \gg \xi \gg \delta \gg \epsilon > 0,
\end{equation}
by which we mean that the constants are chosen from left to right, so that each is sufficiently small depending on all previous constants. Later in the section we shall introduce two further sequences of constants. The relationships between these new constants and those in~\eqref{eq:constants1} will be set out explicitly in~\eqref{eq:constants2} and~\eqref{eq:constants3}, below.

Next we state the induction hypothesis.

\begin{definition}\label{de:ih}
For each $\beta_1,\beta_2\in\N$ with $\beta_1+\beta_2 \leq 2\alpha+1$, let $\ih(\beta_1,\beta_2)$ be the following statement:
\begin{itemize}
\item [] \emph{Let $D$ be an $\stabu$-droplet such that
\[
w(D) \leq p^{-\beta_1(1-2\eta)-\eta} \quad \text{and} \quad h(D) \leq p^{-\beta_2(1-2\eta)-\eta}.
\]
Then
\begin{equation}\label{eq:ihbound}
\Pr_p\big(\ispan(D)\big)\leq p^{\delta\max\{w(D), \, h(D)\}},
\end{equation}
where $\delta = \delta(\beta_1+\beta_2)$.}
\end{itemize}
\end{definition}

We mention briefly that we would prefer the width and height conditions in Definition~\ref{de:ih} to be $w(D)\leq p^{-\beta_1(1-\eta)}$ and $h(D)\leq p^{-\beta_2(1-\eta)}$ respectively, but for technical reasons we cannot quite square the bound on the width between $\beta_1=1$ and $\beta_1=2$; this is why the conditions take the slightly less elegant form above.

The specific induction statements that we shall prove are:
\begin{align*}
\ih(\beta,\beta) \, &\Rightarrow \, \ih(\beta+1,\beta) && \text{for all } 1\leq\beta\leq\alpha; \\
\ih(\beta,\beta) \, &\Rightarrow \, \ih(\beta,\beta+1) && \text{for all } 1\leq\beta\leq\alpha; \\
\big( \ih(\beta+1,\beta) \, \wedge \, \ih(\beta,\beta+1) \big) \, &\Rightarrow \, \ih(\beta+1,\beta+1) && \text{for all } 1\leq\beta\leq\alpha-1.
\end{align*}
Note that $\ih(1,1)$ is an immediate consequence of Lemma~\ref{le:smalldroplet}, and therefore together these statements will be enough to prove the following lemma.

\begin{lemma}\label{le:induction}
The assertions $\ih(\alpha+1,\alpha)$ and $\ih(\alpha,\alpha+1)$ both hold.
\end{lemma}

Lemma~\ref{le:induction} alone is not enough to give the bound~\eqref{eq:critunbal} that we want on internally spanned critical droplets. However, the techniques and lemmas that we use to prove Lemma~\ref{le:induction} will be the same as those that we use in the proof of Lemma~\ref{le:critunbal} to deduce~\eqref{eq:critunbal}.

The steps in the induction are of two types: \emph{horizontal steps} of the form
\[
\ih(\beta_1,\beta_2)\Rightarrow\ih(\beta_1+1,\beta_2)\text,
\]
and \emph{vertical steps} of the form
\[
\ih(\beta_1,\beta_2)\Rightarrow\ih(\beta_1,\beta_2+1)\text.
\]
Common to both is the key idea of \emph{crossings}. Roughly speaking, these are events that say that it is possible to `cross' a parallelogram of sites from one side to the other with `help' from one of the sides in the form of an infected half-plane. The events should be thought of in the context of a growing droplet: a combination of crossing events, one for each side of the droplet, enable an internally filled droplet to grow into a larger internally spanned droplet. We obtain bounds for the probabilities of crossings by showing that, to a certain level of precision, the most likely way these events could occur is via the droplet (or half-plane) advancing row-by-row, rather than via the merging of many smaller droplets. One could think of this as saying that the growth mechanism we used to prove the upper bound for unbalanced families in Theorem~\ref{th:upperunbal}, which was indeed row-by-row, was essentially the `correct' mechanism. For vertical crossings in the case of models with drift, our proof will make use of the results of Section~\ref{se:iceberg} on the iceberg algorithm to bound the range of the $\U$-bootstrap process in directions close to $\pm u^*$. Full statements and proofs of the crossing lemmas, together with precise definitions, are given in Section~\ref{se:crossings}.

For the horizontal steps, in addition, we require the use of `hierarchies' to bound the extent of sideways growth at any given step. These are by now a standard tool in the bootstrap percolation literature, so we omit some of the details.
 
There are nine subsections in this section, which deal with the following aspects of the proof: in the first we establish the hierarchies framework; in the second we derive a bound on the range of the $\U$-bootstrap process in the geometry of the $u$-norm; in the third we prove the crossing lemmas; in the fourth we use the crossing lemmas to bound the probability that a droplet is crossed with help from both above and below; in the fifth we deduce a bound on the probability that a droplet `grows' into a slightly larger droplet; in the sixth we apply the method of hierarchies to bound the probability that a `wide' droplet is internally spanned; in the seventh we assemble the different parts of the induction statement and prove Lemma~\ref{le:induction}; in the eighth we deduce Lemma~\ref{le:critunbal}, which is the bound for internally spanned critical droplets; and in the ninth and final subsection we complete the proof of Theorem~\ref{th:lowerunbal}.

\subsection{Hierarchies}\label{se:hierarchies}

The use of hierarchies to control the formation of critical droplets was introduced by Holroyd in~\cite{Hol} and has since developed into a standard technique in the study of bootstrap percolation, see for example~\cite{BBDCM,BBM3D,DCvE,DCH,GravHolMor}. In this subsection we recall some of the standard definitions and lemmas, making only minor adaptations along the way to suit the general model. We are relatively brief with the details, referring the reader instead to~\cite{Hol}, and the more recent refinements in~\cite{DCvE,GravHolMor}, for a more extensive introduction to the method.

The key result of this subsection is Lemma~\ref{le:boundoverH}, which gives an upper bound for the probability that a droplet $D$ is internally spanned in terms of the family of hierarchies of $D$.

Given a directed graph $G$ and a vertex $v\in V(G)$, we write $N_G^\to(v)$ for the set of out-neighbours of $v$ in $G$. Moreover, given a rooted tree, we consider it as a directed graph by orienting all of the edges so that they are pointing away from the root.

\begin{definition}\label{de:hierarchy}
Let $D$ be an $\stabu$-droplet. A \emph{hierarchy} $\hier$ for $D$ is an ordered pair $\hier=(G_\hier,D_\hier)$, where $G_\hier$ is a rooted tree, and $D_\hier$ is a function that assigns to each vertex of $G_\hier$ an $\stabu$-droplet, such that the following conditions are satisfied:
\begin{enumerate}
\item the root vertex corresponds to $D$, i.e., $D_\hier(\root)=D$;\smallskip
\item each vertex has out-degree at most 2;\smallskip
\item if $v \in N_{G_\hier}^\to(u)$ then $D_\hier(v) \subset D_\hier(u)$;\smallskip
\item if $N_{G_\hier}^\to(u) = \{v,v'\}$ then $\< D_\hier(v)\cup D_\hier(v') \> = \{ D_\hier(u) \}$.
\end{enumerate}
\end{definition}
Condition (iv) is equivalent to the statement that $D_\hier(v)$ and $D_\hier(v')$ are strongly connected to each other, and that $D_\hier(u)$ is the smallest droplet containing their union. We shall usually abbreviate $D_\hier(u)$ to $D_u$. 

The next definition controls the absolute and relative sizes of the droplets corresponding to vertices of $G_\hier$, which in turn allows us to control the number of hierarchies. In order to limit the number of hierarchies as much as possible, we choose the step size to be as large as possible, subject to the condition that we can control the probability of each step. 

\begin{definition}\label{de:hier2}
Fix $\beta\in\N$. A hierarchy $\hier$ for an $\stabu$-droplet $D$ is \emph{good} if it satisfies the following conditions for each $u \in V(G_\hier)$:
\begin{enumerate}
\setcounter{enumi}{4}
\item $u$ is a leaf if and only if $w(D_u) \le p^{-\beta(1-2\eta)-\eta}$;\smallskip
\item if $N_{G_\hier}^\to(u)=\{v\}$ and $|N_{G_\hier}^\to(v)| = 1$ then
\[
\xi \cdot p^{-\beta(1-2\eta)-\eta} \leq w(D_u) - w(D_v) \leq 2\xi \cdot p^{-\beta(1-2\eta)-\eta};
\]
\item if $N_{G_\hier}^\to(u)=\{v\}$ and $|N_{G_\hier}^\to(v)| \ne 1$ then $w(D_u) - w(D_v) \leq 2\xi \cdot p^{-\beta(1-2\eta)-\eta}$;\smallskip
\item if $N_{G_\hier}^\to(u)=\{v,v'\}$ then $w(D_u)-w(D_v) \geq \xi \cdot p^{-\beta(1-2\eta)-\eta}$.
\end{enumerate}
\end{definition}

Next we relate the abstract family of good hierarchies defined above to the initial set $A$ of infected sites and to the $\U$-bootstrap process. Given nested $\stabu$-droplets $D \subset D'$, we write $\Delta(D,D')$ for the event that $D'$ is internally spanned by $D$ together with the infected sites in $D' \setminus D$. That is,
\begin{equation}\label{def:deltaDD}
\Delta(D,D') := \big\{ D' \in \< D \cup (D' \cap A) \> \big\}. 
\end{equation}
The final two conditions below ensure that a good hierarchy for an internally spanned droplet $D$ accurately represents the growth of the initial sites $D\cap A$.

\begin{definition}\label{def:satisfied}
A hierarchy $\hier$ for an $\stabu$-droplet $D$ is \emph{satisfied} by $A$ if the following events all occur \emph{disjointly}:
\begin{enumerate}
\setcounter{enumi}{8}
\item if $v$ is a leaf then $D_v$ is internally spanned by $A$;\smallskip
\item if $N_{G_\hier}^\to(u) = \{v\}$ then $\Delta(D_v,D_u)$ occurs.
\end{enumerate}
\end{definition}

Having established all of the properties of hierarchies that we need, we now show that there exists a good and satisfied hierarchy for every internally spanned droplet. The proof is almost identical to Propositions~31 and~33 of~\cite{Hol}, which deal with the 2-neighbour setting, except that here we use the spanning algorithm in place of the rectangles process. We are therefore rather brief with the details.


\begin{lemma}\label{le:goodsat}
Let $D$ be an $\stabu$-droplet internally spanned by $A$. Then there exists a good and satisfied hierarchy for $D$.
\end{lemma}

\begin{proof}
In order to prove the lemma we consider a suitable `contraction' of the tree given by the spanning algorithm. To that end, let $\D = \< D \cap A \>$, and note that $D \in \D$ by Lemma~\ref{le:span}, since $D$ is internally spanned. The proof will be by induction on $w(D)$, so note first that if $w(D) \leq p^{-\beta(1-2\eta)-\eta}$ then we may take $V(G_\hier) = \{\root\}$. 

For the induction step, first we claim that there exists a pair of sequences, 
\[
D \cap A \supset K_0 \supset K_1 \supset \cdots \supset K_m \quad \text{and} \quad D = D_0 \supset D_1 \supset \cdots \supset D_m,
\]
such that $|K_m| = 1$ and such that for every $1\leq i\leq m$,
\[
D_i = D_U([K_i]) \quad \text{and} \quad [K_i] \cup [K_{i-1} \setminus K_i] \text{ is strongly connected}.
\]
To construct these sequences, the idea is to run the spanning algorithm backwards, choosing at each step the larger of the two droplets. We make this idea precise using Lemma~\ref{le:penultimate}. Indeed, since $D \in \< D \cap A \>$, there exists a set $K_0 \subset D \cap A$ such that $[K_0]$ is strongly connected and $D = D_U([K_0])$. Now, given $K_{i-1}$ such that $[K_{i-1}]$ is strongly connected, Lemma~\ref{le:penultimate} gives a (non-trivial) partition $K_i \cup K_i'$ of $K_{i-1}$ such that the sets $[K_i]$, $[K_i']$ and $[K_i] \cup [K_i']$ are each strongly connected. Set $D_i = D_U([K_i])$ and $D_i' = D_U([K_i'])$, where $w(D_i) \geq w(D_i')$.

Now, let $s \geq 1$ be minimal such that either
\[
w(D_s) \leq p^{-\beta(1-2\eta)-\eta} \qquad \text{or} \qquad w(D) - w(D_s) \geq \xi \cdot p^{-\beta(1-2\eta)-\eta},
\]
and attach a vertex $u$ corresponding to $D_s$ to the root. First, if $w(D_s) \leq p^{-\beta(1-2\eta)-\eta}$ and $w(D) - w(D_s) \leq 2\xi \cdot p^{-\beta(1-2\eta)-\eta}$, then our construction of $\hier$ is complete, since $D_s$ is internally spanned (by $K_s$), the event $\Delta(D_s,D)$ occurs (since $D$ is internally spanned), and the events $\ispan(D_s)$ and $\Delta(D_s,D)$ depend on disjoint sets of sites.

Next, suppose that $\xi \cdot p^{-\beta(1-2\eta)-\eta} \leq w(D) - w(D_s) \leq 2\xi \cdot  p^{-\beta(1-2\eta)-\eta}$. In this case we use the induction hypothesis to construct a good and satisfied (by $K_s$) hierarchy $\hier'$ for $D_s$, and identify $u$ with the root of $\hier'$. Once again, our construction is complete, since the events required for $\hier'$ to be satisfied depend on sites inside $D_s$, whereas the event $\Delta(D_s,D)$ depends on sites outside $D_s$.

Finally, if $w(D) - w(D_s) \geq 2\xi \cdot p^{-\beta(1-2\eta)-\eta}$ then, by the minimality of $s$, we have
\[
w(D_{s-1}) - w(D'_s) \geq w(D_{s-1}) - w(D_s) \geq \xi \cdot p^{-\beta(1-2\eta)-\eta}.
\]
In this case we add a vertex $v$ between $u$ and the root, corresponding to $D_{s-1}$, and add another vertex $u'$ attached to $v$, corresponding to $D_s'$. Now, using the induction hypothesis, we construct good and satisfied (by $K_s$ and $K_{s-1} \setminus K_s$ respectively) hierarchies $\hier'$ and $\hier''$ for $D_s$ and $D_s'$, and identify $u$ and $u'$ with the roots of $\hier'$ and $\hier''$. Since the events that $\hier'$ and $\hier''$ are satisfied are witnessed by sites in $K_s$ and $K_{s-1} \setminus K_s$ respectively, and the event $\Delta(D_{s-1},D)$ depends on sites outside $D_{s-1}$, the hierarchy thus constructed satisfies conditions (i)--(x), as required.
\end{proof}

\begin{remark}
We emphasize that the existence of a good and satisfied hierarchy for $D$ does \emph{not} imply that $D$ is internally spanned, since the intersection of the events $\ispan(D_v)$ and $\Delta(D_v,D_u)$ does not imply that $D_u$ is internally spanned, and since we do not insist that $[( D_v \cup D_{v'}) \cap A]$ is strongly connected whenever $N_{G_\hier}^\to(u) = \{ v, v' \}$. 
\end{remark}

We will use the following fundamental bound on the probability that a droplet is internally spanned (cf.~\cite[Section~10]{Hol} or~\cite[Lemma~20]{BBM3D}) in the proof of Lemma~\ref{le:induction}, and then again to prove~\eqref{eq:critunbal} for type $(L)$ critical droplets.

Let us write $\hier_D$ for the set of all good hierarchies for $D$, and $L(\hier)$ for the set of leaves of $G_\hier$. We will write $\sum_{u \to v}$ and $\prod_{u \to v}$ for the sum and product (respectively) over all pairs $\{u,v\} \subset V(G_\hier)$ such that $N_{G_\hier}^\to(u) = \{v\}$.

\begin{lemma}\label{le:boundoverH}
Let $D$ be an $\stabu$-droplet. Then
\begin{equation}\label{eq:boundoverH}
\Pr_p\big(\ispan(D)\big) \le \sum_{\hier \in \hier_D} \bigg( \prod_{u \in L(\hier)} \Pr_p\big(\ispan(D_u)\big) \bigg)\bigg( \prod_{u \to v} \Pr_p\big(\Delta(D_v,D_u)\big) \bigg).
\end{equation}
\end{lemma}

\begin{proof}
If $D$ is internally spanned then by Lemma~\ref{le:goodsat} there exists a good and satisfied hierarchy for $D$. Taking the union bound over good hierarchies, and noting that for a fixed good hierarchy $\hier$ the events $\ispan(D_u)$ (for $u\in L(\hier)$) and $\Delta(D_v,D_u)$ (for pairs $(u,v)$ such that $N_{G_\hier}^\to(u) = \{v\}$) are increasing and occur disjointly, the result follows from the van den Berg--Kesten inequality (Lemma~\ref{le:bk}).
\end{proof}


In order to use Lemma~\ref{le:boundoverH} we must bound the various probabilities that appear on the right-hand side of~\eqref{eq:boundoverH}, and the number of good hierarchies for $D$. If $u$ is a leaf of $V(G_\hier)$ then (following Holroyd~\cite{Hol}) we call $D_u$ a \emph{seed} of $\hier$. A sufficiently strong bound on the probability that a seed is internally spanned will follow immediately from the induction hypothesis, but bounding the probability of the events $\Delta(D_v,D_u)$ will be more challenging (see below). To count the good hierarchies, we will partition the set $\hier_D$ according to the number of seeds. 

Let us denote by $b(\hier)$ the number of seeds in a hierarchy $\hier$, and by $\hier_D^b$ the set of all good hierarchies for $D$ that have exactly $b$ seeds.


\begin{lemma}\label{le:counthier}
Let $D$ be an $\stabu$-droplet, and let $\hier \in \hier_D$ be a good hierarchy for $D$ with at least two vertices. Then
\begin{equation}\label{eq:vertices:hierarchy}
\big|V(G_\hier)\big| \le  4\xi^{-1} \cdot b(\hier) \cdot w(D) \cdot p^{\beta(1-2\eta)+\eta}.
\end{equation}
In particular, if $w(D) + h(D) \le p^{-O(1)}$, then
$$\big|\hier_D^b\big| \leq \exp\bigg[ O\bigg( \xi^{-1} b \cdot w(D) \cdot p^{\beta(1-2\eta)+\eta} \log\frac{1}{p} \bigg) \bigg]$$
for every $b \in \N$.
\end{lemma}

\begin{proof}
Let $\hier \in \hier_D$ be a good hierarchy with at least two vertices, and observe that
$$\big|V(G_\hier)\big| \le b(\hier) \big( d(\hier) + 1 \big) \le 2 \cdot b(\hier) \cdot d(\hier),$$
where $d(\hier)$ is the depth of the tree $G_\hier$ (i.e., the maximum length of a path from the root to a leaf in $G_\hier$), so $d(\hier) \ge 1$. We claim that 
$$d(\hier) \le 2 \xi^{-1} \cdot w(D) \cdot p^{\beta(1-2\eta)+\eta}.$$
Indeed, this follows when $d(\hier) = 1$ because $D$ is not a seed of $\hier$, and hence $w(D) > p^{-\beta(1-2\eta)-\eta}$. For $d(\hier) \ge 2$ it follows from the fact that every two steps up $G_\hier$, the width of the corresponding droplet increases by at least $\xi \cdot p^{-\beta(1-2\eta)-\eta}$. We therefore have~\eqref{eq:vertices:hierarchy} for every $\hier \in \hier_D$ with at least two vertices.


Now, the number of choices for the tree $G_\hier$ is at most $e^{O(N)}$, where $N$ is our bound on $|V(G_\hier)|$. Moreover, for each $u \in V(G_\hier)$, there are at most $p^{-O(1)}$ possible droplets $D_u$. Hence
$$\big|\hier_D^b\big| \leq \exp\bigg[ O\bigg( \xi^{-1} b \cdot w(D) \cdot p^{\beta(1-2\eta)+\eta} \log\frac{1}{p} \bigg) \bigg],$$
as claimed.
\end{proof}

We will need the following lower bound on the size of a seed, which is an easy consequence of Lemma~\ref{lem:subadditivity} and Definition~\ref{de:hier2}.

\begin{lemma}\label{obs:all:seeds:big}
If $\hier$ is a good hierarchy with at least two vertices, then 
$$w(D_v) \ge \frac{\xi}{2} \cdot p^{-\beta(1-2\eta)-\eta}$$
for every $v \in L(\hier)$. 
\end{lemma}

\begin{proof}
Let $v \in L(\hier)$, and observe that (since $\hier$ has at least two vertices) there exists $u \in V(G_\hier)$ with $v \in N_{G_\hier}^\to(u)$. Since $u \notin L(\hier)$ (and $\hier$ is good), we have $w(D_u) \ge p^{-\beta(1-2\eta)-\eta}$. Now, if $N_{G_\hier}^\to(u) = \{v\}$, then, by Definition~\ref{de:hier2} (vi) and (vii), 
$$w(D_v) = w(D_u) - \big( w(D_u) - w(D_v) \big) \ge \big( 1 - 2\xi \big) p^{-\beta(1-2\eta)-\eta}.$$
On the other hand, if $N_{G_\hier}^\to(u) = \{v,v'\}$ then, by Lemma~\ref{lem:subadditivity} (applied with $u = w^\perp$), 
$$w(D_v) \ge w(D_u) - w(D_{v'}) - O(1) \ge \xi \cdot p^{-\beta(1-2\eta)-\eta} - O(1).$$
Indeed, by Definition~\ref{de:hierarchy} (iv) the droplets $D_v$ and $D_{v'}$ are strongly connected to each other, so the first inequality follows from Lemma~\ref{lem:subadditivity}, while the second follows from Definition~\ref{de:hier2} (viii). Since $p$ is sufficiently small, the claimed bound follows.
\end{proof}

As mentioned above, bounding the probability of the event $\Delta(D_v,D_u)$ will be considerably more difficult (see Sections~\ref{se:crossings}--\ref{sec:growing:droplets}); moreover, we will only be able to prove a non-trivial bound when $|N_{G_\hier}^\to(v)| = 1$. We will use the following lemma to deduce a sufficiently strong bound on $\Pr_p\big(\ispan(D)\big)$. Let us write $d_\hier(v) := |N_{G_\hier}^\to(v)|$ for the out-degree of a vertex $v \in V(G_\hier)$ in the graph $G_\hier$.

\begin{lemma}\label{lem:sumofwidths}
Let $D$ be an $\stabu$-droplet, and let $b \le \xi^2 \cdot p^{-\beta(1-2\eta) + \eta}$. Then
$$\sum_{u \in L(\hier)} w(D_u) + \sum_{\substack{u \to v \\ d_\hier(v) = 1}} \big( w(D_u) - w(D_v) \big) \ge \frac{w(D)}{10}$$
for every $\hier \in \hier^b_D$. 
\end{lemma}

\begin{proof}
Note that the inequality holds trivially if $\hier$ has only one vertex, so we may assume that $\hier$ has at least two vertices. We first claim that 
\begin{equation}\label{eq:subadditivity:app}
\sum_{u \in L(\hier)} w(D_u) + \sum_{u \to v} \big( w(D_u) - w(D_v) \big) \ge w(D) - O\big( |V(G_\hier)| \big).
\end{equation}
Indeed, this follows by combining Definition~\ref{de:hierarchy} with the geometric inequality
$$w\big( D_U(D_1\cup D_2) \big) \le w(D_1) + w(D_2) + O(1),$$
which, by Lemma~\ref{lem:subadditivity}, holds for any pair of $\stabu$-droplets $D_1$ and $D_2$ that are strongly connected to each other. Next, observe that, by Lemma~\ref{le:counthier}, we have 
\begin{equation}\label{eq:number:of:vertices:bound}
\big|V(G_\hier)\big| \le 4\xi^{-1} \cdot b \cdot w(D) \cdot p^{\beta(1 - 2\eta) + \eta} \le 4\xi \cdot w(D).
\end{equation}
We next claim that
\begin{equation}\label{eq:sum:of:bad:pairs}
\sum_{\substack{u \to v \\ d_\hier(v) \ne 1}} \big( w(D_u) - w(D_v) \big) \le 4\xi \cdot b \cdot p^{-\beta(1-2\eta)-\eta} \le 8 \sum_{u \in L(\hier)} w(D_u).
\end{equation}
To prove the first inequality in~\eqref{eq:sum:of:bad:pairs}, recall that $w(D_u) - w(D_v) \le 2\xi \cdot p^{-\beta(1-2\eta)-\eta}$ for each pair $u,v$ in the sum, by Definition~\ref{de:hier2}, and observe that there are at most $2 \cdot b(\hier)$ pairs $u,v \in V(G_\hier)$ such that $N_{G_\hier}^\to(u) = \{v\}$ and $|N_{G_\hier}^\to(v)| \ne 1$. Indeed, $v$ must be either a leaf or a vertex with out-degree~$2$, and there are $b(\hier)$ seeds and $b(\hier) - 1$ vertices with out-degree~$2$. To prove the second inequality in~\eqref{eq:sum:of:bad:pairs}, recall that $w(D_u) \ge (\xi / 2) \cdot p^{-\beta(1-2\eta)-\eta}$ for every $u \in L(\hier)$, by Lemma~\ref{obs:all:seeds:big}.  

Now, combining~\eqref{eq:subadditivity:app},~\eqref{eq:number:of:vertices:bound} and~\eqref{eq:sum:of:bad:pairs}, it follows that
\begin{align*}
9 \sum_{u \in L(\hier)} w(D_u) + \sum_{\substack{u \to v \\ d_\hier(v) = 1}} \big( w(D_u) - w(D_v) \big) & \, \ge \sum_{u \in L(\hier)} w(D_u) + \sum_{u \to v} \big( w(D_u) - w(D_v) \big)\\
& \, \ge \big( 1 - O(\xi) \big) w(D).
\end{align*}
Since the implicit constant factor depends only on $\U$, and the constant $\xi$ was chosen to be sufficiently small, the claimed bound follows.
\end{proof}


\subsection{The range of unbalanced models with drift}\label{se:range}

In this short section we assume that $\U$ is an unbalanced model with drift and we use the results about $u$-icebergs from Section~\ref{se:iceberg} to prove a bound (see Lemma~\ref{le:rhobound}) on the range of the $\U$-bootstrap process helped by a half-plane $\H_u$.

Recall from Section~\ref{se:iceberg} that if $\alpha^-(u^*) = \infty$ then we choose $u_0 \in S^1$ to the left of and sufficiently close to $u^*$, so in particular $\alpha^-(u) = \infty$ for every $u \in [u^*,u_0]$. We also choose $u_1 \in (u^*,u_0)$. Similarly, if $\alpha^+(-u^*) = \infty$ then we choose a corresponding $u'_0 \in S^1$ to the right of and sufficiently close to $-u^*$, and a corresponding $u_1'$. Set
\[
\stabu^+ :=
\begin{cases}
  \hfill [u^*,u_1] \hfill & \text{if $\alpha^-(u^*) = \infty$} \\
  \hfill \{u^*\} \hfill & \text{otherwise} \\
\end{cases}
\quad \text{and} \quad \stabu^- :=
\begin{cases}
  \hfill [u_1',-u^*] \hfill & \text{if $\alpha^+(-u^*) = \infty$} \\
  \hfill \{-u^*\} \hfill & \text{otherwise,} \\
\end{cases}
\]
and set 
\[
\stabu^\pm := \big( \stabu^+\cup\stabu^- \big) \cap \Q_1 \qquad \text{and} \qquad \stabu' := \stabu^\pm \cup \{u^l,u^r\}.
\]
Recall also that for each $u \in \stabu^+$, we defined $\sigma(u)$ to be the angle between $u$ and $u^*$, and similarly for each $u\in\stabu^-$. Recall from Section~\ref{sec:measuring} that $w^\perp$ denotes one of the two unit vectors that are orthogonal to~$u^*$. We define a norm $\|\cdot\|_u$ on $\R^2$ as follows.

\begin{definition}\label{def:unorm}
For each $u\in\stabu'$, define
\begin{equation}\label{eq:unorm}
\|x\|_u :=
\begin{cases}
\hfill |\< x,u^* \> | + \sigma(u) |\< x,w^\perp \>| \hfill & \text{if } u\in\stabu^\pm\setminus\{u^*,-u^*\}, \\
\hfill \|x\| \hfill & \text{if } u \in \{u^*,-u^*,u^l,u^r\},
\end{cases}
\end{equation}
where, as always, the unadorned norm $\|\cdot\|$ denotes the Euclidean norm on $\R^2$.
\end{definition}

We record for later use the inequalities\footnote{If $u \in \{u^*,-u^*,u^l,u^r\}$ then both inequalities are trivial. If $u\in\stabu^\pm\setminus\{u^*,-u^*\}$ then note that the left-hand side is at most $\cos \sigma \cdot |\<x,u^*\>| + \sin \sigma \cdot |\<x,w^\perp\>|$, which implies the first inequality, and that $\sigma(u) < 1$, since $u_0$ was chosen sufficiently close to $u^*$, which implies the second.}
\begin{equation}\label{eq:normbound}
|\<x,u\>| \leq \|x\|_u \leq 2 \cdot \|x\|,
\end{equation}
which hold for every $u \in \S_U$ and $x \in \R^2$. Let $\rho\colon \stabu'\times\N\to\R$ be the function
\begin{equation}\label{eq:rhoweak}
\rho(u,\gamma) := \sup \Big\{ \big\|y-Y\big\|_u : |Y|=\gamma-1, \, y \in \big[ \H_u \cup Y \big] \setminus \H_u \Big\}.
\end{equation} 
The key property that we need for the vertical crossings lemma, and the main result of this section, is the following bound on $\rho(u,\gamma)$, which is uniform in $u$.

\begin{lemma}\label{le:rhobound}
Let $u \in \stabu'$ and $\gamma \in \N$, with $\gamma \leq \alphabar(u)$. Then $\rho(u,\gamma) \le f(\gamma)$, for some function $f(\cdot)$ depending only on\/ $\U$, the set $\S_U$, and the directions $u_0$ and $u_1$.
\end{lemma}


\begin{proof}
If $u \in \S_U$ then $\rho(u,\gamma) < \infty$ follows from the proof of Lemma~\ref{le:rhofinite}. Indeed, we proved there that $\rho_i(u)$ is finite for all $u \in \Q_1$ and $i \in \N$ such that $i < \alphabar(u) \le \infty$, and it follows from the definition that $\rho_i(u)$ depends only on $\U$, $u$ and $i$. Since $|\S_U| = 4 < \infty$, we may therefore assume that $u \in \stabu^\pm\setminus\{u^*,-u^*\}$, and hence, by symmetry, that $u\in\stabu^+\setminus\{u^*\}$. Note that this implies that $\alpha^-(u^*)=\infty$. 

Let $Y \subset\Z^2$ be a set of size $\gamma-1$, and let $\W$ be a $u$-iceberg cover of $Y$. The set
\[
\H_u \cup \bigcup_{W \in \W} W
\]
contains $Y$ and is closed, by Lemma~\ref{lem:iceberg:covers:are:closed}. Hence, if $y \in \big[\H_u\cup Y\big] \setminus \H_u$, then $y \in W$ for some $W \in \W$. Recall that $W$ is either an $\stabu$-droplet, or a $u$-iceberg that is $u$-iceberg covered by $W \cap Y$. Thus, by Lemma~\ref{le:extremaliceberg} and Definition~\ref{def:unorm}, there exists $x \in W \cap Y$ such that $\| x - y \|_u = O(\gamma)$,  where the implicit constant depends only on $\U$ and the fixed directions $u^*$, $u_0$ and $u_1$, as required.
\end{proof}

\subsection{Crossing lemmas}\label{se:crossings}

In this subsection we will prove two key lemmas, which bound the probabilities of certain `crossing' events (see Definition~\ref{def:ucrossed}, below), and which will allow us to bound (in Section~\ref{sec:tall:droplets}) the probability that an $\S_U$-droplet is `vertically crossed', and (in Section~\ref{sec:growing:droplets}) the probability of the events $\Delta(D,D')$ that appear in Lemma~\ref{le:boundoverH}. We will refer to these two lemmas as the `horizontal crossings lemma' (Lemma~\ref{le:horizcrossings}), which allows us to control growth in directions $u^l$ and $u^r$, and the `vertical crossings lemma' (Lemma~\ref{le:vertcrossings}), which deals with growth in directions close to $\pm u^*$. These bounds are the key step in the proof of Theorem~\ref{th:lowerunbal}, and this subsection is the most technical of the paper.

Since there is significant overlap between the proofs for horizontal and vertical crossings, it will be convenient to work in the following (slightly) more general framework.

\begin{definition}\label{de:ustrip}
Let $u \in \stabu'$. A finite set is a \emph{$u$-strip} if it is a $\T$-droplet, where $\T = \{u,-u,v,-v\}$ and either
\begin{itemize}
\item $u\in\{u^l,u^r\}$ and $v=u^*$ (a \emph{horizontal strip}), or
\item $u \in \stabu^\pm = \stabu' \setminus \{u^l,u^r\}$ and $v=w^\perp$ (a \emph{vertical strip}).
\end{itemize}
\end{definition}

Although it is convenient to define $u$-strips in terms of $\T$-droplets, we stress again that \emph{all} sets described in this section as `droplets' without reference to a set $\T$ are assumed to be $\stabu$-droplets.

Recall from~\eqref{eq:kappa} that $\kappa = 2\rho + 3\nu$, 
that we denote by $G_\kappa$ the graph with vertex set $\Z^2$ and edge set $E = \big\{ xy : \|x-y\| \leq \kappa \big\}$, and that a strongly connected component is defined to be a component in this graph. Recall also that the $u$-projection $\pi(K,u)$ of a finite set $K \subset \Z^2$ was defined in~\eqref{eq:uproj} by
\[
\pi(K,u) = \max\big\{ \<x-y,u\> \,:\, x,y\in K \big\},
\]
and that if $D$ is a $\T$-droplet and $u \in \T$, then the $u$-side $\partial(D,u)$ of $D$ was defined in~\eqref{def:uside} to be the set $D \cap \ell_u(i)$, where $i$ is maximal so that this set is non-empty. 

\begin{definition}\label{def:ucrossed}
Let $u \in \stabu'$, let $S$ be a $u$-strip, and let $x \in \partial(S,-u)$. We say that $S$ is \emph{$u$-crossed} if there exists a strongly connected set in $[ \H_u(x) \cup (S \cap A) ]$ that intersects both $\H_u(x)$ and $\partial(S,u)$. 
\end{definition}

Note that the half-plane $\H_u(x)$ does not depend on the choice of~$x \in \partial(S,-u)$. Unless the precise position of the $u$-strip is important, we will usually assume that the $(-u)$-side of the $u$-strip is a subset of $\ell_u$. 

Before continuing with the results of this subsection, we give a more complete account of the relationships between the different constants of this section than that given in~\eqref{eq:constants1}. We mentioned that, during the course of the inductive proof of Lemma~\ref{le:induction}, two sequences of constants would be defined, in addition to the constants $\delta(\beta)$ already introduced in Definition~\ref{de:ih}. These sequences are $\delta'(2),\dots,\delta'(2\alpha+1)$, which will appear in the statements of Lemmas~\ref{le:horizcrossings} and~\ref{le:vertcrossings}, and $\weak(2),\dots,\weak(2\alpha+1)$, which will appear in Definition~\ref{de:weakconn}. These constants will be chosen to have the following relative sizes. First, for each $2\leq\beta\leq 2\alpha$,
\begin{equation}\label{eq:constants2}
1 \gg \delta(\beta) \gg \frac{1}{\weak(\beta)} \gg \delta'(\beta) \gg \delta(\beta+1) > 0,
\end{equation}
and second,
\begin{equation}\label{eq:constants3}
\delta(2\alpha+1) \gg \frac{1}{\weak(2\alpha+1)} \gg \delta'(2\alpha+1) \gg \xi \gg \delta \gg \epsilon > 0.
\end{equation}
We emphasize again that these statements mean that the constants are chosen from left to right, and that each is chosen to be sufficiently small depending on all previously chosen constants. Note that these two sets of relations subsume those in~\eqref{eq:constants1}.

The main results of this subsection are Lemmas~\ref{le:horizcrossings} and~\ref{le:vertcrossings}, below. One may think of the lemmas as exchanging bounds on the probability that a droplet is internally spanned for bounds on the probability that similarly sized $u$-strips are $u$-crossed, for some $u$. (It may be helpful, therefore, to think of the $\delta'(\beta)$ as being to crossed $u$-strips as the $\delta(\beta)$ are to internally spanned $\stabu$-droplets.) 

The first lemma bounds the probability of crossing a $u$-strip when $u\in\{u^l,u^r\}$.

\begin{lemma}\label{le:horizcrossings}
Let $S$ be a $u$-strip, where $u\in\{u^l,u^r\}$.
\begin{itemize}
\item[(i)] Let $1 \le \beta_1 \le \beta_2 \le \alpha$, and suppose that $\ih(\beta_1,\beta_2)$ holds. If 
\[
\xi^{-1} \le \pi(S,u) \leq p^{-\beta_1(1-2\eta)-\eta} \quad\text{and}\quad h(S) \leq p^{-\alpha(1-2\eta)-\eta},
\]
then $S$ is $u$-crossed with probability at most $p^{\delta' \pi(S,u)}$, where $\delta'=\delta'(\beta_1+\beta_2)$.\smallskip
\item[(ii)] Suppose that $\ih(\alpha,\alpha+1)$ holds. If
\[
\xi^{-1} \le \pi(S,u)\leq p^{-\alpha(1-2\eta)-\eta} \quad\text{and}\quad h(S)\leq\frac{\xi}{p^\alpha}\log\frac{1}{p},
\]
then $S$ is $u$-crossed with probability at most
\[
\pi(S,u) \cdot \exp\Big(-p^{O(\xi)}\cdot \pi(S,u) \Big),
\]
where the implicit constant depends on $\weak(2\alpha+1)$.
\end{itemize}
\end{lemma}


The second of our two main crossing lemmas deals with crossings $u$-strips when $u$ is close to $\pm u^*$. Recall that $u^*$ has difficulty at least $\alpha+1$, and therefore either $\alphabar(u^*) \geq \alpha+1$ or $\alpha^-(u^*) = \infty$. The behaviour of the $\U$-bootstrap process differs markedly depending on which of these two cases we are in. Note that, while the lemma is stated only for $u \in \stabu^+$, it is plain by symmetry that a similar statement holds for $u\in\stabu^-$. Moreover, for simplicity we assume that $p$ is chosen so that if $\sigma(u) = p^{1-\eta}$, then $u \in \Q_1$.

\begin{lemma}\label{le:vertcrossings}
Let $u \in \stabu^+ \cap \Q_1$ be such that either
\begin{align*}
u = u^* \quad &\text{and} \quad \alphabar(u^*) \geq \alpha + 1, \quad \text{or} \\
\sigma(u) = p^{1-\eta} \quad &\text{and} \quad \alpha^-(u^*)=\infty.
\end{align*}
Let $\beta_1,\beta_2 \in \N$ with $\beta_1 + \beta_2 \le 2\alpha + 1$, set $\beta := \min\{\beta_1,\beta_2\}$, and suppose that $\ih(\beta_1,\beta_2)$ holds. If $S$ is a $u$-strip with 
$$w(S) \le p^{-(\alpha + 1)(1-2\eta)-\eta} \qquad \text{and} \qquad \xi^{-1} \le \pi(S,u) \le p^{-\beta(1-2\eta)-\eta},$$ 
then $S$ is $u$-crossed with probability at most $p^{\delta' \pi(S,u)}$, where $\delta'=\delta'(\beta_1+\beta_2)$.
\end{lemma}

Observe that if $u^*$ is not a drift direction (that is, if $\alpha^-(u^*) < \infty$) then the lemma says it is unlikely that a $u^*$-strip of an appropriate size is $u^*$-crossed -- this is what one would expect. If $u^*$ is a drift direction, on the other hand, then instead the lemma is stated in terms of crossing $u$-strips, where $\sigma(u)=p^{1-\eta}$. Why might this be the natural direction in which to bound growth? Since $u^*$ is a drift direction, $\alpha^+(u^*)$ may be as small as $1$, and therefore one would expect a triangle of sites of slope $p$ to form on the $u^*$-side of the droplet, similarly to the set $T$ in Figure~\ref{fi:unbaldroplets}. By rotating $u^*$ through an angle of $p^{1-\eta}$, we are `giving away' more sites than one would expect to become infected, but not so many more that it adversely affects the bound. We expand on these remarks before the proof of the lemma.

The first step towards proving Lemmas~\ref{le:horizcrossings} and~\ref{le:vertcrossings} is a deterministic description of the structure of $S\cap A$ when $S$ is $u$-crossed, which is given by Lemma~\ref{le:generalpartitions}. We partition the $u$-strip into consecutive $u$-strips $S_1,\dots,S_m$ of constant $u$-projection, and we consider how the infection could spread from the $(-u)$-side of $S$ (and the adjacent half-plane $\H_u$) to the $u$ side of $S$. One of the key concepts we use will be that of being `$u$-weakly connected', defined as follows. 

\begin{definition}\label{de:weakconn}
Fix $\beta_1,\beta_2 \in \N$ with $\beta_1 + \beta_2 \le 2\alpha + 1$, and let $\weak=\weak(\beta_1+\beta_2)$. For each $u \in \stabu'$ and $\gamma \in \N$, we say that:
\begin{itemize}
\item[$(a)$] A set $Z \subset \Z^2$ is \emph{$u$-weakly connected} if it is connected in the graph $G_{u,\weak}$ with vertex set $\Z^2$ and edge set $E(G_{u,\weak}) = \big\{ xy : \|x-y\|_u \leq \weak \big\}$.\smallskip 
\item[$(b)$] A \emph{$u$-weak $\gamma$-cluster} is a set of $\gamma$ sites that is $u$-weakly connected.
\end{itemize}
\end{definition}

Note that we suppress the dependence on the pair $(\beta_1,\beta_2)$ in the definition of a $u$-weak $\gamma$-cluster; we trust that this will not cause any confusion. We remark that in what follows we will always take $\gamma \leq \alphabar(u)$, so if $Y$ is a $u$-weak $(\gamma-1)$-cluster, then (by Lemma~\ref{le:rhobound}) taking the closure of $\H_u \cup Y$ only causes `local' new infections, measured in the $u$-norm.

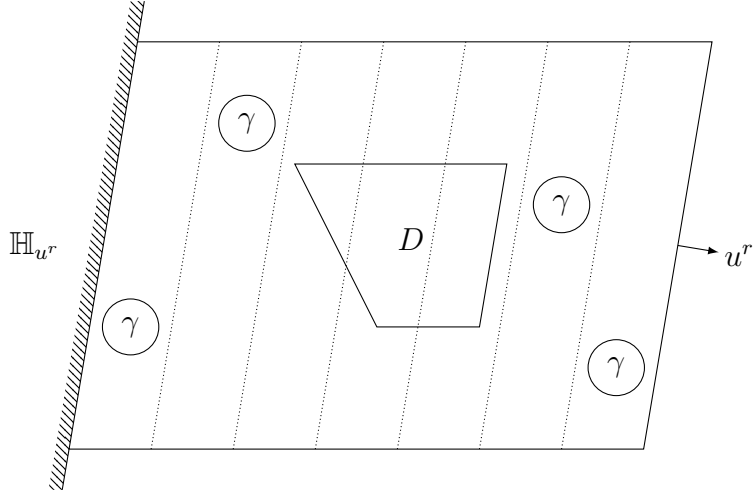
\begin{figure}[ht]
  \centering
  \begin{tikzpicture}[>=latex,scale=0.9]
    \draw (-0.1,-0.6) -- (0,0) (1,6) -- (1.1,6.6);
    \fill [pattern=north west lines] (-0.3,-0.6) -- (-0.1,-0.6) -- (1.1,6.6) -- (0.9,6.6) -- cycle;
    \draw (8.4,0) -- (0,0) -- (1,6) -- (9.4,6) -- cycle;
    \foreach \x in {1,2,...,6}
      \draw [densely dotted] (1.2*\x,0) -- (1.2*\x+1,6);
    \node [circle,draw] at (0.9,1.8) {$\gamma$};
    \node [circle,draw] at (2.6,4.8) {$\gamma$};
    \node [circle,draw] at (7.2,3.6) {$\gamma$};
    \node [circle,draw] at (8,1.2) {$\gamma$};
    \draw (3.3,4.2) -- (6.4,4.2) -- (6,1.8) -- (4.5,1.8) -- cycle;
    \node at (5,3.1) {$D$};
    \draw [->] (8.9,3) -- ++(0.6,-0.1);
    \path (8.9,3) ++(0.9,-0.15) node {$u^r$};
    \node at (-0.5,3) {$\H_{u^r}$};
    \path [name path=P1] (2.6,1.2) -- ++(1.8,-0.3);
    \path [name path=P2] (3.6,0) -- (4,2.4);
    \path [name intersections={of=P1 and P2,by=I}];
  \end{tikzpicture}
  \caption{A $u^r$-crossed $u^r$-strip $S$ together with a possible $(u^r,\gamma)$-partition for $S\cap A$ in which $a_1=a_2=a_4=a_5=1$ and $a_3=3$.}\label{fi:partitions}
\end{figure}

We can now define the deterministic structural property that we shall prove (in Lemma~\ref{le:generalpartitions}, below) is implied by the event that $S$ is $u$-crossed by $A$. The definition is illustrated in Figure~\ref{fi:partitions}.

\begin{definition}\label{de:partition}
Fix $\beta_1,\beta_2 \in \N$ with $\beta_1 + \beta_2 \le 2\alpha + 1$, and let $\weak=\weak(\beta_1+\beta_2)$. Let $u \in \stabu'$ and $\gamma \in \N$, and suppose that $S$ is a $u$-strip. Let $S_1 \cup \dots \cup S_{m+1}$ be a partition of $S$ into $u$-strips, with $S_i$ adjacent to $S_{i+1}$ for each $i \in [m]$,
\begin{equation}\label{eq:partition:strip:size}
3 \weak \gamma \leq \pi(S_i,u) = \pi(S_j,u) \leq 4 \weak \gamma
\end{equation}
for each $i,j \in [m]$, and $\pi(S_{m+1},u) < 4 \weak \gamma$.

A \emph{$(u,\gamma)$-partition} for $S \cap A$ is a sequence $(a_1,\ldots,a_k)$ of positive integers with $a_1 + \dots + a_k = m$, such that for each $1\leq j\leq k$, setting $t_j = a_1 + \cdots + a_j$, either
\begin{itemize} 
\item $a_j = 1$ and $S_{t_j} \cap A$ contains a $u$-weak $\gamma$-cluster, or \smallskip
\item there exists an $\stabu$-droplet $D$ internally spanned by $\big(S_{t_{j-1} + 1} \cup \cdots \cup S_{t_j}\big)\cap A$, where
\[
\max\big\{ w(D), h(D) \big\} \geq \frac{a_j \weak}{5}.
\]
\end{itemize}
\end{definition}

In the following lemma we need an upper bound on $\gamma$ when $\alphabar(u) = \infty$. For this purpose, recalling~\eqref{def:eta}, fix 
\begin{equation}\label{eq:lambda}
\lambda := 5\alpha/\eta = 50 \alpha^2
\end{equation}
throughout this section, and observe that 
\begin{equation}\label{eq:partrho}
\weak(2) \gg \rho(u,\gamma)
\end{equation}
for every $u \in \stabu'$ and $\gamma \in \N$ such that $\gamma \le \min\{ \alphabar(u), \lambda \}$. Indeed, $\rho(u,\gamma)$ is bounded above by a constant that depends only on $\,\U$ and the fixed directions $u^*$, $u_0$ and $u_1$, by Lemma~\ref{le:rhobound}, and $\weak(2)$ was chosen in~\eqref{eq:constants2} to be sufficiently large.

The following deterministic lemma, which says that every $u$-crossed strip has a $(u,\gamma)$-partition, is the key step in the proof of Lemmas~\ref{le:horizcrossings} and~\ref{le:vertcrossings}.

\begin{lemma}\label{le:generalpartitions}
Let $\beta_1,\beta_2 \in \N$ with $\beta_1 + \beta_2 \le 2\alpha + 1$, and set $\weak = \weak(\beta_1+\beta_2)$. Suppose that $u \in \stabu'$ and $\gamma \in \N$ satisfy $\gamma \leq \min\big\{ \alphabar(u),\lambda \big\}$, and let $S$ be a $u$-strip. If $S$ is $u$-crossed by $A$, then there exists a $(u,\gamma)$-partition for $S\cap A$.
\end{lemma}

Roughly speaking, the proof of the lemma is as follows. We shall show that if $S_1 \cap A$ does not contain a $u$-weak $\gamma$-cluster, then $S_1$ cannot itself be $u$-crossed. Since $S$ is $u$-crossed, this will allow us to deduce that there exists a droplet $D$ internally spanned by $S\cap A$ such that $D\cap S_1\neq\emptyset$, and moreover such that $D$ extends at least halfway across $S_1$. We call such a droplet $D$ a \emph{saver}. Letting $a_1$ be maximal such that $D \cap S_{a_1}\neq\emptyset$, the result follows by induction.

\begin{proof}[Proof of Lemma~\ref{le:generalpartitions}]
Let $S = S_1 \cup \dots \cup S_{m+1}$ be a partition of $S$ into $u$-strips as in Definition~\ref{de:partition}, so $S_i$ adjacent to $S_{i+1}$ for each $i \in [m]$,~\eqref{eq:partition:strip:size} holds for each $i,j \in [m]$, and $\pi(S_{m+1},u) < 4 \weak \gamma$. We prove the lemma by induction on $m$. If $m = 0$ there is nothing to prove, so let $m \geq 1$ and assume that the result holds for every smaller non-negative value of $m$. If $S_1 \cap A$ contains a $u$-weak $\gamma$-cluster then we are done, since we may set $a_1 = 1$ and observe that $S \setminus S_1$ is $u$-crossed by $A$. 

So assume that $S_1 \cap A$ does not contain a $u$-weak $\gamma$-cluster, let $Y_1,\dots,Y_s$ be the collection of $u$-weakly connected components in $S \cap A$ that are each also $u$-weakly connected to $\H_u$, and set 
\[
Y := Y_1 \cup \dots \cup Y_s \qquad \text{and} \qquad  Z := [\H_u \cup Y] \setminus \H_u. 
\]
We claim that $|Y_i| \leq \gamma-1$ for each $1\leq i\leq s$. Indeed, if $|Y_i| \geq \gamma$ then there exists a $u$-weak $\gamma$-cluster $Y' \subset Y_i$ such that $\| y - \H_u \|_u \le \weak \gamma$ for every $y \in Y'$. Recalling from~\eqref{eq:normbound} that $\< x, u \> \le \| x \|_u$ for every $x \in \Z^2$, and that $\pi(S_1,u) \ge 3\weak \gamma$, it follows that $Y' \subset S_1$. This contradicts our assumption that $S_1 \cap A$ does not contain a $u$-weak $\gamma$-cluster, and thus proves that $|Y_i| \leq \gamma-1$ for each $1\leq i\leq s$, as claimed.

We next claim that
\begin{equation}\label{eq:partZ}
\H_u \cup Z = [\H_u\cup Y_1] \cup \dots \cup [\H_u\cup Y_s]. 
\end{equation}
To prove this, let
\[
z_i\in[\H_u\cup Y_i]\setminus\H_u \qquad \text{and} \qquad z_j\in[\H_u\cup Y_j]\setminus\H_u,
\]
and note that $\|z_i-Y_i\|_u\leq\rho(u,\gamma)$ and $\|z_j-Y_j\|_u\leq\rho(u,\gamma)$, by the definition~\eqref{eq:rhoweak} of $\rho(u,\gamma)$. Hence
\[
\| z_i - z_j \| \geq \frac{\| z_i - z_j \|_u}{2} \geq \frac{\weak - 2\rho(u,\gamma)}{2} > \nu,
\]
where the first inequality follows from~\eqref{eq:normbound}, the second by the triangle inequality, and the third from~\eqref{eq:partrho}, since $\gamma \leq \min\{\alphabar(u),\lambda\}$. Therefore, the set
\[
[\H_u\cup Y_1]\cup\dots\cup[\H_u\cup Y_s]
\]
is closed (and contains $Y$), which proves \eqref{eq:partZ}. Note that it follows from the above argument that moreover
\begin{equation}\label{eq:partZ:moreover}
\| z - Y \|_u \le \rho(u,\gamma)
\end{equation} 
for every $z \in Z$. 
 
We are now ready to prove our key claim, which says that, under our assumption that $S_1$ does not contain a $u$-weak $\gamma$-cluster, there exists a droplet that is internally spanned by $S \cap A$, and has large intersection with $S_1$. 

\begin{claim}\label{cl:partitions}
There exists a droplet $D$ internally spanned by $S \cap A$ such that
\[
|\<D-\H_u,u\>| \le 2\weak\gamma \qquad \text{and} \qquad \max\big\{ w(D), h(D) \big\} \geq \frac{\weak}{5},
\]
where $\<D-\H_u,u\> = \min\big\{ \< x - y, u \> : x \in D, \, y \in \H_u \big\}$. 
\end{claim}

\begin{proof}[Proof of Claim~\ref{cl:partitions}]\let\qednow\qedsymbol\renewcommand{\qedsymbol}{}
The first step is to show that there exist $z \in \H_u \cup Z$ and $w \in [ S \cap A \setminus Y ]$ with
\begin{equation}\label{eq:partw}
\|w-z\|_u \leq 2\cdot\|w-z\| \leq 2\kappa.
\end{equation}
Note that the first inequality follows by~\eqref{eq:normbound}, so we just need to prove the second. To do so, recall that $S$ is $u$-crossed by $A$, which means that there exists a strongly connected component $L \subset [ \H_u \cup (S \cap A) ]$ that intersects both $\H_u$ and $\partial(S,u)$. Note that $\H_u \cup Z$ does not intersect $\partial(S,u)$, since~\eqref{eq:normbound},~\eqref{eq:partrho},~\eqref{eq:partZ:moreover} and $m \ge 1$ imply that
$$|\<z - \H_u, u \>| \le \|z - \H_u\|_u \le \rho(u,\gamma) + \weak \gamma < 3 \weak \gamma \le \pi(S,u)$$
for every $z \in Z$. Now, either $[ \H_u \cup (S \cap A) ] = [\H_u \cup Y] \cup [ S \cap A \setminus Y ]$, in which case there must exist $z \in [\H_u \cup Y] = \H_u \cup Z$ and $w \in [ S \cap A \setminus Y ]$ with $\|w - z\| \leq \kappa$, or $[ \H_u \cup (S \cap A) ] \ne [\H_u \cup Y] \cup [ S \cap A \setminus Y ]$, in which case there must exist $z \in \H_u \cup Z$ and $w \in [ S \cap A \setminus Y ]$ with $\|w - z\| \le \nu$. Since $\kappa \ge 3\nu$, in either case~\eqref{eq:partw} holds.


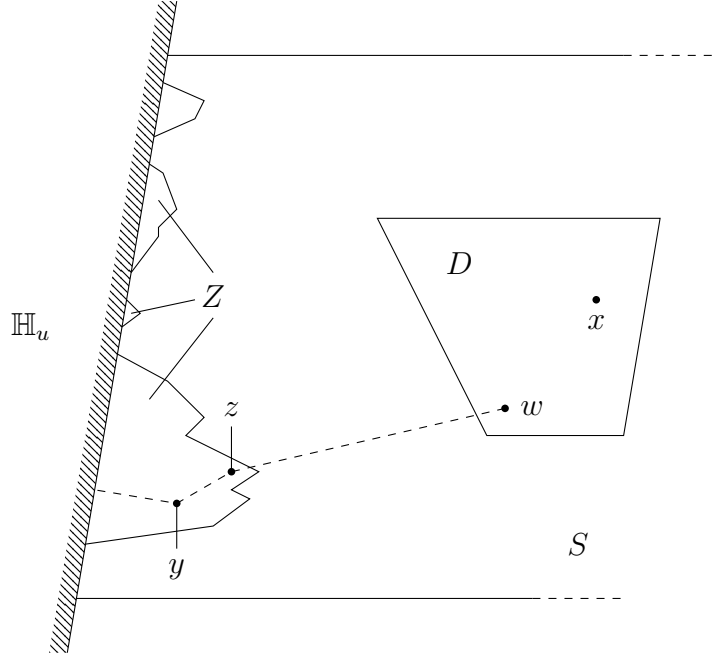
\begin{figure}[ht]
  \centering
  \begin{tikzpicture}[>=latex,scale=1.2]
    \draw (-0.1,-0.6) -- (0,0) (1,6) -- (1.1,6.6);
    \fill [pattern=north west lines] (-0.3,-0.6) -- (-0.1,-0.6) -- (1.1,6.6) -- (0.9,6.6) -- cycle;
    \draw (5,0) -- (0,0) -- (1,6) -- (6,6);
    \draw [dashed] (5,0) -- (6,0) (6,6) -- (7,6);
    \draw (3.3,4.2) -- (6.4,4.2) -- (6,1.8) -- (4.5,1.8) -- cycle;
    \draw (0.1,0.6) -- (1.5,0.8) -- (1.9,1.1) -- (1.7,1.2) -- (2,1.4) -- (1.2,1.8) -- (1.4,2) -- (1,2.4) -- (0.45,2.7);
    \draw (0.5,3) -- (0.7,3.15) -- (0.55,3.3);
    \draw (0.6,3.6) -- (0.9,4) -- (0.9,4.1) -- (1.1,4.3) -- (0.95,4.7) -- (0.8,4.8);
    \draw (0.85,5.1) -- (1.3,5.3) -- (1.4,5.5) -- (0.95,5.7);
    \node [inner sep=0pt,minimum size=0.1cm,circle,fill,label=below:$x$] at (5.7,3.3) {};
    \node [inner sep=0pt,minimum size=0.1cm,circle,fill,label=right:$w$] at (4.7,2.1) {};
    \node [inner sep=0pt,minimum size=0.1cm,circle,fill] at (1.7,1.4) {};
    \draw (1.7,1.4) -- ++(0,0.5) node [above] {$z$};
    \node [inner sep=0pt,minimum size=0.1cm,circle,fill] at (1.1,1.05) {};
    \draw (1.1,1.05) -- ++(0,-0.5) node [below] {$y$};
    \draw [dashed] (4.7,2.1) -- (1.7,1.4) -- (1.1,1.05) -- (0.2,1.2);
    \draw (0.8,2.2) -- (1.5,3.1) node [above] {$Z$};
    \draw (0.6,3.15) -- (1.3,3.3) (0.9,4.4) -- (1.5,3.6);
    \node at (4.2,3.7) {$D$};
    \node at (5.5,0.6) {$S$};
    \node at (-0.5,3) {$\H_u$};
  \end{tikzpicture}
  \caption{The situation in the proof of Claim~\ref{cl:partitions} is depicted assuming $z\in Z$. The size of the projection $\big|\<D-\H_u,u\>\big|$ is at most the total length of the dashed line in the $u$-norm.}
  \label{fi:partclaim}
\end{figure}

Now, let $\D$ be the output of the spanning algorithm with input $S \cap A \setminus Y$, and let $D \in \D$ be the droplet spanned by the strongly connected component of $[S\cap A\setminus Y]$ containing $w$. If $z \in \H_u$, then it follows by~\eqref{eq:partw} and the $u$-norm bound in~\eqref{eq:normbound} that
$$\big|\<D-\H_u,u\>\big| \leq |\< w - z, u \>| \leq \|w-z\|_u \leq 2\kappa \leq 2\weak \gamma.$$
On the other hand, if $z \in Z$ then $\|z - Y\|_u \leq \rho(u,\gamma) \ll \weak$, by~\eqref{eq:partrho} and~\eqref{eq:partZ:moreover}. Therefore, recalling that every $y \in Y$ is within distance at most $\weak \gamma$ of $\H_u$ in the $u$-norm, the triangle inequality and~\eqref{eq:normbound} gives
$$|\< D - \H_u, u \>| \leq \|w-z\|_u + \|z-Y\|_u + \weak \gamma \le 2\weak \gamma,$$
as required.

To bound the dimensions of $D$, let $x \in D \cap A \setminus Y$, and observe that
\[
\big\|x-(\H_u\cup Y)\big\|_u > \weak
\]
by the definition of $Y$. Using~\eqref{eq:partZ:moreover}, it follows that
\[
\big\| x - (\H_u\cup Z) \big\|_u > \weak - \rho(u,\gamma),
\]
and hence, by~\eqref{eq:normbound},
\[
\big\| x - (\H_u\cup Z) \big\| \geq \frac{\big\| x - (\H_u\cup Z) \big\|_u}{2}  \geq \frac{\weak-\rho(u,\gamma)}{2}.
\]
However, by~\eqref{eq:partw} and our choice of $w$, we also have
\[
\big\| w - (\H_u\cup Z) \big\| \leq \| w - z \| \leq \kappa.
\]
Since $x,w \in D$, it follows that 
\begin{align*}
\max\big\{ w(D), h(D) \big\} &\geq \frac{\|w-x\|}{2} \geq \frac{\big\|x-(\H_u\cup Z)\big\| - \big\|w-(\H_u\cup Z)\big\|}{2} \\
&\geq \frac{\weak-\rho(u,\gamma)-2\kappa}{4} \geq \frac{\weak}{5},
\end{align*}
by the triangle inequality and~\eqref{eq:partrho}, as required. \qednow
\end{proof}

To complete the proof of the lemma, simply set $a_1 = \max\{ i : D\cap S_i\neq\emptyset \}$, and observe that $S \setminus (S_1 \cup \cdots \cup S_{a_1})$ is $u$-crossed by $A$. It follows easily from Claim~\ref{cl:partitions}, our choice of $a_1$, and the fact that $\pi(S_i,u) \ge 3 \weak \gamma$ for every $i \in [m]$, that $\max\{ w(D), h(D) \} \geq a_1 \weak / 5$, as required.
\end{proof}

We next prove an upper bound, depending on $u \in \S_U'$ and on the size of $S$, on the probability that a $p$-random subset of a $u$-strip $S$ admits a $(u,\gamma)$-partition. In order to (somewhat) simplify the statement, given $u \in \S_U'$ and a $u$-strip $S$, let us write $g_u(S)$ for the maximum number of $u$-weak $\gamma$-clusters in a sub-strip $S' \subset S$ of $u$-projection at most $4\weak\gamma$ (that is, the number in a maximal such sub-strip). Note that  in this notation we suppress the dependence on $\gamma$ and on the pair $(\beta_1,\beta_2)$.

\begin{lemma}\label{le:partitionscalc}
Let $\beta_1,\beta_2 \in \N$ with $\beta_1 + \beta_2 \le 2\alpha + 1$, set $\beta:=\min\{\beta_1,\beta_2\}$ and $\weak = \weak(\beta_1+\beta_2)$, and assume that $\ih(\beta_1,\beta_2)$ holds. Let $u \in \stabu'$ and $\gamma \in \N$, with $\gamma \leq \min\big\{ \alphabar(u),\lambda \big\}$, and let $S$ be a $u$-strip with $|S| \le p^{-3\alpha}$, and
\begin{equation}\label{eq:Sconds}
\xi^{-1} \leq \pi(S,u)\leq p^{-\beta(1-2\eta)-\eta}.
\end{equation}
Then the probability that $S\cap A$ admits a $(u,\gamma)$-partition is at most
\begin{equation}\label{eq:calc}
\pi(S,u) \cdot \max_{0\,\leq \,j\,\leq \,m} \Big(1-\big(1-p^\gamma\big)^{g_u(S)}\Big)^{m-j} \big( \pi(S,u)^2 \cdot p^{2\alpha}\big)^j,
\end{equation}
where $m= \lfloor \pi(S,u) / 5\weak \gamma \rfloor$. 
\end{lemma}

\begin{proof}
We first deal with a technicality: the saver droplets need only be \emph{internally spanned} by the sites in $S\cap A$; they do not have to be \emph{contained} in $S$, and therefore their dimensions may be too large to use $\ih(\beta_1,\beta_2)$. Moreover, even if the savers \emph{are} contained in $S$, they may still have dimensions too large to use $\ih(\beta_1,\beta_2)$. However, neither of these is a problem, as we now show. Let $D$ be any saver droplet (so $D$ is internally spanned by $S\cap A$) such that either
\begin{equation}\label{eq:Dconds}
w(D) \geq p^{-\beta_1(1-2\eta)-\eta} \quad \text{or} \quad h(D) \geq p^{-\beta_2(1-2\eta)-\eta}.
\end{equation}
Then by Lemma~\ref{le:ALspan}, applied once with $u=u^*$ and again if necessary with $u=w^\perp$, there exists a droplet $D'\subset D$, also spanned by sites in $S \cap A$, such that either
\begin{align*}
&w(D')\leq p^{-\beta_1(1-2\eta)-\eta} \quad \text{and} \quad p^{-\beta_2(1-2\eta)-\eta}/3 \leq h(D') \leq p^{-\beta_2(1-2\eta)-\eta}, \\
\text{or} \quad &p^{-\beta_1(1-2\eta)-\eta}/3 \leq w(D')\leq p^{-\beta_1(1-2\eta)-\eta} \quad \text{and} \quad h(D') \leq p^{-\beta_2(1-2\eta)-\eta}.
\end{align*}
Therefore, by $\ih(\beta_1,\beta_2)$, we have $\Pr_p\big(\ispan(D')\big)\leq p^{\delta k/3}$, where $\delta=\delta(\beta_1+\beta_2)$ and
\[
k := \min\big\{ p^{-\beta_1(1-2\eta)-\eta}, p^{-\beta_2(1-2\eta)-\eta} \big\}.
\]
But $k \geq \pi(S,u)$, since $S$ satisfies \eqref{eq:Sconds}, and therefore
\[
\Pr_p\big(\ispan(D')\big) \leq p^{\delta \pi(S,u)/3}.
\]
Now, observe that there are at most $p^{-13\alpha}$ distinct $\stabu$-droplets spanned by sites in~$S$. Indeed, although $D$ might not be contained in $S$, any $\stabu$-droplet spanned by sites in $S$ is contained in $D_U(S)$ (the smallest $\stabu$-droplet that contains $S$), which has size $O(|S|^2)$. Also, an $\stabu$-droplet is determined uniquely by two opposite corners, for each of which there are at most $O(|D_U(S)|) = O(|S|^2) = O(p^{-6\alpha})$ choices. 

Now, recalling from~\eqref{eq:constants2} and~\eqref{eq:constants3} that 
$$\pi(S,u) \geq \xi^{-1} \gg \weak(\beta_1+\beta_2) \gg \delta(\beta_1+\beta_2)^{-1},$$
it follows that the probability that $S$ admits a $(u,\gamma)$-partition containing a saver droplet~$D$ satisfying \eqref{eq:Dconds} is at most 
$$p^{-13\alpha} \cdot p^{\delta(\beta_1+\beta_2) \pi(S,u)/3} \leq p^{\delta(\beta_1+\beta_2) \pi(S,u)/4} \leq p^{2\alpha m},$$
and this is at most~\eqref{eq:calc} (with $j = m$).


Let us therefore assume from now on that if $D$ is a saver droplet in a $(u,\gamma)$-partition of $S$, then the dimensions of $D$ satisfy
\begin{equation}\label{eq:Dconds2}
w(D) \le p^{-\beta_1(1-2\eta)-\eta} \qquad \text{and} \qquad h(D) \le p^{-\beta_2(1-2\eta)-\eta}.
\end{equation}
Let $S_1 \cup \dots \cup S_{m'+1}$ be a partition of $S$ into $u$-strips as in Definition~\ref{de:partition}, and note that we have $m' \ge m$, since $\| S_i - S_{i+1} \| \le 1$ for each $i \in [m']$. Note also that $m \ge 1$, since $\pi(S,u) \geq \xi^{-1} \gg \weak(\beta_1+\beta_2)$, and $\gamma \le \lambda$. 

Next, observe that, for each $1 \leq i \leq m$, the probability that $S_i \cap A$ contains a $u$-weak $\gamma$-cluster is at most
\begin{equation}\label{eq:partcluster}
1-\big(1-p^\gamma\big)^{g_u(S)},
\end{equation}
by Harris's inequality, since by definition there are at most $g_u(S)$ such sets in $S_i$. 

Now, as noted above, there are at most $p^{-13\alpha}$ distinct $\S_U$-droplets that are internally spanned by a subset of $S$, and by $\ih(\beta_1,\beta_2)$ and~\eqref{eq:Dconds2}, each such droplet $D$ is internally spanned with probability at most
\[
p^{\delta(\beta_1+\beta_2) \max\{ w(D), h(D) \}}.
\]
Thus, for each $a \in [m]$ and $0 \leq t \leq m - a$, the probability that there is a droplet $D$ with $\max\{ w(D), h(D) \} \ge a \cdot \weak(\beta_1+\beta_2) / 5$ such that $D$ is internally spanned by $(S_{t+1} \cup \cdots \cup S_{t+a}) \cap A$ is at most 
\begin{equation}\label{eq:partsaver}
p^{\delta(\beta_1+\beta_2) \weak(\beta_1+\beta_2) a/5 - 13\alpha} \leq p^{2\alpha a},
\end{equation}
since $\weak(\beta_1+\beta_2) \gg \delta(\beta_1+\beta_2)^{-1}$.

Finally, note that there are at most $\pi(S,u)^{2j}$ partitions of $m'$ containing at least $m' - j$ ones. By~\eqref{eq:partcluster} and~\eqref{eq:partsaver}, and taking a union bound over $j$, it follows that $S$ admits a $(u,\gamma)$-partition with probability at most
$$\pi(S,u) \cdot \max_{0\,\leq \,j\,\leq \,m} \Big(1-\big(1-p^\gamma\big)^{g_u(S)}\Big)^{m-j} \big( \pi(S,u)^2 \cdot p^{2\alpha}\big)^j,$$
as claimed.
\end{proof}

We shall now apply Lemmas~\ref{le:generalpartitions} and~\ref{le:partitionscalc} three times: once to prove Lemma~\ref{le:horizcrossings} for horizontal crossings, and twice to prove Lemma~\ref{le:vertcrossings} for vertical crossings, once each for drift and non-drift directions. We begin with horizontal crossings.


\begin{proof}[Proof of Lemma~\ref{le:horizcrossings}]
Suppose first that $\ih(\beta_1,\beta_2)$ holds, where $1 \leq \beta_1 \le \beta_2 \leq \alpha$, and let $S$ be a $u$-strip, where $u\in\{u^l,u^r\}$ and
\[
\xi^{-1} \leq \pi(S,u) \leq p^{-\beta_1(1-2\eta)-\eta} \qquad\text{and}\qquad h(S) \leq p^{-\alpha(1-2\eta)-\eta}.
\]
If $S$ is $u$-crossed by $A$, then, recalling from~\eqref{eq:choosing:stab:U} that $\alphabar(u) \geq \alpha$, it follows by Lemma~\ref{le:generalpartitions} that there exists a $(u,\alpha)$-partition for $S \cap A$.

Let $S_1 \cup \dots \cup S_{m'+1}$ be a partition of $S$ into $u$-strips as in Definition~\ref{de:partition}, and note that there are at most $O\big(h(S)\big)$ $u$-weak $\alpha$-clusters in each sub-strip $S_i$, where the implicit constant depends on $\weak=\weak(\beta_1+\beta_2)$. It therefore follows from Lemma~\ref{le:partitionscalc} that $S\cap A$ admits a $(u,\alpha)$-partition with probability at most
\begin{equation}\label{eq:horizpart}
\pi(S,u) \cdot \max_{0\,\leq \,j\,\leq \,m} \Big(1-\big(1-p^\alpha \big)^{O(h(S))}\Big)^{m-j} \big( \pi(S,u)^2 \cdot p^{2\alpha}\big)^j,
\end{equation}
where $m= \lfloor \pi(S,u) / 5\weak \alpha \rfloor$. Since $h(S) \le p^{-\alpha(1-2\eta)-\eta}$, we have\footnote{Here, and below, we use the inequality $1-xy \le (1-x)^y$, which is valid if $x \le 1 \le y$.}
$$1-(1-p^\alpha)^{O(h(S))} = O\big( p^{\alpha - \alpha(1-2\eta)-\eta} \big) \le p^{\eta/2}.$$
Also, since $\pi(S,u) \leq p^{-\beta_1(1-2\eta)-\eta}$ and $1 \le \beta_1 \le \alpha$, we have $\pi(S,u)^2 \cdot p^{2\alpha} \leq p^\eta$. Therefore, recalling that $m \gg 1$ (since $\pi(S,u) \geq \xi^{-1} \gg \weak(\beta_1+\beta_2)$), it follows that~\eqref{eq:horizpart} is at most
$$\pi(S,u) \cdot \max_{0\,\leq \,j\,\leq \,m} p^{\eta (m - j)/2} \cdot p^{\eta j} \le \pi(S,u) \cdot p^{\eta m / 2} \le p^{\delta' \pi(S,u)},$$
as required, where $\delta' = \delta'(\beta_1+\beta_2) \ll 1 / \weak(\beta_1+\beta_2)$.

Now suppose that $\ih(\alpha,\alpha+1)$ holds, and let $S$ be a $u$-strip, where $u\in\{u^l,u^r\}$ and 
$$\xi^{-1} \le \pi(S,u) \le p^{-\alpha(1-2\eta)-\eta} \qquad \text{and} \qquad h(S) \le \frac{\xi}{p^\alpha}\log\frac{1}{p}.$$
Then, exactly as above, it follows that~\eqref{eq:horizpart} is an upper bound on the probability that $S$ is $u$-crossed by $A$. Since $h(S)\leq\frac{\xi}{p^\alpha}\log\frac{1}{p}$, we have
\[
1-(1-p^\alpha)^{O(h(S))} \leq 1-\exp\Big(-O\big( p^\alpha \cdot h(S) \big)\Big) \leq 1-p^{O(\xi)} \leq e^{-p^{O(\xi)}},
\]
where the implicit constant depends on $\weak(2\alpha+1)$. Since $\pi(S,u) \le p^{-\alpha(1-2\eta)-\eta}$, we also have $\pi(S,u)^2 \cdot p^{2\alpha} \leq p^\eta$, as before. Since $\pi(S,u) = O(m)$, it follows that~\eqref{eq:horizpart} is at most
\[
\pi(S,u) \cdot \max_{0\,\leq \,j\,\leq \,m} \big( e^{-p^{O(\xi)}} \big)^{m-j} p^{\eta j} \le \pi(S,u) \cdot \exp\Big( - p^{O(\xi)} \cdot \pi(S,u) \Big),
\]
where the implicit constants depend on $\weak(2\alpha+1)$, as required.
\end{proof}

Finally, let us prove Lemma~\ref{le:vertcrossings}, which bounds the probability of a vertical crossing. When $\alphabar(u^*) \geq \alpha+1$ (the `non-drift' case), in which case $u = u^*$, the proof is straightforward; indeed, in this case the application of Lemma~\ref{le:generalpartitions} is the same as in the proof of Lemma~\ref{le:horizcrossings}. When $\alpha^-(u^*) = \infty$ (the `drift' case), and $u \in \stabu^+$ is such that $\sigma(u) = p^{1-\eta}$, on the other hand, this naive approach no longer works, and the proof in this case is conceptually a little more difficult, since it requires us to use the stretched geometry of the $u$-norm in order to control the unbounded sideways growth of small sets. This is the only point in the proof of Theorem~\ref{th:lowerunbal} where we specifically need the $u$-norm.

\begin{proof}[Proof of Lemma~\ref{le:vertcrossings}]
Let $u \in \stabu^+ \cap \Q_1$ and $\beta,\beta_1,\beta_2 \in \N$ be as in the statement of the lemma, and suppose that $\ih(\beta_1,\beta_2)$ holds. Let $S$ be a $u$-strip with 
$$w(S) \le p^{-(\alpha + 1)(1-2\eta)-\eta} \qquad \text{and} \qquad \xi^{-1} \le \pi(S,u) \le p^{-\beta(1-2\eta)-\eta},$$ 
and suppose that $S$ is $u$-crossed. We begin with the (easier) `non-drift' case, for which the proof is almost identical to that of Lemma~\ref{le:horizcrossings}.

\medskip
\noindent {\bf Case 1:} $u=u^*$ and $\alphabar(u^*) \ge \alpha + 1$. 
\medskip

Since $S$ is $u^*$-crossed by $A$, there exists a $(u^*,\alpha+1)$-partition for $S \cap A$, by Lemma~\ref{le:generalpartitions}. Let $S_1 \cup \dots \cup S_{m'+1}$ be a partition of $S$ into $u^*$-strips as in Definition~\ref{de:partition}, and note that $h(S) = \pi(S,u^*)$, and that there are at most $O\big(w(S)\big)$ $u^*$-weak $(\alpha+1)$-clusters in each sub-strip $S_i$, where the implicit constant depends on $\weak=\weak(\beta_1+\beta_2)$. It therefore follows from Lemma~\ref{le:partitionscalc} that $S \cap A$ admits a $(u^*,\alpha+1)$-partition with probability at most
\begin{equation}\label{eq:vertpart}
h(S) \cdot \max_{0\,\leq \,j\,\leq \,m} \Big(1-\big(1-p^{\alpha+1} \big)^{O(w(S))}\Big)^{m-j} \big( h(S)^2 \cdot p^{2\alpha}\big)^j,
\end{equation}
where $m= \lfloor h(S) / 5(\alpha+1)\weak(\beta_1+\beta_2) \rfloor$. Since $w(S) \le p^{-(\alpha + 1)(1-2\eta)-\eta}$, we have 
\[
1 - ( 1 - p^{\alpha+1} )^{O(w(S))} = O\big( w(S) \cdot p^{\alpha+1} \big) \le p^{\eta}.
\]
Also, since $h(S) = \pi(S,u) \le p^{-\alpha(1-2\eta)-\eta}$, we have $h(S)^2 \cdot p^{2\alpha} \le p^\eta$.  

\enlargethispage*{\baselineskip}

Thus, recalling that $h(S) \ge \xi^{-1} \gg \weak(\beta_1+\beta_2)$, and hence $h(S) = O(m)$, it follows that~\eqref{eq:vertpart} is at most
\[
h(S) \cdot \max_{0\,\leq \,j\,\leq \,m} p^{\eta (m - j)} \cdot p^{\eta j} = h(S) \cdot p^{\eta m} \le p^{\delta' h(S)} = p^{\delta' \pi(S,u)},
\]
as required, where $\delta' = \delta'(\beta_1+\beta_2) \ll 1 / \weak(\beta_1+\beta_2)$.

We now turn to the `drift' case.

\medskip
\noindent {\bf Case 2:} $\sigma(u) = p^{1-\eta}$ and $\alpha^-(u^*) = \infty$. 
\medskip

Since $S$ is $u$-crossed by $A$, and $\alphabar(u) = \infty$, by Lemma~\ref{le:generalpartitions} there exists a $(u,\lambda)$-partition for $S \cap A$. By Lemma~\ref{le:partitionscalc}, this occurs with probability at most
\begin{equation}\label{eq:vertpart2}
\pi(S,u) \cdot \max_{0\,\leq \,j\,\leq \,m} \Big( 1 - \big( 1 - p^\lambda \big)^{g_u(S)} \Big)^{m-j} \big( \pi(S,u)^2 \cdot p^{2\alpha} \big)^j,
\end{equation}
where $m := \lfloor \pi(S,u)/ 5\lambda \weak(\beta_1+\beta_2) \rfloor$, and $g_u(S)$ is the maximum number of $u$-weak $\lambda$-clusters in a sub-strip $S' \subset S$ of $u$-projection $4\lambda \weak(\beta_1+\beta_2)$. We claim that
$$g_u(S) = O\big( w(S) \cdot p^{-\lambda(1 - \eta)} \big),$$
where the implicit constant depends on $\weak(\beta_1+\beta_2)$. Indeed, 
there are $O\big(w(S)\big)$ choices for the first site in the $u$-weak $\lambda$-cluster, and at most $O(1/\sigma)$ choices for each of the remaining $\lambda - 1$ sites, as required. It follows that
\[
1 - \big( 1 - p^\lambda \big)^{g_u(S)} \leq O\big( w(S) \big) \cdot p^{-\lambda(1 - \eta)} \cdot p^\lambda = O\big( w(S) \big) \cdot p^{5\alpha} \leq p^{3\alpha},
\]
since $w(S) \le p^{-(\alpha + 1)(1-2\eta)-\eta}$, and recalling from~\eqref{eq:lambda} that $\lambda = 5\alpha / \eta$. 

Finally, note that $\pi(S,u)^2 \cdot p^{2\alpha} \leq p^\eta$, since $\pi(S,u) \le p^{-\alpha(1-2\eta)-\eta}$. Thus, recalling that $\pi(S,u) \ge \xi^{-1} \gg \weak(\beta_1+\beta_2)$, and hence $\pi(S,u) = O(m)$, it follows that~\eqref{eq:vertpart2} is at most
$$\pi(S,u) \cdot \max_{0\,\leq \,j\,\leq \,m} p^{3\alpha(m-j)} \cdot p^{\eta j} \leq \pi(S,u) \cdot p^{\eta m} \leq p^{\delta' \pi(S,u)},$$
where $\delta' = \delta'(\beta_1+\beta_2) \ll 1 / \weak(\beta_1+\beta_2)$. This completes the proof of the lemma.
\end{proof}

\subsection{Crossing droplets vertically}\label{sec:tall:droplets}

We will next use Lemma~\ref{le:vertcrossings} to bound the probability that a droplet is `crossed' vertically with help from above and below. 

\begin{definition}\label{def:vert:crossed}
Given an $\S_U$-droplet $D = \bigcap_{u \in \S_U} \H_u(a_u)$, let us write 
$$D^\updownarrow := \H_{u^l}(a_{u^l}) \cap \H_{u^r}(a_{u^r})$$
for the (infinite) set enclosed by the $u^l$- and $u^r$-sides of $D$, and define
$$\Dup := D^\updownarrow \setminus \H_{u^*}(a_{u^*})  \qquad \text{and} \qquad \Ddown := D^\updownarrow \setminus \H_{-u^*}(a_{-u^*})$$
to be the parts of $D^\updownarrow$ above and below $D$, respectively. 

We say that $D$ is \emph{vertically crossed} if there exists a strongly connected set in $[ \Dup \cup \Ddown \cup (D \cap A) ]$ that intersects both $\Dup$ and $\Ddown$. 
\end{definition}

We will use the following lemma several times in the proof below: first to bound the probability of the event $\Delta(D,D')$ (see~\eqref{def:deltaDD} and Definition~\ref{def:satisfied}), then to show that $\ih(\beta,\beta)\Rightarrow\ih(\beta,\beta+1)$ and $\ih(\beta+1,\beta)\wedge\ih(\beta,\beta+1)\Rightarrow\ih(\beta+1,\beta+1)$, and finally to bound the probability that a critical droplet is internally spanned. 

\begin{lemma}\label{le:indvert}
Let $\beta_1,\beta_2 \in \N$ with $\beta_1 + \beta_2 \le 2\alpha + 1$, set $\beta := \min\{\beta_1,\beta_2\}$, and suppose that $\ih(\beta_1,\beta_2)$ holds. Let $D$ be an $\stabu$-droplet such that
$$w(D) \le p^{- (\beta + 1)(1 - 2\eta) - \eta} \quad \text{and} \quad h(D) \ge \xi^{-3} \max\big\{ p^{1-\eta} w(D), \, 1 \big\}. $$
Then 
$$\Pr_p\big( D \text{ is vertically crossed}\, \big) \le p^{\delta' h(D) / 8},$$
where $\delta' = \delta'(\beta_1+\beta_2)$. 
\end{lemma}

Before we can prove Lemma~\ref{le:indvert}, we need one more technical lemma, which says that if a droplet $D$ is vertically crossed (that is, with help from \emph{both} the $u^*$-side and the $(-u^*)$-side) then, in a certain sense (which is made explicit in the lemma), the droplet is at least `half-crossed' with help from just one side. This will allow us to deduce that a  vertically crossed droplet contains a $u$-crossed $u$-strip, for an appropriate $u$. In Figure~\ref{fi:halfcrossed}, below, we illustrate an application of this lemma.

\begin{lemma}\label{le:halfcrossed}
Let $D$ be an $\stabu$-droplet. If $D$ is vertically crossed, then for some $* \in\{ \uparrow, \downarrow\}$ there exists a set $L \subset D$ such that 
$$h(L) \geq h(D)/2 - 2\kappa$$
and $L \cup D^*$ is a strongly connected component of $\big[D^* \cup ( D \cap A ) \big]$. 
\end{lemma}

\begin{proof}
For each $* \in\{ \uparrow, \downarrow\}$, let us write $Z^*$ for the strongly connected component of $\big[ D^* \cup ( D \cap A ) \big]$ containing $D^*$, and note that $Z^* \subset D^* \cup D$, since $D^* \cup D$ is an $\S_U$-droplet. If the set $Z^\uparrow \cup Z^\downarrow$ is strongly connected, then set $L^* := Z^* \cap D$ for each $* \in\{ \uparrow, \downarrow\}$, and observe that 
\[
h(L^\uparrow) + h(L^\downarrow) \geq h(D) - 3\kappa,
\]
as required.

So suppose that $Z^\uparrow \cup Z^\downarrow$ is not strongly connected, and let $\mathcal{Y}$ be the collection of strongly connected components of $\big[ (D \cap A) \setminus (Z^\uparrow \cup Z^\downarrow) \big]$. Then $Z^* \cup Y$ is not strongly connected for any $* \in\{ \uparrow, \downarrow\}$ and $Y \in \mathcal{Y}$, and thus
\[
\Y \cup \big\{ Z^\uparrow, Z^\downarrow \big\}
\]
is precisely the collection of strongly connected components of $[ \Dup \cup \Ddown \cup (D \cap A) ]$. Since $[ \Dup \cup \Ddown \cup (D \cap A) ]$ contains a strongly connected component intersecting both $\Dup$ and $\Ddown$, this is a contradiction, which completes the proof of the lemma.
\end{proof}

\begin{figure}[ht]
  \centering
  \begin{tikzpicture}[>=latex,scale=1]
    \begin{scope}
      \clip (0,0) -- (5,0) -- (5.3,1.8) -- (-0.6,1.8) -- cycle;
      \fill [pattern=north west lines] (-0.6,1.8) -- (5.3,1.8) -- (5.3,1.6) -- (-0.6,1.6) -- cycle;
    \end{scope}
    \begin{scope}
      \clip (6,6) -- (-2,6) -- (-1.4,4.2) -- (5.7,4.2) -- cycle;
      \fill [pattern=north west lines] (-1.6,4.2) -- (5.9,4.2) -- (5.9,4.4) -- (-1.6,4.4) -- cycle;
    \end{scope}
    \draw (0,0) -- (5,0) -- (6,6) -- (-2,6) -- cycle;
    \draw (-0.6,1.8) -- (5.3,1.8) (-1.4,4.2) -- (5.7,4.2);
    \draw (0,1.8) -- (0.3,2.4) -- (0.5,2.1) -- (0.6,2.2) -- (0.9,1.8);
    \draw (1.5,1.8) -- (1.8,2) -- (2.1,1.8);
    \draw (2.5,1.8) -- (3,3.5) -- (3.9,2.7) -- (4.3,1.8);
    \draw [dashed,name path=H] (-2,1.66) -- (7,2.56);
    \path [name path=S1] (0,1.8) -- (0.3,2.4);
    \path [name path=S2] (3.9,2.7) -- (4.3,1.8);
    \path [name intersections={of=H and S1,by=P1}];
    \path [name intersections={of=H and S2,by=P2}];
    \path [name path=T] (0,3.2) -- (3,3.5) -- (5,3.7);
    \path [name path=V1] (P1) -- ++(0,3);
    \path [name path=V2] (P2) -- ++(0,3);
    \path [name intersections={of=T and V1,by=Q1}];
    \path [name intersections={of=T and V2,by=Q2}];
    \draw [dashed] (P1) -- (Q1) -- (Q2) -- (P2);
    \node at (2.5,0.8) {$D \cap D_i^\downarrow$};
    \node at (2.5,5.2) {$D \cap D_i^\uparrow$};
    \node at (-0.4,3.7) {$D_i$};
    \draw (P2) ++(0,1) -- ++(0.5,0) node [right] {$T$};
    \draw (2.9,2.4) -- ++(-1,0.3) node [left] {$L$};
    \draw (0.3,2.1) -- (1.3,2.7) (1.8,1.9) -- (1.7,2.4);
    \draw (-0.6,1.8) node [circle,fill,inner sep=0pt,minimum size=0.1cm] {} -- (-1,1.3) node [below left] {$x$};
    \node at (7,2.2) {$\H_u(x)$};
    \draw [->] (2.5,0) -- ++(0,-0.5);
    \node at (2.5,-0.8) {$-u^*$};
    \draw [->] (2.5,6) -- ++(0,0.5);
    \node at (2.5,6.8) {$u^*$};
    \draw [->] (5.8,4.8) -- ($(5.8,4.8)!0.5cm!(7,4.6)$);
    \node at ($(5.8,4.8)!0.8cm!(7,4.6)$) {$u^r$};
    \draw [->] (-1.6,4.8) -- ($(-1.6,4.8)!0.5cm!(-2.2,4.6)$);
    \node at ($(-1.6,4.8)!0.8cm!(-2.2,4.6)$) {$u^l$};
  \end{tikzpicture}
  \caption{The figure depicts the application of Lemma~\ref{le:halfcrossed} in the proof of Lemma~\ref{le:indvert} assuming $\alpha^-(u^*)=\infty$. $T$ is the minimal $u$-strip such that $L \subset T \cup\H_u(x)$, where $\sigma(u) = p^{1-\eta}$, and is $u$-crossed by $A$.}
  \label{fi:halfcrossed}
\end{figure}

We can now deduce Lemma~\ref{le:indvert} from Lemmas~\ref{le:vertcrossings} and~\ref{le:halfcrossed}. 

\begin{proof}[Proof of Lemma~\ref{le:indvert}]
Suppose that $D$ is vertically crossed, and choose a partition $D = D_1 \cup \cdots \cup D_m$ such that each $D_i$ is a finite $\S_U$-droplet with $D_i^\updownarrow = D^\updownarrow$, and 
$$\xi^{-2} \max\big\{ p^{1-\eta} w(D), \, 1 \big\} \le h(D_i)  \le 2 \cdot \xi^{-2} \max\big\{ p^{1-\eta} w(D), \, 1 \big\}$$
for each $i \in [m]$, noting that such a partition exists by the lower bound on $h(D)$. We claim that $D_i$ is vertically crossed for each $i \in [m]$. Indeed, noting that
$$\Dup \cup \Ddown \cup (D \cap A) \subset D_i^\uparrow \cup D_i^\downarrow \cup (D_i \cap A),$$
and recalling Definition~\ref{def:vert:crossed}, it follows that there exists a strongly connected set in $[ D_i^\uparrow \cup D_i^\downarrow \cup (D_i \cap A) ]$ that intersects both $D_i^\uparrow$ and $D_i^\downarrow$, as required.

By Lemma~\ref{le:halfcrossed}, it follows that there exists $L \subset D_i$ with $h(L) \geq h(D_i)/2 - 2\kappa$ such that $L \cup D_i^*$ is a strongly connected component of $\big[D_i^* \cup( D_i \cap A ) \big]$ for some $* \in\{ \uparrow, \downarrow\}$. Let $E_i$ denote the event that such a set $L$ exists in $D_i$.
 
\begin{claim}\label{cl:vert:crossing}
$\Pr_p(E_i) \leq p^{\delta' h(D_i) / 7}$ for each $1 \le i \le m$, where $\delta' = \delta'(\beta_1+\beta_2)$.
\end{claim}

\begin{proof}[Proof of Claim~\ref{cl:vert:crossing}]\let\qednow\qedsymbol\renewcommand{\qedsymbol}{}
Suppose (without loss of generality) that there exists $L \subset D_i$ with $h(L) \geq h(D_i)/2 - 2\kappa$ such that $L \cup D_i^\downarrow$ is a strongly connected component of $\big[ D_i^\downarrow \cup( D_i \cap A ) \big]$. Note that each site of $D_i \cap A$ used to infect a member of $L$ is itself a member of $L$ (since $L \cup D_i^\downarrow$ is a strongly connected component, and $\kappa > \nu$), so $L \cup D_i^\downarrow$ is also a strongly connected component of $[ D_i^\downarrow \cup( L \cap A ) ]$. We will divide into two cases according to whether or not $u^*$ is a drift direction.

\medskip
\noindent {\bf Case 1:} $\alphabar(u^*) \ge \alpha + 1$. 
\medskip

We first claim that the minimal $u^*$-strip $S$ containing $L$ is $u^*$-crossed. To see this, observe that $L \subset [ \H_{u^*}(x) \cup (S \cap A) ]$ for any $x \in \partial(S,-u^*)$, since $L \cup D_i^\downarrow$ is a strongly connected component of $[ D_i^\downarrow \cup ( L \cap A) ]$, and $D_i^\downarrow \subset \H_{u^*}(x)$. Note also that $w(S) = w(L) \le w(D) \le p^{-(\beta + 1)(1-2\eta)-\eta}$, and hence
$$\pi(S,u^*) \le h(D_i) \le 2 \cdot \xi^{-2} \max\big\{ p^{1-\eta} w(D), 1 \big\} \le p^{-\beta(1-2\eta)-\eta}.$$
Moreover, we have 
$$\pi(S,u^*) = h(L) \geq h(D_i)/3 \geq \xi^{-1},$$
since $h(D_i) \ge \xi^{-2} \gg \kappa$. Therefore, $S$ satisfies the conditions of Lemma~\ref{le:vertcrossings}, and is thus $u^*$-crossed with probability at most $p^{\delta' \pi(S,u^*)} \le p^{\delta' h(D_i)/3}$. 

Now, recall 
that $S$ is a $\T$-droplet, where $\T = \{u^*,-u^*,w^\perp,-w^\perp \}$, and so is defined by two opposite corners. There are at most $w(D) \cdot h(D_i)$ choices for each, and hence at most $p^{-13\alpha}$ choices for $S$. Thus, taking a union bound over $u^*$-strips $S$, and recalling that $h(D_i) \ge \xi^{-2} \gg \delta'(\beta_1+\beta_2)^{-1}$, it follows that 
$$\Pr_p(E_i) \le p^{-13\alpha} \cdot p^{\delta' h(D_i) / 3} \le p^{\delta' h(D_i) / 4},$$
as claimed.

\medskip
\noindent {\bf Case 2:} $\alpha^-(u^*)=\infty$. 
\medskip

Let $x$ be the element of $\R^2$ at the intersection of the $u^l$ and $(-u^*)$-sides of $D_i$, let $u \in \stabu^+ \cap \Q_1$ be such that $\sigma(u) = p^{1-\eta}$, and let $T$ be the minimal $u$-strip such that
$$L \subset \H_u(x)\cup T$$
(see Figure~\ref{fi:halfcrossed}). Observe that $T$ is $u$-crossed; indeed, since $D_i^\downarrow \subset \H_u(x)$ and $L \subset \big[ D_i^\downarrow \cup ( L \cap A) \big]$, it follows that $L \subset [ \H_u(x) \cup (T \cap A) ]$.

We next claim that $w(T)$ and $h(T)$ satisfy the conditions of Lemma~\ref{le:vertcrossings}. Indeed, we have $w(T) \le w(L) \le w(D) \le p^{-(\beta + 1)(1 - 2\eta) - \eta}$, and 
$$\pi(T,u) \le h(D_i) + \sigma(u) \cdot w(D) \leq p^{-\beta(1 - 2\eta) - \eta},$$ 
since $h(D_i) \le 2 \cdot \xi^{-2} \max\big\{ p^{1-\eta} w(D), 1 \big\}$ and $\sigma(u) = p^{1-\eta}$ (cf.~Case~1). Moreover, 
$$\pi(T,u) \ge \big( h(L) - \sigma(u) \cdot w(D) \big) \cos \sigma(u) \ge \frac{h(L)}{2} \ge \xi^{-1},$$
since $h(L) \ge h(D_i)/3 \ge 2 \cdot \xi^{-1} \max\big\{ p^{1-\eta} w(D), \, 1 \big\}$. Thus, by Lemma~\ref{le:vertcrossings}, $T$ is $u$-crossed with probability at most $p^{\delta' \pi(T,u)} \le p^{\delta' h(D_i)/6}$. 

Now, as in Case~1, there are at most $p^{-13\alpha}$ choices for $T$. Thus, taking a union bound over $u$-strips $T$, and recalling again that $h(D_i) \ge \xi^{-2} \gg \delta'(\beta_1+\beta_2)^{-1}$, it follows that 
$$\Pr_p(E_i) \le p^{-13\alpha} \cdot p^{\delta' h(D_i) / 6} \le p^{\delta' h(D_i) / 7},$$
as claimed. \qednow
\end{proof}

Finally, note that the events $E_1,\dots,E_m$ are independent, since $E_i$ depends only on the set $D_i \cap A$. Therefore, by Claim~\ref{cl:vert:crossing}, 
$$\Pr_p\big( D \text{ is vertically crossed}\, \big) \le \exp\bigg( \frac{\delta'}{7}  \log \frac{1}{p} \sum_{i = 1}^m h(D_i) \bigg) \leq p^{\delta' h(D) / 8},$$
and this completes the proof of the lemma.
\end{proof}

\subsection{From crossing events to growing droplets}\label{sec:growing:droplets}

Recall from~\eqref{def:deltaDD} the event
$$\Delta(D,D') = \big\{ D' \in \< D \cup (D' \cap A) \> \big\},$$
which appears in the definition of a good and satisfied hierarchy. In this subsection we will use Lemmas~\ref{le:horizcrossings},~\ref{le:vertcrossings} and~\ref{le:indvert} to deduce upper bounds on the probability of the event $\Delta(D,D')$ which hold when $D'$ and $D$ are the droplets corresponding to adjacent vertices of outdegree~$1$ in a good hierarchy (see Definition~\ref{de:hier2}). We will use these bounds, together with Lemma~\ref{le:boundoverH}, first in Section~\ref{se:induction} to prove the induction steps, and then again in Section~\ref{se:criticaldroplets} to bound the probability that a critical droplet is internally spanned.

\begin{lemma}\label{le:deltagrowth}
Let $D \subset D'$ be nested $\S_U$-droplets.
\begin{itemize}
\item[(i)] Let $1 \le \beta_1 \le \beta_2 \le \alpha$, and suppose that $\ih(\beta_1,\beta_2)$ holds. If
$$p^{-\beta_1(1-2\eta)-\eta} \le w(D') \le p^{-(\beta_1+1)(1-2\eta)-\eta}, \qquad h(D') \le p^{-\alpha(1-2\eta)-\eta},$$
and
$$\xi \cdot p^{-\beta_1(1-2\eta)-\eta} \le w(D') - w(D) \le 2\xi \cdot p^{-\beta_1(1-2\eta)-\eta},$$
then
$$\Pr_p\big(\Delta(D,D')\big) \le p^{\Omega(\delta')(w(D')-w(D))},$$
where $\delta'=\delta'(\beta_1+\beta_2)$, and the implicit constant depends on $\S_U$.\smallskip
\item[(ii)] Suppose that $\ih(\alpha,\alpha+1)$ holds. If 
$$p^{-\alpha(1-2\eta)-\eta} \le w(D') \leq p^{-(\alpha+1)(1-2\eta)-\eta}, \qquad h(D') \le \frac{\xi}{p^{\alpha}}\log\frac{1}{p},$$
and
$$\xi \cdot p^{-\alpha(1-2\eta)-\eta} \le w(D') - w(D) \le 2\xi \cdot p^{-\alpha(1-2\eta)-\eta},$$
then
$$\Pr_p\big(\Delta(D,D')\big) \le \exp\Big(-p^{O(\xi)}\big(w(D')-w(D)\big)\Big),$$
where the implicit constant depends on $\weak(2\alpha+1)$.
\end{itemize}
\end{lemma}

\begin{figure}[ht]
 \begin{minipage}{.55\textwidth}
  \centering
  \centering
  \begin{tikzpicture}[>=latex,scale=1]
     \draw (0,0) -- (5,0) -- (5.7,4.2) -- (-1.4,4.2) -- cycle;
     \draw (0.8,0.6) -- (4.2,0.6) -- (4.7,3.6) -- (-0.2,3.6) -- cycle;
     \draw (4.1,0) -- (4.2,0.6) (4.7,3.6) -- (4.8,4.2);
     \fill [pattern=north west lines] (3.9,0) -- (4.1,0) -- (4.8,4.2) -- (4.6,4.2) -- cycle;
     \node at (2.5,2.1) {$D$};
     \draw (0.4,0.3) -- (0.4,-0.3) node [below] {$D'$};
     \draw (4.6,0.3) -- (4.6,-0.3) node [below] {$S^r$};
  \end{tikzpicture}
 \end{minipage}
 \begin{minipage}{.35\textwidth}
  \centering
  \begin{tikzpicture}[>=latex,scale=1]
     \draw (0,0) -- (0.4,0) -- (1.4,6) -- (-2,6) -- cycle;
     \draw (0.1,2.7) -- (0.5,5.1) -- (-0.7,5.1) -- cycle;
     \draw (-0.35,0) -- (0.1,2.7) (0.5,5.1) -- (0.65,6) (-0.35,0) -- (0,0);
     \fill [pattern=north west lines] (-0.55,0) -- (-0.35,0) -- (0.65,6) -- (0.45,6) -- cycle;
     \draw (-0.1,4.3) -- (-1.8,4.3) node [left] {$D$};
     \draw (-0.7,3) -- (-1.8,3) node [left] {$D'$};
     \draw (0.7,3.6) -- (1.4,3.6) node [right] {$S^r$};
     \draw [dotted] (1,6) -- (1.1,5) (1.1,5) -- (0.5,4) (0.5,4) -- (0.6,3) (0.67,3) -- (0.87,3)(0.6,3) -- (0.25,2) (0.25,1.98) -- (0.3,1) (0.3,0.98) -- (0.045,0.528) (0.04,0.31) -- (-0.11,0.3) (0.06,0.46) -- (0.2,-0.02);
  \end{tikzpicture}
 \end{minipage}
 \caption{Two possible configurations of the droplets in Lemma~\ref{le:deltagrowth}. The hatching indicates the sites that would be assumed to be present for the purposes of the event that $S^r$ is $u^r$-crossed. In the first example, the event $\Delta(D,D')$ implies that $S^r$ is $u^r$-crossed, but in the second it does not, e.g, if $A \cap D'$ is the set of dots in $S^r$.}
 \label{fi:Sr}
\end{figure}

The idea behind the proof of the lemma is that the event $\Delta(D,D')$ `should' imply that the strips to the right and left of $D$ are (respectively) $u^l$-crossed and $u^r$-crossed (see Definition~\ref{def:ucrossed}). Unfortunately, however, this is not always necessarily the case (see Figure~\ref{fi:Sr}), and the proof is therefore somewhat delicate. We are  grateful to the anonymous referee for pointing out errors in the first two versions of this lemma. 

Let $S^r\subset D'$ be the unique maximal $u^r$-strip whose $u^r$-side is equal to that of $D'$ and which does not intersect $D$ ($S^r$ may or may not be a subset of $D'$, depending on the shapes and relative positions of $D$ and $D'$; see Figure~\ref{fi:Sr}). Define $S^l$ similarly on the $u^l$-side of $D'$. The first step in the proof of Lemma~\ref{le:deltagrowth} is the following simple observation.

\begin{obs}\label{obs:growing:width}
If $D \subset D'$ are nested $\S_U$-droplets, then
$$w(D') - w(D) = O\big( \pi(S^r,u^r) + \pi(S^l,u^l) + h(D') - h(D) \big),$$
where the implicit constant depends only on $\S_U$. 
\end{obs}

\begin{proof}
Let $D \subset D'' \subset D'$ be the $\S_U$-droplet with $h(D'') = h(D)$ and $(D'')^\updownarrow = (D')^\updownarrow$. Now, simply observe that $w(D'') - w(D)$ is at most a constant times $\pi(S^r,u^r) + \pi(S^l,u^l)$, and $w(D') - w(D'')$ is at most a constant times $h(D') - h(D)$, where the constants depend on the angles between $u^l$, $u^r$ and $\pm u^*$. 
\end{proof}

If $h(D') - h(D) = \Omega\big( w(D') - w(D) \big)$, then we will apply Lemma~\ref{le:indvert} to bound the probability that each of the droplets forming $D' \setminus D''$ is vertically crossed. On the other hand, if $\pi(S^r,u^r) = \Omega\big( w(D') - w(D) \big)$ (say), then we will use the following deterministic lemma. If $D$ is an $\S_U$-droplet, then we will write $\sh(D)$ for the length of its shorter horizontal side. 

\begin{lemma}\label{lem:det:crossed:or:spanned}
If the event $\Delta(D,D')$ occurs, then either $S^r$ is $u^r$-crossed or $D'$ is internally spanned. Moreover, if $S^r$ is not $u^r$-crossed then $\sh(D') \le w(D') - w(D)$. 
\end{lemma}

\begin{proof}
If the event $\Delta(D,D')$ occurs, i.e. $D' \in \< D \cup (D' \cap A) \>$, then by Lemma~\ref{le:spanoutput} there exists a strongly connected component $L$ of $[D \cup (D' \cap A)]$ such that $D'$ is the smallest $\S_U$-droplet containing $L$. It follows that $L$ contains a point on each of the sides of $D'$, and therefore intersects $\partial(S^r,u^r) = \partial(D',u^r)$. 

Let $x \in \partial(S^r,-u^r)$, and suppose first that $L \cap \H_{u^r}(x)$ is non-empty. We claim that in this case $S^r$ is $u^r$-crossed. Indeed, $L$ is strongly connected and intersects both $\H_{u^r}(x)$ and $\partial(S^r,u^r)$, so (recalling Definition~\ref{def:ucrossed}) we only need to observe that 
$$L \subset [D \cup (D' \cap A)] \subset [ \H_{u^r}(x) \cup (S^r \cap A)],$$
where the second inclusion holds since $D \subset \H_{u^r}(x)$ and $D' \subset \H_{u^r}(x) \cup S^r$.

So suppose from now on that $L \cap \H_{u^r}(x)$ is empty, and therefore $L \subset S^r$. We claim that in this case $D'$ is internally spanned, and $\sh(D') \le w(D') - w(D)$. To bound $\sh(D')$, recall that $L$ contains an element $y \in \partial(D',u^l)$, and observe that this element must lie in $S^r$. Now let $z \in  \partial(S^r,u^r) = \partial(D',u^r)$ be chosen such that $\< y - z, u^* \> = 0$, and observe (see Figure~\ref{fi:Sr}) that 
$$\sh(D') \le \| y - z \| \le w(D') - w(D),$$
as claimed, where the final inequality holds by comparing each quantity to the length of the intersection of a (continuous) horizontal line with $S^r$.   

Finally, to show that $D'$ is internally spanned, recall that $L$ is strongly connected, and that $D'$ is the smallest $\S_U$-droplet containing $L$. By Definition~\ref{de:ispan}, it will therefore suffice to prove that $L \subset [D' \cap A]$. To see this, observe once again that any element that was used to infect an element of $L$ is included in it (since $L$ is a strongly connected component and $\kappa > \nu$), and $L$ contains no element of $\H_{u^r}(x) \supset D$. It follows that $D'$ is internally spanned, as required.
\end{proof}

We can now deduce Lemma~\ref{le:deltagrowth} using Lemmas~\ref{le:horizcrossings},~\ref{le:vertcrossings} and~\ref{le:indvert}.

\begin{proof}[Proof of Lemma~\ref{le:deltagrowth}]
Let $1 \le \beta_1 \le \beta_2 \le \alpha + 1$, and suppose that $\ih(\beta_1,\beta_2)$ holds, and that
$$w(D') \le p^{-(\beta_1+1)(1-2\eta)-\eta} \qquad \text{and} \qquad  w(D') - w(D) \ge \xi \cdot p^{-\beta_1(1-2\eta)-\eta},$$
noting that these conditions are assumed in both parts of Lemma~\ref{le:deltagrowth}. We first consider the case where $h(D') - h(D)$ is large. Let us introduce, just for this proof, a constant $\zeta > 0$ that is chosen to be sufficiently small depending on $\S_U$,\footnote{As in Observation~\ref{obs:growing:width}, this constant only depends on the angles between $u^l$, $u^r$ and $\pm u^*$.} but that does not depend on any of the other constants (so, in particular, $1 \gg \zeta \gg \delta(2) > 0$). 

\begin{claim}\label{claim:Delta:vertical}
If $h(D') - h(D) \ge 2\zeta \cdot \big( w(D') - w(D) \big)$, then 
$$\Pr_p\big( \Delta(D,D') \big) \le p^{\Omega(\delta') ( w(D') - w(D) )},$$
where $\delta' = \delta'(\beta_1+\beta_2)$, and the implicit constant depends on $\S_U$. 
\end{claim}

Note that the bound given by Claim~\ref{claim:Delta:vertical} is exactly the bound required in part~$(i)$ of the lemma, and is stronger than that required in part~$(ii)$. In the proof of the claim we shall only use the assumptions stated above, which hold in both cases.

\begin{proof}[Proof of Claim~\ref{claim:Delta:vertical}]\let\qednow\qedsymbol\renewcommand{\qedsymbol}{}
Define $D \subset D'' \subset D'$ (as in the proof of Observation~\ref{obs:growing:width}) to be the $\S_U$-droplet with $h(D'') = h(D)$ and $(D'')^\updownarrow = (D')^\updownarrow$, and consider the two droplets $D' \cap (D'')^\uparrow$ and $D' \cap (D'')^\downarrow$. Setting $\hat{D}$ to be the element of this pair with larger height (choosing arbitrarily if they have equal height), it follows that
$$h(\hat{D}) \ge \zeta \cdot \big( w(D') - w(D) \big) \ge \xi^2 \cdot p^{-\beta_1(1-2\eta)-\eta} \ge \xi^{-3} \max\big\{ p^{1-\eta} w(\hat{D}), \, 1 \big\},$$
since $w(\hat{D}) \le w(D') \le p^{-(\beta_1+1)(1-2\eta)-\eta}$. 

Now, observe that the event $\Delta(D,D')$ implies that the droplet $\hat{D}$ is vertically crossed. Indeed, the event $\Delta(D,D')$ implies that there exists a strongly connected component $L$ of $[D \cup (D' \cap A)]$ that intersects $\partial(D',u^*)$ and $\partial(D',-u^*)$, and $D \cup (D' \cap A) \subset \hat{D}^\uparrow \cup \hat{D}^\downarrow \cup (\hat{D} \cap A)$. Since $\beta_1 \le \beta_2$, it follows by Lemma~\ref{le:indvert} that 
$$\Pr_p\big( \Delta(D,D') \big)  \le \Pr_p\big( \hat{D} \text{ is vertically crossed}\, \big) \le p^{\delta' h(\hat{D}) / 8} = p^{\Omega(\delta') ( h(D') - h(D) )},$$
where $\delta' = \delta'(\beta_1+\beta_2)$, as claimed. \qednow
\end{proof}

We may therefore assume that $h(D') - h(D) \le 2\zeta \big( w(D') - w(D) \big)$, and hence (without loss of generality) that 
\begin{equation}\label{eq:lower:bound:on:piSr}
\pi(S^r,u^r) \ge \frac{\pi(S^r,u^r) + \pi(S^l,u^l)}{2} \ge \zeta \cdot \big( w(D') - w(D) \big),
\end{equation}
by Observation~\ref{obs:growing:width}. We will first bound the probability that $S^r$ is $u^r$-crossed, using Lemma~\ref{le:horizcrossings}, and then deal with the remaining case using Lemmas~\ref{le:indvert} and~\ref{lem:det:crossed:or:spanned}.

To bound the probability that $S^r$ is $u^r$-crossed, note first that
$$\xi^{-1} \le \zeta \cdot \big( w(D') - w(D) \big) \le \pi(S^r,u^r) \le w(D') - w(D) \le 2\xi \cdot p^{-\beta_1(1-2\eta)-\eta},$$
by~\eqref{eq:lower:bound:on:piSr} and the assumed upper bound on $w(D') - w(D)$, so the bounds on $\pi(S^r,u^r)$ required for Lemma~\ref{le:horizcrossings} hold in either case. Now, if $1 \le \beta_1 \le \beta_2 \le \alpha$, then the bound on $h(D')$ in part~$(i)$ of the lemma being proved implies that 
$$h(S^r) = h(D') \le p^{-\alpha(1 - 2\eta) - \eta}.$$ 
Hence, by Lemma~\ref{le:horizcrossings}, and again using~\eqref{eq:lower:bound:on:piSr}, we have 
$$\Pr_p\big( \text{$S^r$ is $u^r$-crossed} \big) \le p^{\delta' \pi(S^r,u^r)} = p^{\Omega(\delta') ( w(D') - w(D) )},$$
where $\delta' = \delta'(\beta_1+\beta_2)$, as required. On the other hand, if $\beta_1 = \alpha$ and $\beta_2 = \alpha + 1$, then the bound on $h(D')$ in part~$(ii)$ of the lemma being proved implies that 
$$h(S^r) = h(D') \le \frac{\xi}{p^{\alpha}}\log\frac{1}{p}.$$ 
Consequently, by Lemma~\ref{le:horizcrossings}, $S^r$ is $u^r$-crossed with probability at most 
$$\pi(S^r,u^r) \, \exp\Big(-p^{O(\xi)} \pi(S^r,u^r) \Big) = \exp\Big(-p^{O(\xi)}\big(w(D')-w(D)\big) \Big),$$
as required, where in the final step we used the inequality~\eqref{eq:lower:bound:on:piSr}, and the lower bound $w(D')-w(D) \ge \xi \cdot p^{-\alpha(1-2\eta)-\eta}$ assumed in the statement of the lemma. 

Finally, suppose that $S^r$ is not $u^r$-crossed. By Lemma~\ref{lem:det:crossed:or:spanned}, it follows that $D'$ is internally spanned, and 
$$\sh(D') \le w(D') - w(D) \le 2 \xi \cdot p^{-\beta_1(1-2\eta)-\eta}.$$ 
Since $w(D') \ge p^{-\beta_1(1-2\eta)-\eta}$, and since $\zeta$ was chosen sufficiently small depending on $\S_U$, it follows that 
$$h(D') \ge 2\zeta \cdot \big( w(D') - \sh(D') \big) \ge \zeta \cdot p^{-\beta_1(1-2\eta)-\eta},$$
and therefore, recalling that $w(D') \le p^{- (\beta_1 + 1)(1 - 2\eta) - \eta}$, 
$$h(D') \ge w(D') - w(D) \quad \text{and} \quad h(D') \ge \xi^{-3} \max\big\{ p^{1-\eta} w(D'), \, 1 \big\}.$$
Note also that $D'$ is vertically crossed, since it is internally spanned. Hence, by Lemma~\ref{le:indvert}, we have
$$\Pr_p\big( \ispan(D')\big) \le \Pr_p\big( D' \text{ is vertically crossed}\, \big) \le p^{\delta' h(D') / 8} \le p^{\Omega(\delta') ( w(D') - w(D) )},$$
where $\delta' = \delta'(\beta_1+\beta_2)$, as required. 
\end{proof}



\subsection{Wide internally spanned droplets}\label{se:induction}

Having bounded the probability of the event $\Delta(D,D')$ in the previous subsection, we are now in a position to bound the probability that a droplet is internally spanned using the method of hierarchies. The following lemma will be used to deduce the implications $\ih(\beta,\beta) \Rightarrow \ih(\beta+1,\beta)$ and $\ih(\beta+1,\beta) \wedge \ih(\beta,\beta+1) \Rightarrow \ih(\beta+1,\beta+1)$. 

\begin{lemma}\label{le:indhoriz}
Let $1 \leq \beta_1 \leq \alpha$ and $1 \leq \beta_2 \leq \alpha$ with $\beta_1 \le \beta_2 \le \beta_1 + 1$, and suppose that $\ih(\beta_1,\beta_2)$ holds. Let $D$ be an $\stabu$-droplet such that
\[
p^{-\beta_1(1 - 2\eta) - \eta} \leq w(D) \leq p^{-(\beta_1+1)(1 - 2\eta) - \eta} \qquad \text{and} \qquad h(D)\leq p^{-\beta_2(1 - 2\eta) - \eta}.
\]
Then 
\[
\Pr_p\big(\ispan(D)\big) \leq p^{\Omega(\delta') w(D)},
\]
where $\delta' = \delta'(\beta_1+\beta_2)$, and the implicit constant depends only on $\U$.
\end{lemma}

For convenience, we shall occasionally use the notation $\exp_p(x):=p^x$.

\begin{proof}[Proof of Lemma~\ref{le:indhoriz}]
We will use the hierarchies framework from Section~\ref{se:hierarchies} with $\beta=\beta_1$. To begin, recall the bound from Lemma~\ref{le:boundoverH}:
\begin{equation}\label{eq:indbound}
\Pr_p\big(\ispan(D)\big)\le  \sum_{\hier \in \hier_D} \bigg( \prod_{u \in L(\hier)} \Pr_p\big(\ispan(D_u)\big) \bigg)\bigg( \prod_{u \to v} \Pr_p\big(\Delta(D_v,D_u)\big) \bigg).
\end{equation}
To bound the right-hand side of~\eqref{eq:indbound}, we will use the induction hypothesis to bound the probability that a seed is internally spanned, Lemma~\ref{le:deltagrowth} to bound the probability of the event $\Delta(D_v,D_u)$, and Lemma~\ref{le:counthier} to bound the number of good hierarchies for $D$. The details are standard, but the following fact will play a crucial role in the argument: if there are `many' seeds, then our bound on the first term in the sum above will be sufficient (since all seeds are large, see Lemma~\ref{obs:all:seeds:big}), while if there are `few' seeds then we will be able to apply Lemma~\ref{lem:sumofwidths}.


We now provide the details of the calculation. First, recall from Definition~\ref{de:hier2} that we have $w(D_u)\leq p^{-\beta_1(1 - 2\eta) - \eta}$ for each $u\in L(\hier)$. Also, since $D_u\subset D$, we have $h(D)\leq p^{-\beta_2(1 - 2\eta) - \eta}$. Hence, by $\ih(\beta_1,\beta_2)$,
\begin{equation}\label{eq:indseed}
\Pr_p\big(\ispan(D_u)\big) \leq p^{\delta w(D_u)},
\end{equation}
where $\delta = \delta(\beta_1+\beta_2)$.

Next, we bound the probability of the event $\Delta(D_v,D_u)$ using Lemma~\ref{le:deltagrowth}. When $|N_{G_\hier}^\to(v)| \ne 1$ we do not have a lower bound on $w(D_u) - w(D_v)$, and in this case we will use the trivial bound $\Pr_p\big(\Delta(D_v,D_u)\big) \le 1$. Let us therefore assume that $N_{G_\hier}^\to(u) = \{v\}$ and $|N_{G_\hier}^\to(v)| = 1$, and hence that
$$\xi \cdot p^{-\beta_1(1-2\eta)-\eta} \le w(D_u) - w(D_v) \le 2\xi \cdot p^{-\beta_1(1-2\eta)-\eta},$$
by Definition~\ref{de:hier2}. Recalling that $D_u \subset D$ and $\beta_2 \le \alpha$, observe also that
$$p^{-\beta_1(1-2\eta)-\eta} \le w(D_u) \le p^{-(\beta_1+1)(1-2\eta)-\eta} \quad \text{and} \quad h(D_u) \le p^{-\alpha(1-2\eta)-\eta},$$
where the lower bound on $w(D_u)$ follows from Definition~\ref{de:hier2} (since $u$ is not a leaf), and the upper bounds follow from our assumptions on $w(D)$ and $h(D)$. It follows, by Lemma~\ref{le:deltagrowth}, that 
\begin{equation}\label{eq:inddelta}
\Pr_p\big(\Delta(D_v,D_u)\big) \le p^{\Omega(\delta')(w(D_u)-w(D_v))},
\end{equation}
where $\delta' = \delta'(\beta_1+\beta_2)$, and the constant implicit in $\Omega(\cdot)$ depends only on $\S_U$, for those pairs $\{u,v\} \subset V(G_\hier)$ with $N_{G_\hier}^\to(u) = \{v\}$ and $|N_{G_\hier}^\to(v)| = 1$.

We now divide into two cases according to the number of seeds of $\hier$: if there are `few' seeds then we will apply Lemma~\ref{lem:sumofwidths}; if there are `many' seeds, then the contribution to~\eqref{eq:indbound} from the seeds alone outweighs the combinatorial cost of counting the good hierarchies. To be precise, set $B := 
\xi^{-4} \cdot p^{-1+2\eta}$ and let
$$\hier^{(1)} := \big\{ \hier \in \hier_D \,:\, b(\hier) \leq B \big\} \quad \text{and} \quad \hier^{(2)} := \hier_D \setminus \hier^{(1)}.$$

Bounding the sum over $\hier \in \hier^{(2)}$ is the simpler case. Indeed, by Lemma~\ref{obs:all:seeds:big} we have $w(D_u) \ge \xi \cdot p^{-\beta_1(1-2\eta)-\eta} / 2$ for every $u \in L(\hier)$, and therefore, using the notation $\exp_p(x) := p^x$, and recalling that $\xi \ll \delta(\beta_1+\beta_2)$, we have
$$\sum_{\hier \in \hier^{(2)}} \prod_{u \in L(\hier)} \Pr_p\big( \ispan(D_u) \big) \le \sum_{b > B} |\hier_D^b| \cdot \exp_p\Big( \xi^2 \cdot b \cdot p^{-\beta_1(1-2\eta)-\eta} \Big),$$
by~\eqref{eq:indseed}. Moreover, by Lemma~\ref{le:counthier}, for each $b$ we have
$$\big|\hier_D^b\big|  \leq  \exp_p \Big( -O\Big( \xi^{-1} \cdot b \cdot w(D) \cdot p^{\beta_1(1-2\eta) + \eta} \Big) \Big) \le \exp_p \Big( - O\Big( \xi^{-1} \cdot b \cdot p^{- 1 + 2\eta} \Big) \Big),$$
since $w(D) \leq p^{-(\beta_1+1)(1 - 2\eta) - \eta}$. Hence, since $\beta_1 \geq 1$, 
\begin{equation}\label{eq:indhier2}
\sum_{\hier \in \hier^{(2)}} \prod_{u \in L(\hier)} \Pr_p\big( \ispan(D_u) \big) \le \sum_{b > B} \exp_p\Big( \xi^{3} \cdot b \cdot p^{-\beta_1(1-2\eta)-\eta} \Big) \le p^{w(D)},
\end{equation}
where the last step follows since $B = \xi^{-4} \cdot p^{-1+2\eta}$ and $w(D) \le p^{-(\beta_1+1)(1-2\eta)-\eta}$.

To deal with $\hier^{(1)}$, first we use the two estimates~\eqref{eq:indseed} and~\eqref{eq:inddelta} to obtain
\begin{multline}\label{eq:indhier1-1}
\sum_{\hier \in \hier^{(1)}} \bigg( \prod_{u \in L(\hier)} \Pr_p\big(\ispan(D_u)\big) \bigg)\bigg( \prod_{u \to v} \Pr_p\big(\Delta(D_v,D_u)\big) \bigg) \\
\leq \sum_{\hier\in\hier^{(1)}} \exp_p\bigg( \delta\sum_{u\in L(\hier)} w(D_u) \,+\, \Omega(\delta') \sum_{\substack{u \to v \\ d_\hier(v) = 1}} \big(w(D_u)-w(D_v) \big) \bigg),
\end{multline}
where $\delta = \delta(\beta_1+\beta_2)$ and $\delta' = \delta'(\beta_1+\beta_2)$. Now, recalling that $\beta_1 \ge 1$, and therefore $B = \xi^{-4} \cdot p^{-1+2\eta} \le \xi^2 \cdot p^{-\beta_1(1-2\eta) + \eta}$, it follows by Lemma~\ref{lem:sumofwidths} that
$$\sum_{u \in L(\hier)} w(D_u) + \sum_{\substack{u \to v \\ d_\hier(v) = 1}} \big( w(D_u) - w(D_v) \big) \ge \frac{w(D)}{10}.$$
Thus, recalling that $\delta(\beta_1+\beta_2) \ge \delta'(\beta_1+\beta_2)$, the right-hand side of~\eqref{eq:indhier1-1} is at most 
\begin{equation}\label{eq:indhier1-2}
\sum_{\hier\in\hier^{(1)}} p^{\Omega(\delta') w(D)} \le \sum_{b\leq B} \big|\hier_D^b\big| \cdot p^{\Omega(\delta') w(D)}.
\end{equation}
Now, by Lemma~\ref{le:counthier}, we have
$$\big|\hier_D^b\big|  \leq  \exp_p \Big( -O\Big( \xi^{-1} B \cdot w(D) \cdot p^{\beta_1(1-2\eta) + \eta} \Big) \Big) \le e^{w(D)},$$
for every $b \le B = \xi^{-4} \cdot p^{-1+2\eta}$, since $\beta_1 \ge 1$. The right-hand side of~\eqref{eq:indhier1-2} is therefore at most
$$B \cdot e^{w(D)} \cdot p^{\Omega(\delta') w(D)} \leq p^{\Omega(\delta') w(D)},$$
since $B = \xi^{-4} \cdot p^{-1+2\eta} \le e^{w(D)}$. Combining this bound with~\eqref{eq:indbound} and~\eqref{eq:indhier2}, it follows that $\Pr_p\big(\ispan(D)\big) \le p^{\Omega(\delta') w(D)}$, as required.
\end{proof}

\subsection{The induction steps}\label{se:induction}

In this subsection we prove Lemma~\ref{le:induction}, which states that $\ih(\alpha+1,\alpha)$ and $\ih(\alpha,\alpha+1)$ both hold.

\begin{proof}[Proof of Lemma~\ref{le:induction}]
We shall prove by induction on $\beta_1 + \beta_2$ that $\ih(\beta_1,\beta_2)$ holds for every pair $(\beta_1,\beta_2) \in \N^2$ with 
\[
2 \leq \beta_1 + \beta_2 \leq 2\alpha + 1 \qquad \text{and} \qquad |\beta_1 - \beta_2| \leq 1.
\]
Observe first that $\ih(1,1)$ follows from Lemma~\ref{le:smalldroplet}, since $\delta(2)$ was chosen (in~\eqref{eq:constants1}) to be sufficiently small (depending on $\eta$). The induction steps are of three different types, which are dealt with in the following three claims. The first follows from Lemma~\ref{le:indhoriz}, the second from Lemma~\ref{le:indvert}, and the third requires both.  

\begin{claim}\label{cl:indhoriz}
For each $1\leq\beta\leq\alpha$ we have
\[
\ih(\beta,\beta) \, \Rightarrow \, \ih(\beta+1,\beta).
\]
\end{claim}

\begin{proof}\let\qednow\qedsymbol\renewcommand{\qedsymbol}{}
Let $D$ be a droplet with
\[
w(D) \leq p^{-(\beta+1)(1-2\eta) - \eta} \qquad \text{and} \qquad h(D) \leq p^{-\beta(1-2\eta) - \eta}.
\]
We are required to show that $\Pr_p\big(\ispan(D)\big)\leq p^{\delta \max\{w(D), h(D)\}}$, where $\delta = \delta(2\beta+1)$. If $w(D) \leq p^{-\beta(1-2\eta) - \eta}$ then this follows immediately from $\ih(\beta,\beta)$ (since we chose $\delta(2\beta+1) \leq \delta(2\beta)$ in~\eqref{eq:constants1}), so we may assume that $w(D) \geq p^{-\beta(1-2\eta) - \eta}$. Now, applying Lemma~\ref{le:indhoriz} with $\beta_1=\beta_2=\beta$, it follows that
\[
\Pr_p\big(\ispan(D)\big) \leq p^{\Omega(\delta'(2\beta)) w(D)} \leq p^{\delta(2\beta+1) w(D)},
\]
as required, since we chose $\delta(2\beta+1) \ll \delta'(2\beta)$ in~\eqref{eq:constants2}. \qednow
\end{proof}

\begin{claim}\label{cl:indvert}
For each $1\leq\beta\leq\alpha$ we have
\[
\ih(\beta,\beta) \, \Rightarrow \, \ih(\beta,\beta+1).
\]
\end{claim}

\begin{proof}\let\qednow\qedsymbol\renewcommand{\qedsymbol}{}
Let $D$ be a droplet with
\[
w(D)\leq p^{-\beta(1-2\eta) - \eta} \qquad \text{and} \qquad h(D) \leq p^{-(\beta+1)(1-2\eta) - \eta}.
\]
We again need to show that $\Pr_p\big(\ispan(D)\big)\leq p^{\delta \max\{w(D), h(D)\}}$, where $\delta = \delta(2\beta+1)$. Note that if $h(D) \leq p^{-\beta(1-2\eta) - \eta}$ then this follows immediately from $\ih(\beta,\beta)$, as before, so we may assume that $h(D) \geq p^{-\beta(1-2\eta) - \eta}$, which implies that
$$h(D) \ge \xi^{-3} \max\big\{ p^{1-\eta} w(D), \, 1 \big\}.$$
Hence, applying Lemma~\ref{le:indvert} with $\beta_1 = \beta_2 = \beta$, we obtain
$$\Pr_p\big(\ispan(D)\big) \le p^{\delta'(2\beta) h(D) / 8} \leq p^{\delta(2\beta+1) \max\{w(D), h(D)\}},$$
as required, since $h(D) \ge w(D)$ and we chose $\delta(2\beta+1) \ll \delta'(2\beta)$ in~\eqref{eq:constants2}. \qednow
\end{proof}

\begin{claim}\label{cl:indboth}
For each $1\leq\beta\leq\alpha-1$ we have
\[
\big( \, \ih(\beta+1,\beta) \, \wedge \, \ih(\beta,\beta+1) \, \big) \, \Rightarrow \, \ih(\beta+1,\beta+1).
\]
\end{claim}

\begin{proof}\let\qednow\qedsymbol\renewcommand{\qedsymbol}{}
Let $D$ be an $\stabu$-droplet with
\[
w(D) \leq p^{-(\beta+1)(1-2\eta) - \eta} \qquad \text{and} \qquad h(D) \leq p^{-(\beta+1)(1-2\eta) - \eta}.
\]
This time we are required to show that $\Pr_p\big(\ispan(D)\big)\leq p^{\delta(2\beta+2) \max\{w(D), h(D)\}}$. Note that, since we chose $\delta(2\beta+2) \leq \delta(2\beta+1)$ in~\eqref{eq:constants1}, this follows immediately from $\ih(\beta,\beta+1)$ if $w(D) \leq p^{-\beta(1-2\eta) - \eta}$, and from $\ih(\beta+1,\beta)$ if $h(D) \leq p^{-\beta(1-2\eta) - \eta}$. We may therefore assume that 
$$\min\big\{ w(D), h(D) \big\} \geq p^{-\beta(1-2\eta) - \eta}.$$

Suppose first that $w(D) \geq h(D)$. Then, applying Lemma~\ref{le:indhoriz} with $\beta_1= \beta$ and $\beta_2=\beta+1$, it follows that
\[
\Pr_p\big(\ispan(D)\big)  \leq  p^{\Omega(\delta'(2\beta+1)) w(D)} \leq  p^{\delta(2\beta+2) w(D)},
\]
as required, since we chose $\delta(2\beta+2) \ll \delta'(2\beta+1)$ in~\eqref{eq:constants2}.

\pagebreak

On the other hand, if $w(D) \leq h(D)$ then we have 
$$h(D) \ge \xi^{-3} \max\big\{ p^{1-\eta} w(D), \, 1 \big\}.$$
Hence, applying Lemma~\ref{le:indvert} with $\beta_1 = \beta + 1$ and $\beta_2 = \beta$, we obtain
$$\Pr_p\big(\ispan(D)\big) \le p^{\delta'(2\beta) h(D) / 8} \le p^{\delta(2\beta+2) h(D)},$$
as required, since we chose $\delta(2\beta+2) \ll \delta'(2\beta+1)$ in~\eqref{eq:constants2}. \qednow
\end{proof}

Together with $\ih(1,1)$, these claims imply $\ih(\alpha+1,\alpha)$ and $\ih(\alpha,\alpha+1)$, which completes the proof of the lemma.
\end{proof}

\subsection{Internally spanned critical droplets}\label{se:criticaldroplets}

Recall from Definition~\ref{de:critical} that we call an $\S_U$-droplet $D$ \emph{critical} if one of the following holds:\smallskip
\begin{itemize}
\item[$(T)$] $w(D)\leq 3 p^{-\alpha-1/5}$ and $\frac{\xi}{p^\alpha}\log\frac{1}{p}\leq h(D)\leq\frac{3\xi}{p^\alpha}\log\frac{1}{p}$, or \smallskip
\item[$(L)$] $p^{-\alpha-1/5} \leq w(D) \leq 3p^{-\alpha-1/5}$ and $h(D) \leq \frac{\xi}{p^\alpha}\log\frac{1}{p}$,\smallskip
\end{itemize}
where $\xi > 0$ is the (sufficiently small) constant chosen in~\eqref{eq:constants1}. In this subsection we will prove the following bound on the probability that a critical droplet is internally spanned, which we shall show in Section~\ref{se:lowerunbalend} easily implies Theorem~\ref{th:lowerunbal}. Recall from~\eqref{eq:constants3} that $\delta > 0$ is chosen so that $\delta \ll \xi \ll \delta'(2\alpha+1)$.

\begin{lemma}\label{le:critunbal}
If $D$ is a critical droplet then
\[
\Pr_p\big(\ispan(D)\big) \leq \exp\left( -\frac{\delta}{p^\alpha}\left(\log\frac{1}{p}\right)^2 \right).
\]
\end{lemma}

To prove the lemma for `tall' droplets (type $(T)$), we simply apply Lemmas~\ref{le:induction} and~\ref{le:indvert}. For `long' droplets (type $(L)$), on the other hand, we need to apply the method of hierarchies, as in the proof of Lemma~\ref{le:indhoriz}, together with Lemmas~\ref{le:induction} and~\ref{le:deltagrowth}, which we use to bound the probabilities of the events in Definition~\ref{def:satisfied}.

\begin{proof}[Proof of Lemma~\ref{le:critunbal}]
We will prove that the lemma holds with
\[
\delta = \frac{\xi \cdot \delta'(2\alpha+1)}{8}.
\]
Let $D$ be a critical droplet, and suppose first that $D$ is of type $(T)$. Then, since $w(D) \leq 3 p^{-\alpha-1/5}$ and $h(D) \geq \frac{\xi}{p^\alpha}\log\frac{1}{p}$, and recalling that $\eta = (10\alpha)^{-1}$, we have 
\[
w(D) \le p^{- (\alpha + 1)(1 - 2\eta) - \eta} \quad \text{and} \quad h(D)  \ge \xi^{-3} \max\big\{ p^{1-\eta} w(D), \, 1 \big\}.
\]
We may therefore apply Lemma~\ref{le:indvert} with $\beta_1=\alpha+1$ and $\beta_2=\alpha$, since $\ih(\alpha+1,\alpha)$ holds by Lemma~\ref{le:induction}. Since an internally spanned droplet is vertically crossed, this gives
\[
\Pr_p\big(\ispan(D)\big) \leq  p^{\delta'(2\alpha+1) h(D) / 8}  \le \exp\left( -\frac{\delta}{p^\alpha}\left(\log\frac{1}{p}\right)^2 \right),
\]
as required.

So suppose from now on that $D$ is of type $(L)$; in this case we will prove the following much stronger bound:
\begin{equation}\label{eq:critical:stronger:L}
\Pr_p\big(\ispan(D)\big) \leq \exp\Big( - p^{- \alpha - 1/6} \Big).
\end{equation}
We apply the hierarchies framework, as in the proof of Lemma~\ref{le:indhoriz}, but with $\beta = \alpha$. By Lemma~\ref{le:boundoverH}, we have
\begin{equation}\label{eq:critbound}
\Pr_p\big(\ispan(D)\big) \leq \sum_{\hier \in \hier_D} \bigg( \prod_{u \in L(\hier)} \Pr_p\big(\ispan(D_u)\big) \bigg)\bigg( \prod_{u \to v} \Pr_p\big(\Delta(D_v,D_u)\big) \bigg).
\end{equation}
Now, if $u\in L(\hier)$, then since $D$ is of type $(L)$, and by Definition~\ref{de:hier2}, 
$$w(D_u) \leq p^{-\alpha(1 - 2\eta) - \eta} \quad \text{and} \quad h(D_u) \leq h(D) \leq \frac{\xi}{p^\alpha}\log\frac{1}{p} \leq p^{-(\alpha+1)(1 - 2\eta) - \eta}.$$
Since $\ih(\alpha,\alpha+1)$ holds, by Lemma~\ref{le:induction}, it follows that
\begin{equation}\label{eq:critseed}
\Pr_p\big(\ispan(D_u)\big) \leq  p^{\delta(2\alpha+1) w(D_u)},
\end{equation}
for every $u\in L(\hier)$.

Next, if $N_{G_\hier}^\to(u) = \{v\}$ and $|N_{G_\hier}^\to(v)| \ne 1$ then (as in the proof of Lemma~\ref{le:indhoriz}) we use the trivial bound $\Pr_p\big(\Delta(D_v,D_u)\big) \le 1$. On the other hand, if $N_{G_\hier}^\to(u) = \{v\}$ and $|N_{G_\hier}^\to(v)| = 1$ then we will bound the probability of the event $\Delta(D_v,D_u)$ using Lemma~\ref{le:deltagrowth}. To do so, recall that 
$$\xi \cdot p^{-\alpha(1-2\eta)-\eta} \le w(D_u) - w(D_v) \le 2\xi \cdot p^{-\alpha(1-2\eta)-\eta},$$
by Definition~\ref{de:hier2}. Moreover, since $D$ is a type $(L)$ critical droplet, we also have
$$p^{-\alpha(1-2\eta)-\eta} \le w(D_u) \le p^{-(\alpha+1)(1-2\eta)-\eta} \quad \text{and} \quad h(D_u) \leq h(D) \leq \frac{\xi}{p^\alpha}\log\frac{1}{p},$$
where the lower bound on $w(D_u)$ follows from Definition~\ref{de:hier2} (since $u$ is not a leaf). Since $\ih(\alpha,\alpha+1)$ holds, it follows by Lemma~\ref{le:deltagrowth} that 
\begin{equation}\label{eq:critdelta}
\Pr_p\big(\Delta(D_v,D_u)\big) \le \exp\Big(-p^{O(\xi)}\big(w(D')-w(D)\big)\Big),
\end{equation}
where the implicit constant depends on $\weak(2\alpha+1)$, for those pairs $\{u,v\} \subset V(G_\hier)$ with $N_{G_\hier}^\to(u) = \{v\}$ and $|N_{G_\hier}^\to(v)| = 1$.

As in the proof of Lemma~\ref{le:indhoriz}, we divide into two cases according to whether $\hier$ has many or few seeds. To be precise, let $B := p^{-2/3}$ and define
$$\hier^{(1)} := \big\{ \hier \in \hier_D \,:\, b(\hier) \leq B \big\} \quad \text{and} \quad \hier^{(2)} := \hier_D \setminus \hier^{(1)}.$$
Bounding the sum in~\eqref{eq:critbound} over $\hier \in \hier^{(2)}$ is again the easier case. Indeed, by Lemma~\ref{obs:all:seeds:big} we have $w(D_u) \ge \xi \cdot p^{-\alpha(1-2\eta)-\eta} / 2$ for every $u \in L(\hier)$, and therefore, by~\eqref{eq:critseed}, and recalling that $\xi \ll \delta(2\alpha + 1)$, we have
$$\sum_{\hier \in \hier^{(2)}} \prod_{u \in L(\hier)} \Pr_p\big( \ispan(D_u) \big) \le \sum_{b > B} |\hier_D^b| \cdot \exp_p\Big( \xi^2 \cdot b \cdot p^{-\alpha(1-2\eta)-\eta} \Big),$$
where we again use the notation $\exp_p(x) := p^x$. Now, by Lemma~\ref{le:counthier}, we have
$$\big|\hier_D^b\big| \le \exp_p \Big( -O\Big( \xi^{-1} \cdot b \cdot w(D) \cdot p^{\alpha(1-2\eta) + \eta} \Big) \Big) \le \exp_p \Big( - O\Big( \xi^{-1} \cdot b \cdot p^{- 1 + 2\eta} \Big) \Big)$$
for each $b \in \N$, since $w(D) \le p^{-(\alpha+1)(1 - 2\eta) - \eta}$. Since $\alpha \ge 1$, it follows that
\begin{equation}\label{eq:crithier2}
\sum_{\hier \in \hier^{(2)}} \prod_{u \in L(\hier)} \Pr_p\big( \ispan(D_u) \big) \le \sum_{b > B} \exp_p\Big( \xi^{3} \cdot b \cdot p^{-\alpha(1-2\eta)-\eta} \Big) \le \exp\Big( -p^{-\alpha-1/3} \Big),
\end{equation}
where the last step follows since $\eta = (10\alpha)^{-1}$ and $B = p^{-2/3}$. 

For hierarchies with few seeds, observe that by~\eqref{eq:critseed} and~\eqref{eq:critdelta}, we have
\begin{multline}\label{eq:crithier1-1}
\sum_{\hier \in \hier^{(1)}} \prod_{u \in L(\hier)} \Pr_p\big( \ispan(D_u) \big) \prod_{u \to v} \Pr_p\big(\Delta(D_v,D_u)\big) \\
\le \sum_{\hier \in \hier^{(1)}} \exp\bigg[ - p^{O(\xi)} \bigg( \sum_{u \in L(\hier)} w(D_u) + \sum_{\substack{u \to v \\ d_\hier(v) = 1}} \big( w(D_u) - w(D_v) \big) \bigg) \bigg]
\end{multline}
where we used the bound $\delta(2\alpha+1) \log(1/p) > p^{O(\xi)}$. By Lemma~\ref{lem:sumofwidths}, we have
$$\sum_{u \in L(\hier)} w(D_u) + \sum_{\substack{u \to v \\ d_\hier(v) = 1}} \big( w(D_u) - w(D_v) \big) \ge \frac{w(D)}{10}.$$
Thus, the right-hand side of~\eqref{eq:crithier1-1} is at most
\begin{equation}\label{eq:crithier1-2}
\sum_{b \le B} \big|\hier_D^b\big| \cdot \exp\big( - p^{O(\xi)} w(D) \big).
\end{equation}
Now, by Lemma~\ref{le:counthier}, we have
$$\big|\hier_D^b\big| \le \exp_p \Big( -O\big( \xi^{-1} B \cdot w(D) \cdot p^{\alpha(1-2\eta) + \eta} \big) \Big) \le \exp\big( p^{1/5} w(D) \big),$$
for every $b \le B = p^{-2/3}$, since $\alpha \ge 1$ and $\eta = (10\alpha)^{-1}$. Hence, combining this with~\eqref{eq:critbound} and the bound~\eqref{eq:crithier2} proved above for $\hier\in\hier^{(2)}$, it follows that  
\begin{align*}
\Pr_p\big(\ispan(D)\big) & \leq \, \sum_{b = 1}^B \big|\hier_D^b\big| \exp\Big( - p^{O(\xi)} w(D) \Big) + \exp\Big( -p^{-\alpha-1/3} \Big),\\
& \leq B \exp\Big( - p^{- \alpha - 1/5 + O(\xi)} \Big) \leq \exp\Big( - p^{- \alpha - 1/6} \Big),
\end{align*}
since $w(D) \ge p^{- \alpha - 1/5}$ and $\xi > 0$ was chosen to be sufficiently small. This proves~\eqref{eq:critical:stronger:L}, and hence completes the proof of the lemma.
\end{proof}

\subsection{The proof of Theorem~\ref{th:lowerunbal}}\label{se:lowerunbalend}

We need one final lemma in order to deduce Theorem~\ref{th:lowerunbal} from Lemma~\ref{le:critunbal}. It is a simple consequence of Lemma~\ref{le:ALspan}. 

\pagebreak

\begin{lemma}\label{lem:exists:critdroplet}
If $n \ge p^{-3\alpha}$ and $[A] = \Z_n^2$, then there exists a critical droplet that is internally spanned by $A$.
\end{lemma}

\begin{proof}
Run the spanning algorithm on $\Z_n^2$, with initial set $A$. Since $[A] = \Z_n^2$, we will (at some point in the algorithm) obtain an internally spanned droplet $D_0$ with $\max\{ w(D_0), h(D_0) \} \ge p^{-2\alpha}$. Let $D_0$ be the first such droplet to appear in the algorithm, and note that at this scale the wrap-around effect of the torus is irrelevant, since we will only use the fact that $D_0$ is internally spanned. We may therefore work from now on in $\Z^2$, cf.~Remark~\ref{remark:wraparound}.

Suppose first that $w(D_0) \leq 3 p^{-\alpha-1/5}$, so $h(D_0) \geq p^{-2\alpha}$. Since $D_0$ is internally spanned, by applying Lemma~\ref{le:ALspan} with $u = u^*$ and $k = \frac{\xi}{p^\alpha}\log\frac{1}{p}$, we obtain an internally spanned droplet $D \subset D_0$ with 
\begin{equation}\label{eq:crit:height}
\frac{\xi}{p^\alpha} \log \frac{1}{p} \leq h(D) \leq \frac{3\xi}{p^\alpha} \log \frac{1}{p},
\end{equation}
so $D$ is a type $(T)$ critical droplet. 

On the other hand, if $w(D_0) \geq 3 p^{-\alpha-1/5}$ then we can apply Lemma~\ref{le:ALspan} with $u = w^\perp$ and $k = p^{-\alpha-1/5}$ to obtain an internally spanned droplet $D_1 \subset D_0$ with 
\begin{equation}\label{eq:crit:width}
p^{-\alpha-1/5} \leq w(D_1) \leq 3p^{-\alpha-1/5}.
\end{equation}
If $h(D_1) \leq \frac{\xi}{p^\alpha}\log\frac{1}{p}$ then $D_1$ is a type $(L)$ critical droplet, in which case we are done, so assume not. Now, applying Lemma~\ref{le:ALspan} with $u = u^*$ and $k = \frac{\xi}{p^\alpha}\log\frac{1}{p}$, we obtain an internally spanned droplet $D \subset D_1$ such that~\eqref{eq:crit:height} holds. Since $w(D) \le w(D_1) \leq 3p^{-\alpha-1/5}$, by~\eqref{eq:crit:width}, it follows that $D$ is a type $(T)$ critical droplet, as required.
\end{proof}

We now have all the tools we need to complete the proof of Theorem~\ref{th:lowerunbal}, and hence Theorem~\ref{th:main}. 

\begin{proof}[Proof of Theorem~\ref{th:lowerunbal}]
Let $\eps > 0$ be the sufficiently small constant chosen in~\eqref{eq:constants1}, so in particular $\eps \ll \delta$. Set
\[
p = \bigg( \frac{\eps(\log\log n)^2}{\log n} \bigg)^{1/\alpha},
\]
and let $A$ be a $p$-random subset of $\Z_n^2$. We claim that $\Pr_p\big( [A] = \Z_n^2 \big) \to 0$ as $n \to \infty$. 

Indeed, if $[A] = \Z_n^2$ then, by Lemma~\ref{lem:exists:critdroplet}, there exists a critical droplet $D$ that is internally spanned by $A$. By Lemma~\ref{le:critunbal}, the probability that $D$ is internally spanned is at most
\[
\exp\left( -\frac{\delta}{p^\alpha}\left(\log\frac{1}{p}\right)^2 \right),
\]
and there are at most $n^2 p^{-4\alpha} \leq n^3$ critical droplets in $\Z_n^2$. Hence
\[
\Pr_p\big( [A] = \Z_n^2 \big) \leq n^3 \exp\left( -\frac{\delta}{p^\alpha}\left(\log\frac{1}{p}\right)^2 \right)  \to  0
\]
as $n \to \infty$, as required, since $\eps \ll \delta$. This completes the proof of the theorem.
\end{proof}

\section{Conjectures for higher dimensions}\label{se:highd}

We conclude by briefly discussing the $\U$-bootstrap percolation models in higher dimensions. Fix an integer $d\geq 2$ and let $\U$ be a $d$-dimensional update family. The definition of the stable set $\stab = \stab(\U)$ is the natural generalization of the two-dimensional definition:
\[
\stab := \big\{ u\in S^{d-1} \,:\, [\H_u^d]=\H_u^d \big\},
\]
where
\[
\H_u^d := \big\{ x\in\Z^d \,:\, \<x,u\> < 0 \big\}
\]
is the discrete half-space in $\Z^d$ with normal $u\in S^{d-1}$. Observe that, as in two dimensions, it is easy to show that the dichotomy $[\H_u^d]\in\big\{\H_u^d,\Z^d\big\}$ holds for any unit vector $u\in S^{d-1}$.

Given a set $\T \subset S^{d-1}$, we write $\interior(\T)$ for the interior of $\T$ in the usual topology on $S^{d-1}$ (induced by geodesic distance). Generalizing\footnote{Recall from Lemma~\ref{le:stableset} that in two dimensions $\interior(C \cap\stab) = \emptyset$ if and only if $C \cap\stab$ is finite.} Definition~\ref{de:class}, we classify $d$-dimensional update familes as follows.

\begin{definition}
A $d$-dimensional update family is:
\begin{itemize}
\item \emph{subcritical} if $\interior(C\cap\stab) \neq \emptyset$ for every hemisphere $C\subset S^{d-1}$;\vspace{0.1cm}
\item \emph{critical} if there exists a hemisphere $C\subset S^{d-1}$ such that $\interior(C\cap\stab) = \emptyset$ and if $C\cap\stab\neq\emptyset$ for every open hemisphere $C\subset S^{d-1}$;\vspace{0.1cm}
\item \emph{supercritical} if $C\cap\stab = \emptyset$ for some open hemisphere $C\subset\stab$.
\end{itemize}
\end{definition}

As in two dimensions, the subcritical/critical/supercritical trichotomy depends only on the stable set $\S$. We make the following `universality' conjecture, which says that in $d$ dimensions, every update family  behaves (roughly speaking) like the classical $r$-neighbour model for some $1 \le r \le d + 1$. 

\begin{conjecture}\label{con:highd}
Let\/ $\U$ be a $d$-dimensional bootstrap percolation update family.
\begin{itemize}
\item[$(i)$] If\/ $\U$ is supercritical then $p_c(\Z_n^d,\U) = n^{-\Theta(1)}$.\vspace{0.1cm}
\item[$(ii)$] If\/ $\U$ is critical then there exist $r\in\{2,\dots,d\}$ 
such that\/\footnote{Here $\log_{(r)}$ denotes an $r$-times iterated logarithm, so $\log_{(1)} = \log$ and $\log_{(r)} = \log \log_{(r-1)}$.}
\[
p_c(\Z_n^d,\U) = \bigg( \frac{1}{\log_{(r-1)}n}\bigg)^{\Theta(1)}. 
\]
\item[$(iii)$] If\/ $\U$ is subcritical then $p_c(\Z^d,\U) > 0$.
\end{itemize}
\end{conjecture}

We remark that the conjecture for subcritical families was originally made by Balister, Bollob\'as, Przykucki and Smith~\cite{BBPS}. We expect Conjecture~\ref{con:highd} to be extremely difficult, and a proof appears to be far out of reach with current techniques.\footnote{While this sentence was true in 2014, when it was first written, a proof of Conjecture~\ref{con:highd} was recently announced by Balister and three of the current authors~\cite{BBMSlower,BBMSupper,BBMSsub}.}

\addtocontents{toc}{\protect\setcounter{tocdepth}{0}}
\section*{Acknowledgements}
\addtocontents{toc}{\protect\setcounter{tocdepth}{1}}

This work was started during a visit of the second author to IMPA in the spring of 2013. He would like to thank IMPA for its hospitality. The authors would also like to express their sincere gratitude to the anonymous referee, who (twice) read this long paper extremely carefully, pointed out several small gaps in the original proof, and made a large number of additional helpful comments on the paper. 

\enlargethispage*{\baselineskip}

\bibliographystyle{amsplain}
\bibliography{critbprefs}

\end{document}